\documentclass[reqno,11pt]{amsart}
\usepackage{amssymb,mathrsfs}
\usepackage[usenames,dvipsnames]{color}
\usepackage[a4paper,margin=1in]{geometry}
\usepackage{microtype}
\usepackage{enumitem}
\usepackage[colorlinks=true,linkcolor=Red,citecolor=Green,pdfstartview=XYZ]{hyperref}
\usepackage[open,openlevel=2,depth=3,atend]{bookmark}

\newtheorem{theorem}{Theorem}[section]
\newtheorem{lemma}[theorem]{Lemma}
\newtheorem{proposition}[theorem]{Proposition}

\theoremstyle{remark}
\newtheorem{remark}[theorem]{Remark}

\numberwithin{equation}{section}

\renewcommand{\Im}{\operatorname{Im}}
\renewcommand{\Re}{\operatorname{Re}}
\DeclareMathOperator*{\supp}{supp}
\DeclareMathOperator{\rank}{rank}
\DeclareMathOperator{\Tr}{Tr}
\DeclareMathOperator{\ran}{ran}
\DeclareMathOperator{\tr}{tr}
\DeclareMathOperator{\Id}{Id}

\newcommand{\R}{{\mathbb R}}
\newcommand{\eps}{\varepsilon}
\newcommand{\PsiDO}{\Psi_{\mathrm{cl}}}
\newcommand{\cc}{\mathrm{cc}}
\newcommand{\cH}{\mathcal H}
\newcommand{\cD}{\mathcal D}
\newcommand{\cN}{\mathcal N}
\newcommand{\tang}{\mathbf t}
\newcommand{\normal}{\mathbf n}
\newcommand{\restr}{\big|}

\title{A Geometric Calder\'on Problem for Differential Forms}
\author{Yujun Jin}
\address{School of Mathematical Sciences\\
	University of Science and Technology of China\\
	Hefei, 230026\\ P.R. China}
\email{byjyj@mail.ustc.edu.cn}

\author{Zuoqin Wang}
\address{School of Mathematical Sciences\\
	University of Science and Technology of China\\
	Hefei, 230026\\ P.R. China}
\email{wangzuoq@ustc.edu.cn}

\date{\today}

\begin{document}

\begin{abstract}
We establish boundary determination and global partial-data uniqueness for
Karpukhin's Dirichlet-to-Neumann operator on differential $p$-forms.
We derive an explicit recursion for its full symbol and determine the full
boundary jet of the metric, modulo conformal gauge in the critical
dimension $n=2p+2$. For $p\geq1$, partial data determine metrics that are analytic
in the interior and near the measured boundary subset, up to isometry or
conformal equivalence. The global uniqueness relies on a new
projection-kernel embedding that detects the unobserved boundary through
rank loss and small eigenvalues.  We also establish global partial-data uniqueness for smooth metrics with harmonic curvature, using elliptic–transport unique continuation and, in the critical dimension, conformal-geodesic coordinates. A Weyl nondegeneracy condition further upgrades conformal equivalence to isometry.
\end{abstract}
\maketitle

\section{Introduction and main results}
Let $(M^n,g)$ be a compact oriented smooth Riemannian manifold with
nonempty smooth boundary $Y:=\partial M$, outward unit normal $\nu_g$,
and induced boundary metric $h$. The scalar Dirichlet-to-Neumann
operator is defined by
\begin{equation*}
 \mathcal L_g \varphi:=\partial_{\nu_g} u_\varphi\big|_Y,\qquad \Delta_g u_\varphi=0,\quad u_\varphi\big|_Y=\varphi.
\end{equation*}
Its spectrum and
its relation to boundary geometry are discussed in the surveys
\cite{GirouardPolterovich2017,ColboisGirouardGordonSher2024}.
The geometric Calder\'on problem asks whether the interior geometry
can be recovered from this boundary operator
\cite{LassasLiimatainenSalo2020}. We study the corresponding partial-data inverse problem for Karpukhin's Dirichlet-to-Neumann operator on differential forms. Our aim is to determine how much of the interior geometry, including the topology of the manifold, is encoded by measurements on only a nonempty open subset of the boundary. The results combine boundary determination, global reconstruction, and a smooth uniqueness theorem for metrics with harmonic curvature.

There are several nonequivalent Dirichlet-to-Neumann
constructions associated with the Hodge Laplacian \cite{JoshiLionheart,RaulotSavo2012}. In the present paper, we focus on Karpukhin's Hodge-star modification \cite{Karpukhin} of the Dirichlet-to-Neumann operator considered by Belishev--Sharafutdinov \cite{BelishevSharafutdinov}.

The exterior derivative, codifferential, and Hodge Laplacian on $M$ are denoted
by $d$, $\delta$, and $ \Delta=d\delta+\delta d$. Let $i:Y\hookrightarrow M$ denote the
inclusion.  For a differential form $\eta$ on $M$ we
write
\begin{equation*}
 \tang\eta:=i^*\eta,\qquad\normal\eta:=i^*(\iota_{\nu_g}\eta),
\end{equation*}
for its tangential and normal traces.  Thus $\tang\eta$ has the same degree as
$\eta$, whereas $\normal\eta$ has degree one less.

Throughout the paper, we assume $0\leq p\leq n-2$. Let $\Omega^p(Y):=C^\infty(Y;\Lambda^pT^*Y)$ denote the space of smooth $p$-forms on $Y$.
For $\varphi\in\Omega^p(Y)$,  consider the boundary value problem
\begin{equation}\label{eoriginalproblem}
 \Delta\omega=0,\qquad \delta\omega=0,\qquad \tang\omega=\varphi.
\end{equation}
A solution exists and is unique modulo the finite-dimensional space of
Dirichlet harmonic fields
\[
 \cH_D^p(M) :=\{\eta\in\Omega^p(M):d\eta=0,\ \delta\eta=0,\ \tang\eta=0\}
\]
(see \cite[Chapter~3]{Schwarz} and \cite[Proposition~3.11]{Karpukhin}).  Since every
$\eta\in\cH_D^p(M)$ is closed, the boundary form $\normal(d\omega)$ is
independent of the choice of solution.  We therefore consider Karpukhin's Dirichlet-to-Neumann operator $\Lambda: \Omega^p(Y) \to  \Omega^p(Y)$ defined by 
\begin{equation*}
 \Lambda_g\varphi:=\normal(d\omega).
\end{equation*} 
For $p=0$, $\Lambda_g$ is the scalar Dirichlet-to-Neumann operator $\mathcal L_g$ which is elliptic. For $p\geq1$, however, $\Lambda_g$ is non-elliptic: it annihilates exact forms and its range consists of coclosed forms \cite[Theorem~2.3]{Karpukhin}. Consequently,
\begin{equation}\label{Toeidentity}
 \Lambda_g=\Pi_{\cc}\Lambda_g\Pi_{\cc},
\end{equation}
where $\Pi_{\cc}$ denotes the $L^2(Y)$-orthogonal projection onto
$\ker\delta_Y$.
For $\xi\in T_y^*Y\setminus\{0\}$, let $\eps_\xi$ denote exterior multiplication by $\xi$, i.e. $\eps_\xi a:=\xi\wedge a$, and let $\iota_{\xi^\sharp}$ denote contraction with its metric dual $\xi^\sharp\in T_yY$. By \cite[Eq.~(29)]{JakobsonStrohmaier2007}, we have
$\Pi_{\cc}\in\PsiDO^0(Y;\Lambda^pT^*Y)$ with principal symbol
\begin{equation}\label{rankpcc}
 \pi_{\cc}(y,\xi)=\Id-|\xi|_h^{-2}\eps_\xi\iota_{\xi^\sharp},\qquad\rank\pi_{\cc}(y,\xi)=\binom{n-2}{p}.
\end{equation}
Very recently, Chen--Huang--Liao \cite{CHL26} proved that $\Lambda_g$ is a classical pseudodifferential operator with principal symbol 
\[\sigma_1(\Lambda_g)(y,\xi)=|\xi|_h\pi_\cc=|\xi|_h^{-1}\iota_{\xi^\sharp}\eps_\xi.\] 
For related inverse problems involving $\Lambda_g$, see
\cite{Shonkwiler2013,SharafutdinovShonkwiler2013,SharafutdinovLinearized}.

\subsection{Boundary determination and the critical conformal gauge}
The full symbol is the first link between the boundary measurement and
the metric. In dimensions $n\geq3$, Joshi--Lionheart
\cite[Corollary~2.1(ii)]{JoshiLionheart} showed that the full symbol of the complete Dirichlet-to-Neumann map for the Hodge Laplacian $\Delta$ determines the boundary metric jet (this map corresponds to our $\mathcal D_p$ in \eqref{def::dpHodge} up to boundary
Hodge stars and sign conventions).
For the scalar case $p=0$, Lee--Uhlmann
\cite[Proposition~1.3]{LeeUhlmann} established boundary metric
determination to all orders from the full symbol $\sigma(\mathcal L_g)$ in dimensions
$n\geq3$.
We derive a recursive formula for the full symbol $\sigma(\Lambda_g)$
and use it to determine the boundary metric jet to all orders,
uniquely when $n\neq 2p+2$ and modulo conformal gauge in the
critical case $n=2p+2$.

For a smooth tensor $T$, let $j_Y^\infty T$ denote its boundary normal
jet. In fixed boundary normal coordinates $(r,y)$, this means
\begin{equation*}
 j_Y^\infty T_1=j_Y^\infty T_2\quad\Longleftrightarrow\quad\left.\partial_r^k(T_1-T_2)\right|_{r=0}=0\qquad\text{for all }k\geq0.
\end{equation*}
For an open boundary set $\Gamma\subset Y$, let $j_\Gamma^kT$
denote the normal derivatives $\partial_r^qT|_\Gamma$ for
$0\leq q\leq k$, and let $j_\Gamma^\infty T$ denote all normal
derivatives along $\Gamma$.  Taking $k=0$ gives equality of the tensors themselves on $\Gamma$.
Tangential differentiation of the identities for the normal
derivatives then yields equality of all mixed derivatives on $\Gamma$.
Let $\left.\Lambda_{g}\right|_{\Gamma}$ denote the map
with inputs in $\Omega_c^p(\Gamma)$ and extended by zero to 
$Y$. Our first theorem establishes boundary determination from measurements
restricted to an arbitrary nonempty open subset of the boundary.
\begin{theorem}
\label{thboundaryintrodo}
Let $(M_j,g_j)$, $j=1,2$, be compact oriented smooth Riemannian
$n$-manifolds, let $0\leq p\leq n-2$, and let
$\Gamma_j\subset\partial M_j$ be nonempty domains identified
by a smooth diffeomorphism $\kappa:\Gamma_1\to\Gamma_2$.
Assume
\begin{equation}\label{equsntrosmoothingdiff}
   \left.\Lambda_{g_1}\right|_{\Gamma_1} -\kappa^*\left.\Lambda_{g_2}\right|_{\Gamma_2}
       (\kappa^{-1})^*\,\in\PsiDO^{-\infty}(\Gamma_1;\Lambda^pT^*\Gamma_1).
\end{equation}
Then $\kappa$ extends to a smooth diffeomorphism $F$
between collars of $\Gamma_1$ and $\Gamma_2$ such that:
\begin{enumerate}[label=(\roman*)]
\item If $n\neq2p+2$,
\begin{equation}\label{equintrojetisometry}
 j_{\Gamma_1}^\infty(F^*g_2)=j_{\Gamma_1}^\infty g_1.
\end{equation}
\item If $n=2p+2$, for some real-valued smooth $f$ near $\Gamma_1$ with
$f|_{\Gamma_1}=0$,
\begin{equation}\label{jetintroductionconformal}
 j_{\Gamma_1}^\infty(F^*g_2)=j_{\Gamma_1}^\infty(e^{2f}g_1).
\end{equation}
\end{enumerate}
\end{theorem}

\begin{remark}
In the critical dimension $n=2p+2$, a smooth conformal change on $(M,g)$ satisfies (see also \cite[Proposition~1.3]{SharafutdinovLinearized})
\begin{equation}\label{criticconformalinvariance}
 \Lambda_{e^{2f}g}\varphi=e^{-f|_Y}\Lambda_g\varphi,
 \qquad f\in C^\infty(M;\R).
\end{equation}
Indeed, on $(p+1)$-forms,
\begin{equation}\label{confomaldelta}
 \delta_{e^{2f}g}\eta
 =-*_{e^{2f}g}d*_{e^{2f}g}\eta
 =-e^{-2f}*_g d*_g\eta
 =e^{-2f}\delta_g\eta.
\end{equation}
If $\omega$ and $\widehat\omega$ solve \eqref{eoriginalproblem}
for $g$ and $e^{2f}g$ with the same tangential trace, then
$\delta_gd\omega=\delta_gd\widehat\omega=0$.
Green's formula gives
$\|d(\widehat\omega-\omega)\|_{L^2(M,g)}^2=0$.
Equation~\eqref{criticconformalinvariance} follows from
$\nu_{e^{2f}g}=e^{-f|_Y}\nu_g$.
Thus, in the critical dimension, conformal changes $g\mapsto e^{2f}g$
with $f|_Y=0$ leave $\Lambda_g$ unchanged. This invariance obstructs
recovery of the full boundary metric jet, showing that
\eqref{jetintroductionconformal} is optimal.
\end{remark}
The proof of Theorem~\ref{thboundaryintrodo} proceeds by induction, directly comparing the first unmatched
jet at each step.

To describe the mechanism, suppose that the lower normal jets of the two metrics in Theorem~\ref{thboundaryintrodo} agree and that, in corresponding
boundary normal coordinates,
\begin{equation*}
 g_2-g_1=\frac{r^\ell}{\ell!}\mathfrak h+O(r^{\ell+1}),\qquad \ell\geq1.
\end{equation*}
The homogeneous symbol of degree $1-\ell$ of the difference between
the boundary operators associated with $g_2$ and $g_1$, restricted
to $\ker\iota_{\xi^\sharp}$, is
\begin{equation}\label{vanishingsymbolintro}
 \frac{|\xi|_h^{1-\ell}}{2^{\ell+1}}
 \left( \tr_{\xi^\perp}(h^{-1}\mathfrak h)\,\Id
   -2\left((h^{-1}\mathfrak h)|_{\xi^\perp}\right)^{[p]}\right).
\end{equation}
Here and below, for an endomorphism $P$ of a vector space $V$,
we denote by $P^{[k]}$ the operator on $\Lambda^kV^*$ defined by
\begin{equation*}
 (P^{[k]}\eta)(X_1,\ldots,X_k)
 :=\sum_{j=1}^k\eta(X_1,\ldots,PX_j,\ldots,X_k),
 \qquad k\geq1,
\end{equation*}
and set $P^{[k]}=0$ for $k\leq0$.

We show that the vanishing of \eqref{vanishingsymbolintro} for every $\xi\neq0$
forces $\mathfrak h=0$ when $n\neq2p+2$, and
$\mathfrak h=c h$ when $n=2p+2$.
The latter ambiguity can be absorbed by a conformal change
at each step of the induction.

\subsection{Partial-data uniqueness for real-analytic metrics}
For the scalar case $p=0$, the global uniqueness results
are classical. Lassas--Uhlmann \cite{LassasUhlmann2001}
established uniqueness from measurements on a nonempty
open boundary set, up to isometry for real-analytic
metrics in dimensions $n\geq3$ and up to conformal
equivalence for smooth surfaces when $n=2$.
Lassas--Taylor--Uhlmann \cite{LassasTaylorUhlmann2003}
extended the real-analytic result to complete manifolds
with compact boundary. An approach based on Poisson
embeddings was developed by Lassas--Liimatainen--Salo \cite{LassasLiimatainenSalo2020}.
 We obtain an analogous uniqueness result for differential forms by combining the collar identification with a new reconstruction from the projection kernel.

In the analytic setting, boundary jet determination gives an initial
geometric identification near the observed set. To extend this
identification over the whole manifold, we use equality of the partial
operators themselves, rather than merely agreement of their full symbols.
Let $\kappa:\Gamma_1\to\Gamma_2$ identify nonempty open boundary sets.
Equality of the partial data means
\begin{equation}\label{equDNintropar}
 (\Lambda_{g_1}\varphi)|_{\Gamma_1}
 =\kappa^*\left( \left.\Lambda_{g_2}\bigl((\kappa^{-1})^*\varphi\bigr)\right|_{\Gamma_2}\right),\qquad \varphi\in\Omega_c^p(\Gamma_1),
\end{equation}
where each input is extended by zero to the corresponding boundary.

\begin{theorem}
\label{thmsnalyticcalderon}
Let $(M_j,g_j)$, $j=1,2$, be compact connected oriented
smooth Riemannian $n$-manifolds with nonempty smooth boundaries.
Let $1\leq p\leq n-2$, and let
$\Gamma_j\subset\partial M_j$ be nonempty domains.
Assume that $g_j$ are real-analytic in $M_j^\circ$,
and that both the metrics and boundaries are real-analytic
near $\Gamma_j$, in compatible analytic structures.
Let $\kappa:\Gamma_1\to\Gamma_2$ be a smooth
diffeomorphism satisfying \eqref{equDNintropar}.

Then $\kappa$ extends to a smooth diffeomorphism
$F:M_1\to M_2$, real-analytic in the interior and
up to $\Gamma_1$, with the following properties:
\begin{enumerate}[label=(\roman*)]
\item If $n\neq2p+2$, then $F^*g_2=g_1$.
\item If $n=2p+2$, then
\[
F^*g_2=e^{2f}g_1,\qquad
f\in C^\infty(M_1;\R),\qquad f|_{\Gamma_1}=0,
\]
where $f$ is real-analytic in $M_1^\circ$ and
up to $\Gamma_1$.
\end{enumerate}
\end{theorem}
It is instructive to compare our data set with the complete Cauchy data used by Krupchyk--Lassas--Uhlmann \cite{KrupchykLassasUhlmann}. Their data consist, in our conventions, of $(\tang\omega,\normal\omega,\normal (d\omega),\tang(\delta\omega))$ for harmonic forms $\omega$ and are used to construct an embedding through the Green kernel. In contrast, our measurements use only $(\tang \omega,\normal (d\omega))$ for coclosed harmonic forms $\omega$ in \eqref{eoriginalproblem}. Thus, when $n\neq2p+2$, we obtain the isometric uniqueness conclusion of \cite[Theorem~1.2]{KrupchykLassasUhlmann} from fewer measured trace components. When $n=2p+2$, \eqref{criticconformalinvariance} shows that only conformal equivalence can be expected.

Working with fewer measured trace components than in
\cite{KrupchykLassasUhlmann} also changes the geometric argument.
In particular, it leads us to construct a different embedding. A major difficulty is that our reduced data do not directly determine the Green kernels. 
We instead recover the kernel $\mathsf K_j$ of the orthogonal projection onto the
$L^2$-closure of relative exact $(p+1)$-forms.
From $\mathsf K_j$ we construct the operator
$\mathcal T_j$, which takes values in a fixed Hilbert space of self-adjoint Hilbert--Schmidt
operators. The diagonal singularities of $\mathsf K_j$ show that
$\mathcal T_j$ is injective and has injective differential in
the interior. Its analyticity follows from analytic elliptic
regularity.

The reduced data also require a new way to control the unobserved boundary: the projection kernel $\mathsf K_j$ used here need not vanish at the unobserved boundary, so the boundary exclusion argument in \cite[Lemma~3.8]{KrupchykLassasUhlmann} is not available.
The key new boundary mechanism is a rank drop:
\[
 \rank\mathcal T_j(x)=\binom{n}{p+1}\quad (x\in M_j^\circ),\qquad
 \rank\mathcal T_j(x)=\binom{n-1}{p}\quad (x\in\partial M_j\setminus\Gamma_j).
\]
This difference prevents the continued map from sending
interior points to the unobserved boundary. However, interior analyticity does not by
itself ensure smooth extension to the unobserved boundary,
where no analyticity is assumed. To bridge the gap, we use the small
positive eigenvalues of $\mathcal T_j$ near the boundary to construct a smooth boundary defining function.
Together with tangential
coordinates obtained from $\mathcal T_j$, this function gives boundary coordinates in which both $F$ and its inverse are smooth.

\subsection{Partial-data uniqueness for metrics with harmonic curvature}

The preceding continuation result separates the global argument
from the construction of the observed collar. This allows us to
replace boundary analyticity by a geometric equation.
Recall that a metric has \emph{harmonic curvature} if its Levi--Civita
connection satisfies the Yang--Mills equation
$\nabla^aR_{abcd}=0$. By the contracted Bianchi identity, this is
equivalent to its Ricci tensor being Codazzi.
\begin{remark}
The harmonic-curvature condition includes the following standard
geometries (see also \cite[Sections~2.6--2.8]{DerdzinskiHarmonicCurvature1985}).
\begin{enumerate}[label=(\roman*)]
\item Einstein metrics, including metrics of constant sectional
curvature: $\operatorname{Ric}=\lambda g$ with constant $\lambda$
implies $\nabla\operatorname{Ric}=0$.
\item Locally symmetric metrics, since $\nabla R=0$.
\item Riemannian products of harmonic-curvature metrics: the
product connection acts separately on each factor and the
mixed curvature components vanish. In particular, products of
Einstein metrics are included even when their Einstein constants
differ.
\item Locally conformally flat metrics with constant scalar
curvature. Local conformal flatness makes the Schouten tensor
Codazzi, and constant scalar curvature then makes the Ricci
tensor Codazzi.
\end{enumerate}
\end{remark}
By a theorem of DeTurck--Goldschmidt
\cite{DeTurckGoldschmidt1989}, metrics with harmonic curvature
are real-analytic in suitable interior coordinates.
This regularity does not, however, imply analyticity of the
metric and boundary near the observed set. Thus
Theorem~\ref{thmsnalyticcalderon} cannot be applied directly,
and a separate boundary rigidity argument is needed. 
Guillarmou--S\'a Barreto \cite[Theorem~1.1]{GuillarmouSaBarreto2009}
proved partial-data uniqueness for smooth Einstein
manifolds from the scalar Dirichlet-to-Neumann map. The result
below extends this type of theorem in two directions: from scalar functions to the operator $\Lambda_g$ on differential forms, and from Einstein metrics to the substantially broader class of metrics with harmonic curvature.
\begin{theorem}
\label{thmharmoniccurvaturein}
Let $(M_j,g_j)$, $j=1,2$, be compact connected oriented smooth
$n$-manifolds with nonempty smooth boundaries, where $n\geq3$.
Assume that the metrics are smooth up to the boundary and have
harmonic curvature. Let $0\leq p\leq n-2$, let
$\Gamma_j\subset\partial M_j$ be nonempty domains, and let
$\kappa:\Gamma_1\to\Gamma_2$ be a smooth diffeomorphism satisfying  \eqref{equDNintropar}.
\begin{enumerate}[label=(\roman*)]
\item If $n\neq2p+2$, then $\kappa$ extends to a smooth isometry
$F:(M_1,g_1)\to(M_2,g_2)$.
\item If $n=2p+2$, then $\kappa$ extends to a smooth
diffeomorphism $F:M_1\to M_2$ with
\[
 F^*g_2=e^{2f}g_1,\qquad f\in C^\infty(M_1;\R),\qquad f|_{\Gamma_1}=0.
\]
Moreover, if the Weyl tensor $W_{g_1}$ associated with $(M_1,g_1)$ satisfies 
\begin{equation}\label{ecriticalisometriccondition}
  \iota_{\nu_{g_1}}W_{g_1}\big|_{x_0}\neq0
\end{equation}
for some point $x_0\in\Gamma_1$, then $F$ is an isometry.
\end{enumerate}
\end{theorem}

\begin{remark}\label{remark4dimensionweyl}
  In dimension four, we have
(see \cite[Eq.~(31)]{EdgarHoglund2002})
\begin{equation*}
 (W_{g_1})_{ik\ell m}(W_{g_1})_q{}^{k\ell m}
 =\frac14|W_{g_1}|_{g_1}^2(g_1)_{iq}.
\end{equation*}
Hence \eqref{ecriticalisometriccondition} is equivalent to
$W_{g_1}(x_0)\neq0$.
If $g_1$ is not locally conformally flat,
Lemma~\ref{lemmaweylboundaryuc} ensures that such a point exists
in every nonempty observed boundary open set.
\end{remark}
We briefly outline the proof of
Theorem~\ref{thmharmoniccurvaturein}.
The main difficulty is to pass from matching boundary jets to an exact collar map
without boundary analyticity. In the noncritical case, we
overcome this difficulty by exploiting the harmonic curvature
equation. We couple an elliptic equation for the curvature with normal transport
equations for an orthonormal frame and its connection.
A Carleman estimate then establishes collar rigidity from
metric jets (Proposition~\ref{lemmacurvatureboundaryuc}), replacing analytic
continuation at the observed boundary.

The critical dimension $n=2p+2\geq4$ presents a further
barrier. Boundary determination identifies the metric jets
only up to a conformal factor, while a general conformal transform 
does not preserve harmonic curvature. Thus the preceding isometric rigidity result is not directly
applicable in general. We distinguish two cases.
If the conformal factor vanishes to all orders at the boundary,
Proposition~\ref{lemmacurvatureboundaryuc} applies directly.
Otherwise, the Cotton and scalar curvature identities reveal
an additional constraint on a smaller observed boundary patch:
\begin{equation}\label{constraintWeylintro}
 W_{g_j}(\partial_{r_j},\cdot,\cdot,\cdot)
 =O(r_j^\infty),\qquad j=1,2,
\end{equation}
where $r_j$ are boundary normal coordinates.

Our main construction turns this constraint into an explicit
reconstruction of the geometry from boundary values.
Conformal geodesics provide coordinates and transported frames
that remain unchanged under conformal rescaling with compatible
initial data
(Proposition~\ref{lemmacriticalcoordinates}).
The matching conformal jets then identify the tangential
boundary values of the associated connection matrices with a
common matrix $\mathbf B_0$.
Proposition~\ref{propcriticalradialmodel} proves that, in these
coordinates, the connection matrices in the collar satisfy
\[
 \mathbf B_j
 =\mathbf X_0\,dt
  +e^{-t\mathbf X_0}\mathbf B_0e^{t\mathbf X_0},
 \qquad j=1,2,
\]
where $\mathbf X_0$ is the fixed matrix in
\eqref{ecriticalradialvalues}.
We prove Proposition~\ref{propcriticalradialmodel} by comparing
the connection matrix with the explicit expression above.
The key step is to couple an elliptic equation for the difference
of the Weyl tensors with normal transport equations.
The infinite-order boundary constraint \eqref{constraintWeylintro}
gives zero Cauchy data for this coupled system, so boundary unique
continuation forces the connection difference to vanish
throughout a collar.

The reconstruction therefore makes the two connection matrices,
and hence their coframes, coincide. This converts matching conformal boundary jets at a smooth
boundary into an exact conformal collar map
(Proposition~\ref{propcriticalconformalcollar}), supplying the
local identification needed for global continuation.

\subsection*{Plan of the paper}
Section~\ref{sectiontangeormal} derives the full symbol recursion for
$\Lambda_g$. Section~\ref{boundaryinverse} extracts the successive normal
metric jets and identifies the critical conformal kernel, proving
Theorem~\ref{thboundaryintrodo}.
Section~\ref{sectionanalyticinverse} develops the projection-kernel
embedding and proves the global uniqueness result in
Theorem~\ref{thmsnalyticcalderon}, including smooth extension
across the unobserved boundary. Section~\ref{sectioncurvatureinverse}
establishes boundary unique continuation and shows that matching metric
jets through order two determine a metric with harmonic curvature near
the boundary, providing the local rigidity needed to prove
Theorem~\ref{thmharmoniccurvaturein}(i).
Section~\ref{sectioncriticalinverse} uses conformal geodesics
and connection matrices to construct conformal collar maps in every
critical dimension and completes the proof of
Theorem~\ref{thmharmoniccurvaturein}.

\section{Boundary reduction and full-symbol recursion}
\label{sectiontangeormal}

\subsection{Tangential--normal decomposition}Throughout this section, exterior powers and spaces of differential forms of negative degree are understood to be zero.

For sufficiently small $\varepsilon>0$, the inward normal collar map
$\Phi(r,y):=\exp_y^g(-r\nu_g(y))$ identifies a neighborhood of $Y$ with
$[0,\varepsilon)\times Y$. Writing $Y_r:=\Phi(\{r\}\times Y)$, we have
$g=dr^2+h_r$, $h_0=h$, and $\nu_g=-\partial_r$ on $Y$.
Tangential tensors on $Y_r$ may be identified with tensors on $Y$
through $\Phi$, with $r$ regarded as a parameter. Every smooth $p$-form $\omega$ in the collar admits a unique
decomposition
\begin{equation*}
\omega=\alpha+dr\wedge\beta,
\end{equation*}
where, for every fixed $r$,
\begin{equation*}
\alpha|_{Y_r}\in\Omega^p(Y_r),\qquad\beta|_{Y_r}\in\Omega^{p-1}(Y_r).
\end{equation*}
Equivalently,
\begin{equation*}
\beta=\iota_{\partial_r}\omega,\qquad\alpha=\omega-dr\wedge\iota_{\partial_r}\omega.
\end{equation*}

The shape operator of $Y_r$ is defined by
\begin{equation*}
 S_r:TY_r\longrightarrow TY_r,\qquad S_r(X):=-\nabla_X\partial_r.
\end{equation*}
The induced metrics $h_r$ satisfy
\begin{equation}\label{emetricshaperelation}
 \partial_r h_r(X,Z)=-2h_r(S_rX,Z).
\end{equation}
The corresponding unnormalized mean curvature is
\begin{equation*}
 H_r:=\operatorname{tr}S_r=-\frac12\operatorname{tr}_{h_r}(\partial_r h_r).
\end{equation*}
The normal covariant derivative on tangential
$(p-1)$-forms takes the form
$\nabla_{\partial_r}=\partial_r+S_r^{[p-1]}$.

Let $\eta(r)$ and $\zeta(r)$ be smooth families of tangential
$(p-1)$-forms. Metric compatibility and self-adjointness of
$S_r^{[p-1]}$ imply
\[
 \partial_r\langle\eta,\zeta\rangle_{h_r}=\langle\partial_r\eta,\zeta\rangle_{h_r}+\langle\eta,\partial_r\zeta\rangle_{h_r}+2\langle S_r^{[p-1]}\eta,\zeta\rangle_{h_r}.
\]
Moreover, \eqref{emetricshaperelation} gives
\[
 \partial_r dV_{h_r}=\frac12\operatorname{tr}_{h_r}(\partial_r h_r)\,dV_{h_r}=-H_r\,dV_{h_r}.
\]
Taking $\eta$ compactly supported in the interior of the collar,
integration by parts yields
\[
 \int\langle\partial_r\eta,\zeta\rangle_{h_r}\,dr\,dV_{h_r}=\int\left\langle\eta,
 -\partial_r\zeta+(H_r\Id-2S_r^{[p-1]})\zeta\right\rangle_{h_r}\,dr\,dV_{h_r}.
\]
Thus the formal adjoint of $\partial_r$ is
$\partial_r^*=-\partial_r+H_r\Id-2S_r^{[p-1]}$.
The exterior derivative and its formal adjoint have the collar
expressions
\begin{align}
 d\omega &=d_{Y_r}\alpha +dr\wedge\bigl(\partial_r\alpha-d_{Y_r}\beta\bigr),\label{ecollard}\\
 \delta\omega&=\delta_{Y_r}\alpha-\partial_r\beta +\bigl(H_r\Id-2S_r^{[p-1]}\bigr)\beta -dr\wedge\delta_{Y_r}\beta.\label{ecollardelta}
\end{align}
Equivalently, the second formula follows from
\cite[Lemma~18(ii)]{RaulotSavo} after replacing the normal covariant
derivative by $\partial_r+S_r^{[p-1]}$.
In particular, $\tang\omega=\alpha(0)$ and
$\normal\omega=-\beta(0)$, and
\begin{align}
 \normal(d\omega) &=-\partial_r\alpha(0)+d_Y\beta(0),\label{boundarynormald}\\
 \tang(\delta\omega)&=\delta_Y\alpha(0)-\partial_r\beta(0)
   +\bigl(H_0\Id-2S_0^{[p-1]}\bigr)\beta(0).\label{boundarytangentialdelta}
\end{align}
\subsection{Relative and full Dirichlet boundary problems}

The system \eqref{eoriginalproblem} used to define $\Lambda_g$ can be replaced by the standard relative
boundary problem
\begin{equation}\label{equationrelativebvp}
 \Delta\omega=0,\qquad \tang\omega=\varphi,\qquad \tang(\delta\omega)=0.
\end{equation}
\begin{proposition}[\cite{Karpukhin}, Proposition~3.9]
\label{lemmadeltazero}
Any solution of $\Delta\omega=0$ and $\tang(\delta\omega)=0$
satisfies $\delta\omega=0$. Hence \eqref{equationrelativebvp} is
equivalent to \eqref{eoriginalproblem}.
\end{proposition}

Set
\[
 V^p:=\Lambda^pT^*M\restr_Y \cong\Lambda^pT^*Y\oplus dr\wedge\Lambda^{p-1}T^*Y.
\]
Given $(\alpha_0,\beta_0)\in\Omega^p(Y)\oplus\Omega^{p-1}(Y)$, let
$u=\alpha+dr\wedge\beta$ solve the following full Dirichlet boundary problem:
\begin{equation}\label{eulldirichlet}
 \Delta u=0,\qquad \alpha(0)=\alpha_0,\qquad \beta(0)=\beta_0.
\end{equation}
 Standard elliptic theory
provides a unique smooth solution of \eqref{eulldirichlet} and a
continuous Dirichlet inverse on smooth interior $p$-forms
\cite{Schwarz,Grubb}. Define the outward normal derivative of the coefficient vector by
\begin{equation}\label{eqNdef}
 \cN_p(\alpha_0,\beta_0):=\bigl(-\partial_r\alpha(0),-\partial_r\beta(0)\bigr).
\end{equation}

\subsection{Normal factorization of the Hodge Laplacian}

In boundary normal coordinates $(r,y)$, it follows from \eqref{ecollard} and
\eqref{ecollardelta} that
\begin{equation*}
 \Delta=-\partial_r^2\Id+E(r,y)\partial_r+Q(r,y,D_y),\qquad \sigma_2(Q)(r,y,\xi)=|\xi|_{h_r}^2\Id,
\end{equation*}
where $E$ is a smooth bundle endomorphism and $Q$ is a family of
second-order tangential differential operators. With respect to
$\omega=\alpha+dr\wedge\beta$, their blocks are
\begin{equation}\label{eEcoordinateblocks}
 E=\begin{pmatrix}
 H_r\Id-2S_r^{[p]}&0\\
 0&H_r\Id-2S_r^{[p-1]}
 \end{pmatrix},
\end{equation}
and
\begin{equation}\label{eQcoordinateblocks}
 Q=
 \begin{pmatrix}
 \Delta_{Y_r}&d_{Y_r}\circ E_{22}-E_{11}\circ d_{Y_r}
 \\[2pt]
 \partial_r(\delta_{Y_r}) &\Delta_{Y_r}+\partial_rE_{22}
 \end{pmatrix}.
\end{equation}
Here $\partial_r(\delta_{Y_r})$ and $\partial_rE_{22}$ mean
derivatives of operator coefficients in the fixed coordinate
coframe.
We compute the symbols of operator compositions using
Kohn--Nirenberg quantization:
\begin{equation}\label{esharpproduct}
 a\#b\sim\sum_{\gamma\in\mathbb N^{n-1}}
 \frac1{\gamma!}\partial_\xi^\gamma a\,D_y^\gamma b,\qquad D_y:=\frac1{\sqrt{-1}}\partial_y.
\end{equation}
The following factorization, due to Joshi and Lionheart, is a bundle-valued
extension of the expansion carried out by Lee and Uhlmann for the scalar Laplace--Beltrami
operator \cite{LeeUhlmann}. We include the argument in order to derive explicitly the symbol equation and the recursive determination of the full symbol needed below.
\begin{proposition}[\cite{JoshiLionheart}, Proposition~3.1]

There is a smooth family $A(r)\in\PsiDO^1(Y;V^p)$ with
\begin{equation*}
 a_1(r,y,\xi)=|\xi|_{h_r}\Id_{V^p}
\end{equation*}
such that
\begin{equation}\label{efactorization}
 \Delta=(\partial_r-E-A)(-\partial_r-A)+\mathcal{R},
 \qquad \mathcal{R}(r)\in\PsiDO^{-\infty}(Y;V^p),
\end{equation}
where $\mathcal{R}(r)$ depends smoothly on $r$.
The full symbol of $A$ is determined recursively by the full symbol of
$\Delta$ in the collar.
\end{proposition}

\begin{proof}
Expanding the product in \eqref{efactorization} and using
$
\partial_r\circ A=(\partial_r A)+A\partial_r,
$ reduces the construction to the Riccati equation
\begin{equation}\label{ericcatioperator}
 A^2-\partial_rA+EA=Q
 \quad\bmod\PsiDO^{-\infty}.
\end{equation}
Let $a\sim\sum_{\ell\geq0}a_{1-\ell}$ and $q=q_2+q_1+q_0$
be the full symbols of $A$ and $Q$. The symbol equation corresponding to
\eqref{ericcatioperator} is
\begin{equation}\label{ericcatisymbol}
 a\#a-\partial_ra+Ea\sim q.
\end{equation}
At degree two it reads
\[
 a_1^2=q_2=|\xi|_{h_r}^2\Id.
\]
We choose the positive root
$
 a_1=|\xi|_{h_r}\Id.
$ For $\ell\geq1$, assume that $a_1,a_0,\ldots,a_{2-\ell}$ have already been
determined. Then the terms of degree $2-\ell$ in $a\#a$
containing the unknown $a_{1-\ell}$ are exactly
$a_1a_{1-\ell}+a_{1-\ell}a_1=2|\xi|_{h_r}a_{1-\ell}$.  Hence
\begin{align}\label{eArecursion}
 2|\xi|_{h_r}a_{1-\ell}=q_{2-\ell}+\partial_ra_{2-\ell}-Ea_{2-\ell}
 -\sum_{\substack{k_1+k_2+|\gamma|=\ell\\
       (k_1,k_2,|\gamma|)\neq(\ell,0,0),(0,\ell,0)}}
 \frac1{\gamma!}\partial_\xi^\gamma a_{1-k_1}D_y^\gamma a_{1-k_2}.
\end{align}
Here $k_1,k_2\geq0$ and $q_{2-\ell}=0$ for $\ell>2$.
Since $|\xi|_{h_r}>0$ for $\xi\neq0$, the recursion determines
every homogeneous coefficient, smoothly in $r$. Borel summation then gives a
classical symbol $a$ satisfying \eqref{ericcatisymbol} modulo
$S^{-\infty}$, with remainder estimates uniform on compact $r$-intervals and
valid after any number of $r$-derivatives. Thus $\mathcal{R}(r)=Q-A^2+\partial_r A-EA$ is a smooth family of smoothing operators. This completes the proof.
\end{proof}
\begin{proposition}\label{czprNispsi}
  The operator $\cN_p$ in \eqref{eqNdef} is a classical
pseudodifferential operator of order one on $V^p$, and
\begin{equation*}
 \cN_p=A(0)\quad\bmod\PsiDO^{-\infty}(Y;V^p).
\end{equation*}
In particular,
\begin{equation*}
 \sigma_1(\cN_p)(y,\xi)=|\xi|_h\Id_{V^p}.
\end{equation*}
\end{proposition}
\begin{proof}
For $p=0$, the conclusion follows from
\cite[Propositions~1.1--1.2]{LeeUhlmann}.
For $p\geq1$, our assumption $p\leq n-2$ implies $n\geq3$,
and the conclusion follows from
\cite[Theorem~1.1(i) and Proposition~3.1]{JoshiLionheart}.
\end{proof}

\subsection{Boundary matrix for relative Hodge data}
For the solution $u=\alpha+dr\wedge\beta$ with boundary value
$(\alpha_0,\beta_0)$ of \eqref{eulldirichlet}, define
\begin{equation}\label{def::dpHodge}
 \cD_p\binom{\alpha_0}{\beta_0}:=\binom{\normal(du)}{\tang(\delta u)}.
\end{equation}
Write $\cN_p=(N_{ij})_{1\leq i,j\leq2}$ with respect to $\Lambda^pT^*Y\oplus dr\wedge\Lambda^{p-1}T^*Y$.
By \eqref{eqNdef}, we have
\[
-\partial_r\alpha(0)=N_{11}\alpha_0+N_{12}\beta_0,\qquad -\partial_r\beta(0)=N_{21}\alpha_0+N_{22}\beta_0.
\]
Combining equations \eqref{boundarynormald} and \eqref{boundarytangentialdelta}, we obtain the exact operator
identity
\begin{equation}\label{eDblocks}
 \cD_p= \begin{pmatrix}
 D_{11}&D_{12}\\
 D_{21}&D_{22}
 \end{pmatrix} = \begin{pmatrix}
 N_{11}&N_{12}+d_Y\\
 N_{21}+\delta_Y&N_{22}+E_{22}(0)
 \end{pmatrix},
\end{equation}
where $E_{22}(0)=H_0\Id-2S_0^{[p-1]}$.
By Proposition~\ref{czprNispsi} and the standard formulas
\[
 \sigma_1(d_Y)=\sqrt{-1}\eps_\xi, \qquad \sigma_1(\delta_Y)=-\sqrt{-1}\iota_{\xi^\sharp},
\]
the operator $\cD_p$ is a classical pseudodifferential operator of order one with principal symbol
\begin{equation}\label{eDprincipalmatrix}
 \sigma_1(\cD_p)(y,\xi)
 =\begin{pmatrix}
 |\xi|_h\Id_{\Lambda^p}&\sqrt{-1}\eps_\xi\\
 -\sqrt{-1}\iota_{\xi^\sharp}&|\xi|_h\Id_{\Lambda^{p-1}}
 \end{pmatrix}.
\end{equation}

\subsection{Elimination of the normal boundary trace}

We now impose the boundary condition
$\tang(\delta u)=0$ in \eqref{equationrelativebvp}.
For a fixed tangential boundary value $\alpha_0$, $\tang(\delta u)=0$ is equivalent to
\begin{equation*}
 D_{21}\alpha_0+D_{22}\beta_0=0.
\end{equation*}
By \eqref{eDprincipalmatrix}, $D_{22}$ is an elliptic classical
pseudodifferential operator of order one.
Choose a classical parametrix
$G\in\PsiDO^{-1}(Y;\Lambda^{p-1}T^*Y)$ with
\begin{equation*}
 D_{22}G=\Id-\mathcal{R}_1,\qquad GD_{22}=\Id-\mathcal{R}_2,
 \qquad \mathcal{R}_1,\mathcal{R}_2\in\PsiDO^{-\infty}.
\end{equation*}
\begin{lemma}
Choose the solution of \eqref{equationrelativebvp} orthogonal in
$L^2(M)$ to $\cH_D^p(M)$, and write its normal coefficient as
$\beta_0=B\alpha_0$. Then there exist $\mathcal{R}_B,\mathcal{R}_\Lambda\in\PsiDO^{-\infty}$ such that
\begin{equation}\label{eBformula}
 B=-GD_{21}+\mathcal{R}_B\in\PsiDO^0(Y;\Lambda^pT^*Y,\Lambda^{p-1}T^*Y),
\end{equation}
\begin{equation}\label{eLambdaSchur}
 \Lambda_g
 =D_{11}-D_{12}GD_{21}+\mathcal{R}_\Lambda\in\PsiDO^1(Y;\Lambda^pT^*Y).
\end{equation}
\end{lemma}
\begin{proof}
This is essentially the content of \cite[Proposition~3.3]{CHL26}.
\end{proof}

\subsection{Full-symbol recursion and normal metric jets}
\label{sec:fullsymbol}
We recursively describe the full symbol of $\Lambda_g$, tracking
normal metric derivatives at each degree. The block formulas
below apply to $p\geq1$. For $p=0$, $\Lambda_g=\cN_0$, so its full symbol follows directly from
\eqref{eArecursion} at $r=0$.

Let $(d)_{ij}(y,\xi)\sim (d_1)_{ij}+(d_0)_{ij}+(d_{-1})_{ij}+\cdots$
be the full symbols of the four blocks in \eqref{eDblocks}.  Write
\[
a_{1-\ell}(0,y,\xi)=\begin{pmatrix}
(a_{1-\ell})_{11}&(a_{1-\ell})_{12}\\
(a_{1-\ell})_{21}&(a_{1-\ell})_{22}
\end{pmatrix}.
\]
By Proposition~\ref{czprNispsi}, we have $\sigma(\cN_p)(y,\xi)\sim a(0,y,\xi)$.
Combining this with \eqref{eDblocks}, the degree-one symbols of
the four blocks are
\begin{equation}\label{eDsymboldegreeone}
d_1=
\begin{pmatrix}
|\xi|_h\Id_{\Lambda^p}&\sqrt{-1}\eps_\xi\\
-\sqrt{-1}\iota_{\xi^\sharp}&|\xi|_h\Id_{\Lambda^{p-1}}
\end{pmatrix}.
\end{equation}
The degree-zero symbols are
\begin{equation}\label{eDsymboldegreezero}
d_0
=
\begin{pmatrix}
(a_0)_{11}&(a_0)_{12}\\
(a_0)_{21}+\sigma_0(\delta_Y)&(a_0)_{22}+E_{22}(0)
\end{pmatrix}_{r=0}.
\end{equation}

It follows from \eqref{eArecursion} that
\begin{equation*}
 a_0=\frac1{2|\xi|_{h_r}}\left(
 q_1+\partial_r|\xi|_{h_r}\Id-E|\xi|_{h_r}
 -\sum_{|\gamma|=1}
   \partial_\xi^\gamma|\xi|_{h_r}\,
   D_y^\gamma|\xi|_{h_r}\Id
 \right).
\end{equation*}
Thus the matrix-valued coefficient $E$ enters by left
multiplication.

The additional terms $d_Y,\delta_Y\in\PsiDO^1$ and
$E_{22}(0)\in\PsiDO^0$ in \eqref{eDblocks} have
no homogeneous components of negative degree.  Hence, for every $\ell\geq2$,
\begin{equation}\label{eDsymbolnegativedegrees}
 d_{1-\ell}=a_{1-\ell}|_{r=0}.
\end{equation}
Thus \eqref{eArecursion} determines successively all the homogeneous
symbols $(d_{1-\ell})_{ij}$.

We next construct a symbolic inverse of the block $(d)_{22}$.
Let $\sigma(G)=:e(y,\xi)\sim e_{-1}+e_{-2}+e_{-3}+\cdots$
satisfy
\[
(d)_{22}\# e\sim\Id,\qquad  e\# (d)_{22}\sim\Id.
\]
By \eqref{eDsymboldegreeone}, we have $e_{-1}=|\xi|_h^{-1}\Id$.
Suppose that $e_{-1},\ldots,e_{-\ell}$ have already been determined.  Taking the
homogeneous component of degree $-\ell$ in $(d)_{22}\# e\sim\Id$ gives
\begin{equation}\label{eerecursion}
e_{-1-\ell}=-(d_1)_{22}^{-1}\mathop{\sum_{k_1+k_2+|\gamma|=\ell}}_{(k_1,k_2,|\gamma|)\neq(0,\ell,0)}
\frac{1}{\gamma!}\partial_\xi^\gamma (d_{1-k_1})_{22}\,D_y^\gamma e_{-1-k_2},
\end{equation}
where $(d_1)_{22}^{-1}=|\xi|_h^{-1}\Id$.

By \eqref{eLambdaSchur},
\[
\sigma(\Lambda_g)(y,\xi)\sim(d)_{11}-(d)_{12}\#\bigl(e\#(d)_{21}\bigr).
\]
Expanding the two sharp products by \eqref{esharpproduct}, we obtain the
following direct recursion for the full symbol of $\Lambda_g$:
\begin{equation}\label{eLambdasymbolrecursion}
\sigma_{1-\ell}(\Lambda_g)=(d_{1-\ell})_{11}-\!\!\sum_{k_1+k_2+k_3+|\alpha|+|\beta|=\ell}\frac{1}{\alpha!\beta!}
\partial_\xi^\alpha (d_{1-k_1})_{12}D_y^\alpha
\left(\partial_\xi^\beta e_{-1-k_2}D_y^\beta (d_{1-k_3})_{21}\right),
\end{equation}
where $k_1,k_2,k_3\geq0$ and
$\alpha,\beta\in\mathbb N^{n-1}$.

Thus the full symbol is obtained recursively in the order
\[
a_{1-\ell}\longrightarrow d_{1-\ell}\longrightarrow
e_{-1-\ell}\longrightarrow\sigma_{1-\ell}(\Lambda_g).
\]
We now record the corresponding dependence on the normal metric jets.
\begin{proposition}
\label{czprLambdajetstructure}
For every $\ell\geq0$, the homogeneous symbol
$\sigma_{1-\ell}(\Lambda_g)$ depends only on the boundary values of
$h_r,\partial_rh_r,\ldots,\partial_r^{\ell}h_r$ and their tangential
derivatives.  For $\ell\geq1$, with the lower normal jets fixed,
$\sigma_{1-\ell}(\Lambda_g)$ depends affinely on
$\partial_r^{\ell}h_r|_{r=0}$.
\end{proposition}

\begin{proof}
 By \eqref{eEcoordinateblocks},
$E$ depends linearly on $\partial_rh_r$ when $h_r$ is fixed since $S_r=-\frac12h_r^{-1}\partial_rh_r$.
The symbol $q_1$ of \eqref{eQcoordinateblocks} uses normal jets through order one and is
affine in $\partial_rh_r$ and its tangential derivatives.
The symbol $q_0$ uses normal jets through order two. Its
second-jet terms come only from $\partial_rE_{22}$ and are
linear in $\partial_r^2h_r$.

We first establish the same jet dependence for the symbol
coefficients $a_{1-\ell}$ of $A$. Since $a_1=|\xi|_{h_r}\Id$
depends only on $h_r$, the case $\ell=0$ holds.
At $\ell=1$, \eqref{eArecursion} expresses $a_0$ in terms of $q_1$,
$\partial_ra_1$, $Ea_1$, and tangential derivatives of $a_1$.
Thus $a_0$ depends affinely on the first normal jet and its
tangential derivatives. At $\ell=2$, the second jet occurs only
in $q_0+\partial_ra_0$, and occurs linearly there. For
$\ell\geq3$, assume the assertion for all preceding coefficients.
Then $q_{2-\ell}=0$, while $Ea_{2-\ell}$ and all terms in the sum
in \eqref{eArecursion} involve normal jets of order at most
$\ell-1$. The $\ell$-th jet can therefore enter only through
$\partial_ra_{2-\ell}$, where it occurs linearly, together with
its tangential derivatives. This proves the claimed filtration
and affine dependence for every $a_{1-\ell}$.

Equations \eqref{eDsymboldegreeone}--\eqref{eDsymbolnegativedegrees}
therefore give the same statement for every $(d_{1-\ell})_{ij}$.  The recursion
\eqref{eerecursion} then proves it for $e_{-1-\ell}$.  Indeed, a term
containing the $\ell$-th normal metric derivative can only have $k_1=\ell$ in
\eqref{eerecursion}, since the term with $k_2=\ell$ is excluded. Hence this
highest jet still occurs linearly.

Finally, in \eqref{eLambdasymbolrecursion} the indices satisfy
$k_1+k_2+k_3+|\alpha|+|\beta|=\ell$.  Thus the $\ell$-th normal metric derivative can
occur in at most one of the factors $(d_{1-k_1})_{12}$, $e_{-1-k_2}$, and
$(d_{1-k_3})_{21}$, and if it does, all the remaining indices vanish.  Hence
$\sigma_{1-\ell}(\Lambda_g)$ depends only on the normal metric jet through order
$\ell$, and its dependence on the highest jet is affine.
\end{proof}

\section{Boundary jet rigidity and the critical conformal gauge}
\label{boundaryinverse}
We use the symbol recursions of Subsection~\ref{sec:fullsymbol}
to isolate the first normal metric derivative not already
determined. This determines the boundary metric to all orders
when $n\neq2p+2$, and up to conformal changes when $n=2p+2$,
proving Theorem~\ref{thboundaryintrodo}.
\subsection{Metric jets and the symbol filtration}
\label{section::commoncollarset}
Let $(M_j,g_j)$, $j=1,2$, be compact oriented smooth Riemannian
$n$-manifolds, and let $\kappa:\Gamma_1\to\Gamma_2$ be a smooth
diffeomorphism between open sets $\Gamma_j\subset\partial M_j$.
Assume that
\[
 g_1|_{\Gamma_1}=\kappa^*(g_2|_{\Gamma_2})=:h.
\]
On each coordinate patch in $\Gamma_1$, choose boundary coordinates
$y$ and use $y\circ\kappa^{-1}$ on the corresponding patch in
$\Gamma_2$. Extend each coordinate system along its inward unit
normal geodesics. In these corresponding boundary normal
coordinates, write
\begin{equation}\label{etwonormalmetrics}
 g_1=dr^2+h_r,\qquad g_2=dr^2+\widehat h_r,
 \qquad h_0=\widehat h_0=h.
\end{equation}
Here $h_r$ and $\widehat h_r$ are expressed in the corresponding
coordinate coframes, so their coefficients and the associated
symbols are compared on $\Gamma_1$.
\begin{lemma}\label{lemmanormaljetfiltration}
Suppose, for some $\ell\geq1$, that
\begin{equation}\label{agreetoell1}
 \partial_r^k h_r|_{r=0}=\partial_r^k\widehat h_r|_{r=0},\qquad 0\leq k\leq\ell-1,\quad\text{on }\Gamma_1.
\end{equation}
Then
\begin{equation*}
 \kappa^*\left.\Lambda_{g_2}\right|_{\Gamma_2}(\kappa^{-1})^*
 -\left.\Lambda_{g_1}\right|_{\Gamma_1}\in\PsiDO^{1-\ell}(\Gamma_1;\Lambda^pT^*\Gamma_1),
\end{equation*}
and its principal symbol is linear in
$\partial_r^\ell(\widehat h_r-h_r)|_{r=0}$.
\end{lemma}

\begin{proof}
Apply Proposition~\ref{czprLambdajetstructure} in the corresponding
boundary normal coordinates. The lower normal jets and their
tangential derivatives agree on $\Gamma_1$, so the symbol terms
of degree greater than $1-\ell$ cancel. The degree $1-\ell$ term
is the difference of two affine expressions with the same lower
jets, and is therefore linear in the first unmatched jet.
\end{proof}

\subsection{The first unmatched metric jet}

The following proposition computes the contribution of the first unmatched
normal metric derivative to the boundary symbol.

\begin{proposition}\label{czprfirstunmatchedjet}
Assume \eqref{agreetoell1} and write
\begin{equation*}
 \widehat h_r-h_r
 =\frac{r^\ell}{\ell!}\mathfrak h+O(r^{\ell+1}).
\end{equation*}
For $(y,\xi)\in T^*\Gamma_1\setminus0$, the leading difference symbol satisfies
\begin{align}
 &\left.\sigma_{1-\ell}\!\left(
   \kappa^*\left.\Lambda_{g_2}\right|_{\Gamma_2}(\kappa^{-1})^*-\left.\Lambda_{g_1}\right|_{\Gamma_1}\right)(y,\xi)\right|_{\ker\iota_{\xi^\sharp}}\notag \\&\quad\quad\quad=\frac{|\xi|_h^{1-\ell}}{2^{\ell+1}}
 \left[ \tr_{\xi^\perp}(h^{-1}\mathfrak h)\Id-2\left((h^{-1}\mathfrak h)|_{\xi^\perp}\right)^{[p]}\right].\label{efirstunmatchedsymbol}
\end{align}
Here the restriction to $\xi^\perp$ denotes orthogonal compression
with respect to $h$, namely
\[
 \left((h^{-1}\mathfrak h)|_{\xi^\perp}\right)X
 =(h^{-1}\mathfrak h)X
 -\frac{\mathfrak h(X,\xi^\sharp)}{|\xi|_h^2}\,\xi^\sharp,
 \qquad X\in\xi^\perp.
\]
\end{proposition}

\begin{proof}
\smallskip\noindent\emph{Step 1. The highest normal jet in the tangential factorization block.}
By \eqref{eEcoordinateblocks} and \eqref{eQcoordinateblocks},
\[
 E_{11}=-\frac12\tr(h_r^{-1}\partial_rh_r)\Id  +(h_r^{-1}\partial_rh_r)^{[p]},\qquad Q_{11}=\Delta_{Y_r}.
\]
In particular, every symbol coefficient of $Q_{11}$ involves
only $h_r$ and its tangential derivatives.
The $(1,1)$ block of \eqref{eArecursion} at $\ell=1$ is
\[
 2|\xi|_{h_r}(a_0)_{11}=(q_1)_{11}+\partial_r|\xi|_{h_r}\Id-E_{11}|\xi|_{h_r}
-\sum_{|\gamma|=1} \partial_\xi^\gamma|\xi|_{h_r}\, D_y^\gamma|\xi|_{h_r}\Id.
\]
The first and last terms use no normal metric derivatives.
Differentiating $|\xi|_{h_r}^2=h_r^{ab}\xi_a\xi_b$ gives
\[
 \partial_r|\xi|_{h_r}=-\frac{\langle(h_r^{-1}\partial_rh_r)\xi^\sharp,\xi^\sharp\rangle_{h_r}}{2|\xi|_{h_r}}.
\]
Consequently, the part of $(a_0)_{11}$ involving the first
normal metric derivative is
\begin{equation}\label{ea0highestnormaljet}
 \frac14\left[ \tr_{\xi^\perp}(h_r^{-1}\partial_rh_r)\Id -2(h_r^{-1}\partial_rh_r)^{[p]}
 \right].
\end{equation}
Here $\xi^\sharp$ and $\xi^\perp$ are taken with respect to
$h_r$ on the slice.

For every $k\geq2$, the $(1,1)$ block of $q_{2-k}$ contains
no normal metric derivative: at $k=2$ it is the degree-zero
symbol of $\Delta_{Y_r}$, and it vanishes for $k>2$.
Moreover, $E$ is block diagonal, so
$(Ea_{2-k})_{11}=E_{11}(a_{2-k})_{11}$.
The jet filtration in the proof of
Proposition~\ref{czprLambdajetstructure} shows that this
product and all symbol products in \eqref{eArecursion}
use normal jets of order at most $k-1$. Thus the $k$-th
normal jet in $(a_{1-k})_{11}$ comes only from
\[
 \frac1{2|\xi|_{h_r}}\partial_r(a_{2-k})_{11}.
\]

For $k\geq2$, applying $\partial_r$ to the highest-jet term of
$(a_{2-k})_{11}$ produces $\partial_r^kh_r$ only when the
derivative falls on $\partial_r^{k-1}h_r$ itself.
Differentiating its coefficients contributes only lower normal jets.
After division by $2|\xi|_{h_r}$ in \eqref{eArecursion},
induction from the base case \eqref{ea0highestnormaljet}
therefore identifies, for every $k\geq1$, the terms in
$(a_{1-k})_{11}$ containing $\partial_r^kh_r$:
\begin{equation}\label{eahighestnormaljet}
 \frac{|\xi|_{h_r}^{1-k}}{2^{k+1}}
 \left[
  \tr_{\xi^\perp}(h_r^{-1}\partial_r^kh_r)\Id-2(h_r^{-1}\partial_r^kh_r)^{[p]}
 \right].
\end{equation}
All other terms involve normal derivatives of order at most
$k-1$ and their tangential derivatives.
\smallskip\noindent\emph{Step 2. Comparison and coclosed compression.}
Since $D_{11}=N_{11}$, equations
\eqref{eDsymboldegreezero} and
\eqref{eDsymbolnegativedegrees} give
$(d_{1-k})_{11}=(a_{1-k})_{11}|_{r=0}$ for every $k\geq1$.
The lower jets in \eqref{agreetoell1} agree,
so \eqref{eahighestnormaljet} yields
\begin{equation}\label{eD11firstunmatchedsymbol}
 (\widehat d_{1-\ell})_{11}-(d_{1-\ell})_{11}
 =\frac{|\xi|_h^{1-\ell}}{2^{\ell+1}}
 \left[\tr_{\xi^\perp}(h^{-1}\mathfrak h)\Id-2(h^{-1}\mathfrak h)^{[p]}
 \right].
\end{equation}
Here and below unhatted and hatted symbol coefficients refer to
$g_1$ and $g_2$, respectively, in the corresponding coordinates
specified above.

For $p\geq1$, apply \eqref{eLambdasymbolrecursion}
with index $\ell$ to $g_1$ and $g_2$, and subtract the
resulting identities.
All factors with jet index less than $\ell$ agree, so a nonzero
difference must occur in a factor with index $\ell$. Thus
 \begin{align*}
 \sigma_{1-\ell}\!\left(
   \kappa^*\left.\Lambda_{g_2}\right|_{\Gamma_2}(\kappa^{-1})^*
   -\left.\Lambda_{g_1}\right|_{\Gamma_1}\right) &=(\widehat d_{1-\ell})_{11}-(d_{1-\ell})_{11}\\
 &\quad-\bigl((\widehat d_{1-\ell})_{12}-(d_{1-\ell})_{12}\bigr)
   e_{-1}(d_1)_{21}\\
 &\quad-(d_1)_{12}(\widehat e_{-1-\ell}-e_{-1-\ell})(d_1)_{21}\\
 &\quad-(d_1)_{12}e_{-1}
          \bigl((\widehat d_{1-\ell})_{21}-(d_{1-\ell})_{21}\bigr).
 \end{align*}
The common principal symbols are
\[
 (d_1)_{12}=\sqrt{-1}\eps_\xi,\qquad (d_1)_{21}=-\sqrt{-1}\iota_{\xi^\sharp},\qquad e_{-1}=|\xi|_h^{-1}\Id.
\]
Since $\iota_{\xi^\sharp}\pi_{\cc}=0$ and
$\pi_{\cc}\eps_\xi=0$, we have
\begin{equation}\label{ecompressedfirstunmatchedsymbol}
 \pi_{\cc}\sigma_{1-\ell}\!\left(
   \kappa^*\left.\Lambda_{g_2}\right|_{\Gamma_2}(\kappa^{-1})^*
   -\left.\Lambda_{g_1}\right|_{\Gamma_1}\right)\pi_{\cc}
=\pi_{\cc}\bigl((\widehat d_{1-\ell})_{11}
                       -(d_{1-\ell})_{11}\bigr)\pi_{\cc}.
\end{equation}
For $p=0$, this identity holds with $\pi_{\cc}=\Id$.
The two boundary metrics agree under $\kappa$, so the transported
coclosed projections have the same full symbol on $\Gamma_1$.
It follows from \eqref{Toeidentity} and
Lemma~\ref{lemmanormaljetfiltration} that
\begin{equation}\label{toeplitzsymbol}
 \sigma_{1-\ell}\!\left( \kappa^*\left.\Lambda_{g_2}\right|_{\Gamma_2}(\kappa^{-1})^*-\left.\Lambda_{g_1}\right|_{\Gamma_1}\right)=\pi_{\cc}\sigma_{1-\ell}\!\left(
   \kappa^*\left.\Lambda_{g_2}\right|_{\Gamma_2}(\kappa^{-1})^*-\left.\Lambda_{g_1}\right|_{\Gamma_1}\right)\pi_{\cc}.
\end{equation}

For $\omega\in\ker\iota_{\xi^\sharp}$ and
$X_1,\ldots,X_p\in\xi^\perp$, we have $\pi_{\cc}\omega=\omega$.
Moreover, the component of
$(h^{-1}\mathfrak h)X_s$ along $\xi^\sharp$ contributes zero
because $\iota_{\xi^\sharp}\omega=0$. Hence
\begin{align*}
 \bigl(\pi_{\cc}(h^{-1}\mathfrak h)^{[p]}\pi_{\cc}\omega\bigr)(X_1,\ldots,X_p)
 &=\sum_{s=1}^p
 \omega\bigl(X_1,\ldots,
 ((h^{-1}\mathfrak h)|_{\xi^\perp})X_s,\ldots,X_p\bigr)\\
 &=\left(\left((h^{-1}\mathfrak h)|_{\xi^\perp}\right)^{[p]}
 \omega\right)(X_1,\ldots,X_p).
\end{align*}
Combining this identity with \eqref{eD11firstunmatchedsymbol}, \eqref{ecompressedfirstunmatchedsymbol}, and \eqref{toeplitzsymbol}
proves \eqref{efirstunmatchedsymbol}.
\end{proof}

\subsection{Algebraic rigidity and the critical kernel}
\begin{lemma}\label{lemmaalgebraicrigidity}
Let $V$ be a Euclidean vector space of dimension $n-1$, and let
$K:V\to V$ be self-adjoint. For a unit vector $v$, let $K_v:v^\perp\to v^\perp$ be its
orthogonal compression to $v^\perp$, namely
\begin{equation*}
 K_vX=KX-\langle KX,v\rangle v,
 \qquad X\in v^\perp.
\end{equation*}
Suppose that
\begin{equation}\label{ealgebraicvanishing}
 \tr(K_v)\Id-2K_v^{[p]}=0
 \quad\text{on }\Lambda^p(v^\perp)^*
\end{equation}
for every unit $v$.
\begin{enumerate}[label=(\roman*)]
\item If $n\neq2p+2$, then $K=0$.
\item If $n=2p+2$, then
 $K=c\Id_V$ for some $c\in\R$.
\end{enumerate}
Conversely, every scalar $K$ satisfies
\eqref{ealgebraicvanishing} when $n=2p+2$.
\end{lemma}

\begin{proof}
Choose an orthonormal eigenbasis
$\theta_1,\ldots,\theta_{n-1}$ of $K$,
with $K\theta_k=\lambda_k\theta_k$, and let
$\theta^1,\ldots,\theta^{n-1}$ be the dual covectors.
Suppose $p\leq n-3$. For any $i\neq j$, choose
\begin{equation*}
 I\subset\{1,\ldots,n-1\}\setminus\{i,j\},
 \qquad |I|=p.
\end{equation*}
For both $v=\theta_i$ and $v=\theta_j$, the $p$-form
$\bigwedge_{k\in I}\theta^k$ is a nonzero eigenvector of
$K_v^{[p]}$ with eigenvalue $\sum_{k\in I}\lambda_k$. Since
$\tr K_{\theta_i}=\tr K-\lambda_i$ and
$\tr K_{\theta_j}=\tr K-\lambda_j$,
\eqref{ealgebraicvanishing} gives, respectively,
\begin{equation*}
 \tr K-\lambda_i-2\sum_{k\in I}\lambda_k=0,\qquad\tr K-\lambda_j-2\sum_{k\in I}\lambda_k=0.
\end{equation*}
Thus $\lambda_i=\lambda_j$. If $p=n-2$, the induced action on the top
exterior power of $v^\perp$ is multiplication by $\tr K_v$, so
\eqref{ealgebraicvanishing} reduces to
$\tr K-\lambda_j=0$ for every $j$, again making all eigenvalues equal.
When $n=2$, $V$ is one-dimensional and this conclusion is automatic.
Hence $K=c\Id_V$ in every case. Substitution into
\eqref{ealgebraicvanishing} gives $(n-2-2p)c=0$, proving all claims.
\end{proof}

\subsection{Proof of Theorem~\ref{thboundaryintrodo}}

We combine the symbol comparison with the algebraic lemma to prove
Theorem~\ref{thboundaryintrodo} by induction on the order of the
normal metric derivatives.

\begin{proof}[Proof of Theorem~\ref{thboundaryintrodo}]
The difference in \eqref{equsntrosmoothingdiff} has a smooth
kernel on $\Gamma_1\times\Gamma_1$. Thus all homogeneous symbol
terms of the two operators agree on
$T^*\Gamma_1\setminus0$. We compare them in boundary coordinates
$y$ on $\Gamma_1$ and $y\circ\kappa^{-1}$ on $\Gamma_2$.

The principal symbol first determines the boundary metric. For any
metric $g$ with induced boundary metric $h$, the principal symbol formula in \cite{CHL26}
and \eqref{rankpcc} give
\begin{equation}\label{erecoverboundarynorm}
 \tr\sigma_1(\Lambda_g)(y,\xi)
 =\binom{n-2}{p}|\xi|_h.
\end{equation}
Thus the boundary covector norms, and hence the induced metrics,
agree under $\kappa$ on $\Gamma_1$.
Extend these corresponding coordinates inward along the normal
geodesics and write the metrics as in \eqref{etwonormalmetrics}.
For $\ell\geq1$, assume inductively that the normal jets agree
through order $\ell-1$ on $\Gamma_1$. The induction starts at
$\ell=1$ by \eqref{erecoverboundarynorm}.

If $n\neq2p+2$, since the degree $1-\ell$ difference symbol is zero,
Proposition~\ref{czprfirstunmatchedjet} and
Lemma~\ref{lemmaalgebraicrigidity} imply $\mathfrak h=0$.
Induction therefore gives
\begin{equation*}
 \left.\partial_r^k(\widehat h_r-h_r)\right|_{r=0}=0
 \quad\text{on }\Gamma_1,\qquad k\geq0.
\end{equation*}
Now define
\begin{equation*}
 F(r,y):=(r,\kappa(y)),\qquad y\in\Gamma_1,
\end{equation*}
for sufficiently small $r\geq0$.
Then $F$ is a smooth collar diffeomorphism with
$F|_{\Gamma_1}=\kappa$ and
$j_{\Gamma_1}^\infty(F^*g_2-g_1)=0$, proving
\eqref{equintrojetisometry}.

Suppose now that $n=2p+2$.
At each step, write the current conformal representative of $g_1$
as $dr^2+h_r$ in its boundary normal coordinates, and compare it
with $g_2=dr^2+\widehat h_r$ in the corresponding coordinates.
By Proposition~\ref{czprfirstunmatchedjet} and
Lemma~\ref{lemmaalgebraicrigidity},
\begin{equation*}
    \mathfrak h=c_\ell h.
\end{equation*}
The coefficient is smooth on all of $\Gamma_1$, since
\begin{equation*}
 c_\ell=\frac1{n-1}\tr_h\left(\left.\partial_r^\ell(\widehat h_r-h_r)\right|_{r=0}\right).
\end{equation*}
On any open boundary set $\Gamma_0$ compactly contained in $\Gamma_1$, take a smooth
metric $\widetilde g_1$ on $M_1$, conformal to $g_1$, with
\begin{equation*}
 \widetilde g_1=e^{c_\ell r^\ell/\ell!}(dr^2+h_r)
 \quad\text{near }\Gamma_0,\qquad
 \widetilde g_1|_{\partial M_1}=g_1|_{\partial M_1}.
\end{equation*}
By \eqref{criticconformalinvariance},
$\Lambda_{\widetilde g_1}=\Lambda_{g_1}$, so
\eqref{equsntrosmoothingdiff} also holds for
$\widetilde g_1$.

In the current coordinates $(r,y)$ near $\Gamma_0$, its
Christoffel symbols satisfy
\begin{equation*}
 \widetilde\Gamma^r_{rr}
 =\frac{c_\ell r^{\ell-1}}{2(\ell-1)!},
 \qquad
 \widetilde\Gamma^a_{rr}
 =-\frac{r^\ell}{2\ell!}h_r^{ab}\partial_b c_\ell.
\end{equation*}
The normal geodesic equations, with initial data
$(r,y;\dot r,\dot y)=(0,y_0;1,0)$, give
 \begin{align*}
 \ddot r(t)=O(t^{\ell-1})&\quad\Longrightarrow\quad r(t)=t+O(t^{\ell+1}),\\
 \ddot y(t)=O(t^\ell)&\quad\Longrightarrow\quad y(t)=y_0+O(t^{\ell+2}).
 \end{align*}
Thus $(t,y_0)$ are boundary normal coordinates for
$\widetilde g_1$. In these coordinates, with $g_2$ written
in the corresponding normal coordinates,
\begin{equation*}
 \widetilde g_1
 =dt^2+h_t+\frac{c_\ell t^\ell}{\ell!}h+O(t^{\ell+1})
 =g_2+O(t^{\ell+1}),
\end{equation*}
where $h_t=h_r|_{r=t}$. Hence the normal metric jets agree
through order $\ell$.
Since $\Gamma_0\Subset\Gamma_1$ is arbitrary and $c_\ell$
does not depend on its choice, the accumulated logarithmic
conformal factor is formally
$ \sum_{\ell=1}^{\infty}{c_\ell(y_\ell)r_\ell^\ell}/{\ell!},$
where $(r_\ell,y_\ell)$ are the boundary normal coordinates
used at step $\ell$.
We express each term in the original collar coordinates and
apply Borel's theorem to realize the resulting Taylor series
as
\[
2f\sim \sum_{\ell=1}^{\infty}\frac{c_\ell\big(y_\ell(r,y)\big)}{\ell!}\big(r_\ell(r,y)\big)^\ell,
\]
with $f$ smooth near $\Gamma_1$ and
$f|_{\Gamma_1}=0$. Then $e^{2f}g_1$ and $g_2$ have matching jets
in corresponding normal coordinates.

For $y\in\Gamma_1$ and sufficiently small $r\geq0$, define
\begin{equation*}
 F\!\left(
   \exp_y^{e^{2f}g_1}\!\left(-r\nu_{e^{2f}g_1}(y)\right)
 \right)
 =\exp_{\kappa(y)}^{g_2}
   \!\left(-r\nu_{g_2}(\kappa(y))\right).
\end{equation*}
This defines a smooth
collar diffeomorphism satisfying
\begin{equation*}
 F|_{\Gamma_1}=\kappa,\qquad j_{\Gamma_1}^\infty\bigl(F^*g_2-e^{2f}g_1\bigr)=0.\qedhere
\end{equation*}
\end{proof}

\section{Partial-data uniqueness for real-analytic metrics}
\label{sectionanalyticinverse}
We prove Theorem~\ref{thmsnalyticcalderon} by first constructing
a collar isometry or conformal map from the boundary jets determined
in Section~\ref{boundaryinverse}.
We then establish
Proposition~\ref{rsmoothboundauation},
which extends this local map globally using only interior
analyticity.
\subsection{Local analytic collar rigidity}

By \eqref{equDNintropar}, the difference in
\eqref{equsntrosmoothingdiff} is zero.
The zeroth-order conclusion of Theorem~\ref{thboundaryintrodo} therefore identifies the induced
metrics under $\kappa$ throughout the observed set.
In the noncritical case, it also identifies all metric jets in
corresponding boundary normal coordinates on $\Gamma_1$. Analyticity
then gives an isometric collar map extending $\kappa$.
In the critical dimension $n=2p+2$,
Proposition~\ref{czprfirstunmatchedjet} and
Lemma~\ref{lemmaalgebraicrigidity} show that, once the lower
normal derivatives of the tangential metrics agree, equality
of the full symbols determines the next derivative only
up to a scalar multiple of the boundary metric. We remove this
ambiguity by passing to analytic conformal
representatives satisfying the trace condition in the following
lemma.

\begin{lemma}
\label{lemmaconformalnormalgauge}
Near an analytic boundary patch $\Gamma$, an analytic metric $g$ has a unique
analytic conformal representative $\bar g=e^{2f}g$, with $f=0$ on the
patch, whose boundary normal form satisfies
\begin{equation}\label{econformalnormalform}
 \bar g=d\rho^2+\bar h_\rho,
 \qquad \tr_{\bar h_\rho}(\partial_\rho\bar h_\rho)=0.
\end{equation}
Here $\rho$ is the inward boundary distance for $\bar g$.
\end{lemma}

\begin{proof}We seek a function $\rho$ satisfying
\begin{equation*}
 |d\rho|_{\bar g}=1,\qquad \Delta_{\bar g}\rho=0,
 \qquad \rho|_\Gamma=0,\qquad \partial_{\nu_g}\rho|_\Gamma=-1,
\end{equation*}
where $\nu_g$ is the outward unit normal for $g$. By the first equation,
$f=\log|d\rho|_g$. The conformal transformation formula for the Laplacian reduces the second equation to
\begin{equation}\label{erhoCKequation}
 \Delta_g\rho
 -(n-2)\langle d\log|d\rho|_g,d\rho\rangle_g=0.
\end{equation}
In $g$-normal coordinates the Cauchy data are
$\rho(0,y)=0$ and $\partial_r\rho(0,y)=1$.
At these data the coefficient of $\partial_r^2\rho$ in
\eqref{erhoCKequation} is $-(n-1)$. Indeed, $|d\rho|_g^2=(\partial_r\rho)^2
+h_r^{ab}\partial_a\rho\,\partial_b\rho$. The Cauchy data therefore give
\begin{equation*}
 \left.\langle d\log|d\rho|_g,d\rho\rangle_g\right|_{r=0}=\left.\frac12\partial_r(|d\rho|_g^2)\right|_{r=0}=\left.\partial_r^2\rho\right|_{r=0}.
\end{equation*}
Since $\Delta_g=-\partial_r^2+\cdots$, the coefficient of
$\partial_r^2\rho$ is $-(n-1)\neq0$. Thus $r=0$ is noncharacteristic for
\eqref{erhoCKequation}. Since the metric and the Cauchy data
are analytic, the Cauchy--Kovalevskaya theorem gives a unique
local analytic solution $\rho$. Set $f=\log|d\rho|_g$ and $\bar g=e^{2f}g$.
 Then $f|_\Gamma=0$ and $|d\rho|_{\bar g}=1$.
Transport the boundary coordinates $y$ along the analytic flow of
$\nabla_{\bar g}\rho$. Since $|d\rho|_{\bar g}=1$, the flow parameter
is $\rho$ and $\partial_\rho=\nabla_{\bar g}\rho$. Hence
$\bar g=d\rho^2+\bar h_\rho$, and
\begin{equation*}
 \Delta_{\bar g}\rho
 =-\partial_\rho\log\sqrt{\det\bar h_\rho}
 =-\frac12\tr_{\bar h_\rho}(\partial_\rho\bar h_\rho)=0,
\end{equation*}
which is the trace condition in \eqref{econformalnormalform}.
\end{proof}

Using this trace condition in the critical case, we now obtain the
collar map needed for global continuation.

\begin{proposition}
\label{analyticcollarrigidity}
Let $(M_j,g_j)$, $j=1,2$, be compact oriented smooth
Riemannian $n$-manifolds, and let $0\leq p\leq n-2$.
Let $\Gamma_j\subset\partial M_j$ be nonempty open sets
identified by a smooth diffeomorphism
$\kappa:\Gamma_1\to\Gamma_2$.
Assume \eqref{equDNintropar} and that
$g_j$ and $\partial M_j$ are real-analytic near $\Gamma_j$.
Then $\kappa$ extends to an analytic collar diffeomorphism
$F_0$ near $\Gamma_1$ with the following properties:
\begin{enumerate}[label=(\roman*)]
\item If $n\neq2p+2$, then
\begin{equation}\label{ecollaralternatives}
 F_0^*g_2=g_1.
\end{equation}
\item If $n=2p+2$, then
\begin{equation}\label{ecollarconformalalternative}
 F_0^*g_2=e^{2f_0}g_1,\qquad f_0|_{\Gamma_1}=0,
\end{equation}
for some analytic function $f_0$ near $\Gamma_1$.
\end{enumerate}
\end{proposition}
\begin{proof}
  By Theorem~\ref{thboundaryintrodo}, $\kappa$ is an
isometry between the induced boundary metrics. Since these metrics are real-analytic, $\kappa$ is real-analytic.

In the noncritical case, all normal jets agree in common boundary
normal coordinates, and
\eqref{ecollaralternatives} follows by analyticity.

In the critical case, Lemma~\ref{lemmaconformalnormalgauge}
and uniqueness give analytic representatives
$\bar g_j=e^{2f_j}g_j$ near $\Gamma_j$, with
$f_j|_{\Gamma_j}=0$.
Fix a boundary coordinate patch $\Gamma_0\Subset\Gamma_1$.
Using cutoffs, extend $f_1,f_2$ smoothly to $M_1,M_2$,
unchanged near $\overline{\Gamma_0}$ and
$\kappa(\overline{\Gamma_0})$, respectively, and with
$f_j|_{\partial M_j}=0$.
Then $\Lambda_{\bar g_j}=\Lambda_{g_j}$ by
\eqref{criticconformalinvariance}, and hence
\begin{equation*}
 \left.\Lambda_{\bar g_1}\right|_{\Gamma_0}
 =\kappa^*\left.\Lambda_{\bar g_2}\right|_{\kappa(\Gamma_0)}
      (\kappa^{-1})^*.
\end{equation*}
Thus all homogeneous symbols agree on $\Gamma_0$.
Using $\kappa$ to identify $\Gamma_0$ and $\kappa(\Gamma_0)$,
write $\bar g_j=d\rho^2+\bar h_{j,\rho}$ in corresponding
boundary normal coordinates, with
$\bar h_{1,0}=\bar h_{2,0}=h$.

For $\ell\geq1$, assume that the normal derivatives of
$\bar h_{1,\rho}$ and $\bar h_{2,\rho}$ agree through
order $\ell-1$ on $\Gamma_0$.
Proposition~\ref{czprfirstunmatchedjet} and
Lemma~\ref{lemmaalgebraicrigidity} give
\begin{equation*}
 \left.\partial_\rho^\ell(\bar h_{2,\rho}-\bar h_{1,\rho})\right|_{\rho=0}=c\,h.
\end{equation*}
Writing the trace conditions in
\eqref{econformalnormalform} as
$\bar h_{j,\rho}^{ab}\partial_\rho(\bar h_{j,\rho})_{ab}=0$,
differentiate $\ell-1$ times in $\rho$ and evaluate at $\rho=0$:
\begin{equation*}
 h^{ab}\left.\partial_\rho^\ell (\bar h_{j,\rho})_{ab}\right|_{\rho=0}
 +\left.\sum_{k=1}^{\ell-1}\binom{\ell-1}{k}
 (\partial_\rho^k\bar h_{j,\rho}^{ab})
 \partial_\rho^{\ell-k}(\bar h_{j,\rho})_{ab}
 \right|_{\rho=0}=0.
\end{equation*}
Every term in the sum depends only on normal derivatives
of order at most $\ell-1$, so the two sums agree by induction.
Subtracting gives
\begin{equation*}
 0=h^{ab}\left.\partial_\rho^\ell
       (\bar h_{2,\rho}-\bar h_{1,\rho})_{ab}\right|_{\rho=0}
  =\tr_h(c\,h)=(n-1)c.
\end{equation*}
Hence $c=0$, so the normal derivatives of order $\ell$
also agree. Induction identifies all normal jets on $\Gamma_0$.

For $y\in\Gamma_0$ and sufficiently small $\rho\geq0$, define
\begin{equation*}
 F_0\!\left(
   \exp_y^{\bar g_1}\!\left(-\rho\nu_{\bar g_1}(y)\right) \right)=\exp_{\kappa(y)}^{\bar g_2}
   \!\left(-\rho\nu_{\bar g_2}(\kappa(y))\right).
\end{equation*}
Then $F_0|_{\Gamma_0}=\kappa$, and the matching normal jets
and analyticity give $F_0^*\bar g_2=\bar g_1$ near $\Gamma_0$.
Setting $f_0=f_1-f_2\circ F_0$, we obtain
\begin{equation*}
 F_0^*g_2=e^{-2(f_2\circ F_0)}F_0^*\bar g_2=e^{2f_0}g_1,\qquad f_0|_{\Gamma_0}=0.
\end{equation*}
By uniqueness in Lemma~\ref{lemmaconformalnormalgauge}
and uniqueness of normal geodesics, the maps $F_0$
constructed on different boundary patches agree near
common boundary points, as do the functions $f_0$.
Since $y$ and $\Gamma_0$ are arbitrary,
$F_0$ and $f_0$ are defined near all of $\Gamma_1$.
\end{proof}
\subsection{Exterior extension and the relative exact projection}
\label{section::defiitionsext}
Following \cite{KrupchykLassasUhlmann}, we extend the manifolds
across the observed boundary and construct a kernel, initially
in the smooth setting.

Suppose that a smooth collar
diffeomorphism $F_0$ extending $\kappa:\Gamma_1\to\Gamma_2$ is already given near $\Gamma_1\subset \partial M_1$. Assume \eqref{ecollaralternatives} when
$n\neq2p+2$, or \eqref{ecollarconformalalternative} with smooth
$f_0$ when $n=2p+2$. Note that Proposition~\ref{analyticcollarrigidity}
provides this map under the analytic hypotheses of
Theorem~\ref{thmsnalyticcalderon}.

Denote by $\Gamma\subset \partial M_j$ the observed set $\Gamma_j$ after identification
by $\kappa$.
Choose a boundary coordinate patch $\Gamma_0\Subset\Gamma$ and a
nonzero cutoff function $0\leq\chi\in C_c^\infty(\Gamma_0)$
sufficiently small for the construction below.
Extend $F_0$ and $g_1$ smoothly across $\Gamma_0$, and extend
$f_0$ smoothly as well in the critical case. Define the extension
of $g_2$ across $\Gamma_0$ by the corresponding collar identity
\eqref{ecollaralternatives} or
\eqref{ecollarconformalalternative}.
Use $F_0$ to identify the extended collars, with coordinates $(r,y)$
and $r>0$ inside $M_j$, transporting $\chi$ by $\kappa$. Define
\begin{equation}\label{elocalexteriorextension}
 \begin{aligned}
 \widetilde M_j&:=M_j\cup\{(r,y):y\in\Gamma_0,\ -\chi(y)\leq r\leq0\},\\
 X_j&:=M_j\cup\{(r,y):y\in\Gamma_0,\ -\chi(y)/2\leq r\leq0\}.
 \end{aligned}
\end{equation}
Both manifolds are compact and connected, with smooth boundaries.
Over $\Gamma_0$, their boundaries are $r=-\chi(y)$ and
$r=-\chi(y)/2$, respectively. Outside $\supp\chi$
they coincide with the original boundary. In particular,
\[
 \partial X_j\cap\widetilde M_j^\circ
 =\{(-\chi(y)/2,y):\chi(y)>0\}.
\]
The extended metrics, still denoted by $g_j$, are smooth up to
$\partial\widetilde M_j$.

Set
\[
 U:=\{(r,y):y\in\Gamma_0,\ -\chi(y)<r<0\}.
\]
Choose a nonempty open set $U_0$ with smooth boundary such that
\[
 \overline{U_0}\subset
 \{(r,y):y\in\Gamma_0,\ -\chi(y)<r<-\chi(y)/2\}.
\]
Thus $\overline{U_0}\cap X_j=\varnothing$.

The relative Hodge Laplacian on $(p+1)$-forms has domain
\begin{equation*}
 \{\eta\in H^2(\widetilde M_j;\Lambda^{p+1}T^*\widetilde M_j):
       \tang\eta=0,\; \tang(\delta_{g_j}\eta)=0
       \text{ on }\partial\widetilde M_j\}.
\end{equation*}
Let $G_{\mathrm{rel},j}^{(p+1)}$ be its inverse on
$\cH_D^{p+1}(\widetilde M_j)^\perp$, extended by zero on
$\cH_D^{p+1}(\widetilde M_j)$. 
Hence
\begin{equation}\label{defGrelj}
 \Delta G_{\mathrm{rel},j}^{(p+1)}=\mathrm{Id}-\Pi_{\cH_D^{p+1}(\widetilde M_j)}.
\end{equation}
The Hodge decomposition
\cite[Sections~3.4]{Schwarz} implies that
\begin{equation*}
 \mathsf P_j:=d\delta_{g_j}G_{\mathrm{rel},j}^{(p+1)}
\end{equation*}
is the orthogonal projection onto the $L^2(\widetilde M_j,g_j)$-closure of
\[
 d\{\alpha\in H^1(\widetilde M_j;\Lambda^pT^*\widetilde M_j): \tang\alpha=0\}.
\]

For later use, we record a consequence of the orthogonality
defining $\mathsf P_j$.
Let $\alpha\in L^2(\widetilde M_j)$ with
$\zeta:=d\alpha\in L^2(\widetilde M_j)$.
Assume that $\alpha$ has vanishing tangential trace in the weak
sense, namely
\[
 \langle d\alpha,\psi\rangle_{L^2}=\langle\alpha,\delta_{g_j}\psi\rangle_{L^2},\qquad \psi\in\Omega^{p+1}(\widetilde M_j).
\]
Applying the density result \cite[Lemma~5.3]{Licht2019}
in coordinate charts and using a partition of unity, we obtain $\alpha_k\in\Omega_c^p(\widetilde M_j^\circ)$ such that
\begin{equation*}
 \|\alpha_k-\alpha\|_{L^2(\widetilde M_j,g_j)}
 +\|d\alpha_k-\zeta\|_{L^2(\widetilde M_j,g_j)}
 \longrightarrow0.
\end{equation*}
If, in addition, $\delta_{g_j}\zeta=0$ distributionally,
Green's formula and passage to the limit give
\begin{equation}\label{erelativeenergyvanishing}
 \|\zeta\|_{L^2(\widetilde M_j,g_j)}^2
 =\lim_{k\to\infty} \langle d\alpha_k,\zeta\rangle_{L^2(\widetilde M_j,g_j)}=\lim_{k\to\infty}
   \langle\alpha_k,\delta_{g_j}\zeta \rangle_{L^2(\widetilde M_j,g_j)}
 =0.
\end{equation}
Thus $\zeta=0$. We will use this observation in the proof of
Proposition~\ref{czprrecoverexactkernel}.

Composing the principal symbols of $d$, $\delta_{g_j}$ and the
interior principal symbol
$\sigma_{-2}(G_{\mathrm{rel},j}^{(p+1)})=|\xi|^{-2}_{g_j}\Id$ gives
\begin{equation*}
 \sigma_0(\mathsf P_j)(x,\xi)
 =|\xi|_{g_j}^{-2}\eps_\xi\ \iota_{\xi^{\sharp_{g_j}}}.
\end{equation*}

We write the Schwartz kernel of $\mathsf P_j$ as a double
$(p+1)$-form:
\[
 \mathsf K_j(x,y)\in \Lambda^{p+1}T_x^*\widetilde M_j\otimes\Lambda^{p+1}T_y^*\widetilde M_j,
\]
in the distributional sense, with convention
\begin{equation*}
 (\mathsf P_j u)(x)=\int_{\widetilde M_j}
       \langle\mathsf K_j(x,y),u(y)\rangle_{g_j(y)}\,dV_{g_j}(y),
\end{equation*}
where only the second tensor factor of $\mathsf K_j$ is paired with $u(y)$.
Fix $y\in\widetilde M_j^\circ$ and
$\vartheta\in(\Lambda^{p+1}T_y^*\widetilde M_j)^*$. Define
 \begin{align}
 \Omega_{j,y,\vartheta}(x)
 &:=\delta_{g_j,x}G_{\mathrm{rel},j}^{(p+1)}(x,y)\vartheta,\notag \\
 \eta_{j,y,\vartheta}(x)
 &:=d\Omega_{j,y,\vartheta}(x)=\mathsf K_j(x,y)\vartheta.\label{eOmegaandKcolumn}
 \end{align}
Here $\vartheta$ acts on the second tensor factor, and all
differential operators act in $x$.

Since relative harmonic forms are coclosed, the definition of
$G_{\mathrm{rel},j}^{(p+1)}$ gives
\[
 \delta_{g_j,x}\Delta_{g_j,x}
 G_{\mathrm{rel},j}^{(p+1)}(x,y)=0,\qquad x\neq y.
\]
Thus,
 \begin{align}
 \delta_{g_j}\Omega_{j,y,\vartheta}&=0,\notag\\
 \Delta_{g_j}\Omega_{j,y,\vartheta}
 &=\delta_{g_j,x}\bigl( \Delta_{g_j,x}G_{\mathrm{rel},j}^{(p+1)}(x,y)\vartheta
   \bigr)=0.\label{coclosedKj}
 \end{align}
Consequently,
 \begin{align*}
 d\eta_{j,y,\vartheta}
 &=0,\\
 \delta_{g_j}\eta_{j,y,\vartheta} &=\Delta_{g_j}\Omega_{j,y,\vartheta}-d\delta_{g_j}\Omega_{j,y,\vartheta}=0.
 \end{align*}
Finally, the relative condition
$\tang(\delta_{g_j,x}G_{\mathrm{rel},j}^{(p+1)})=0$ implies
$\tang\Omega_{j,y,\vartheta}=0$. Since tangential trace
commutes with $d$,
\[
 \tang\eta_{j,y,\vartheta}=d_{\partial\widetilde M_j}
      (\tang\Omega_{j,y,\vartheta})=0\qquad\text{on }\partial\widetilde M_j.
\]

Since $\eta_{j,y,\vartheta}(x)=\mathsf K_j(x,y)\vartheta$
and $\vartheta$ is arbitrary,
it follows that
\[
 \Delta_{g_j,x}\mathsf K_j(x,y)=0,\qquad x\neq y.
\]
Self-adjointness of $\mathsf P_j$ makes $\mathsf K_j$ symmetric
under exchange of its tensor factors, so also
$\Delta_{g_j,y}\mathsf K_j(x,y)=0$ for $x\neq y$.
Elliptic boundary regularity ensures smoothness in both variables
up to $\partial\widetilde M_j$ off the diagonal. On an interior
coordinate patch where $g_j$ is analytic, analytic elliptic
regularity implies analyticity in the variable lying in that patch.
If the metric and a boundary patch are analytic up to the
boundary, the same conclusion holds up to that patch by analytic
boundary regularity
\cite[Section~6.6]{Morrey1966}.
Since $\overline{U_0}$ and $X_j$ are disjoint compact sets,
$\mathsf K_j$ and all its derivatives are uniformly bounded
on $\overline{U_0}\times X_j$.

The next lemma allows us to change the metric conformally in the
critical case without changing the projection or its kernel.

\begin{lemma}

Let $g$ be a smooth metric on a compact oriented smooth manifold
$\widetilde M$ with smooth boundary. If $n=2p+2$, then, for every
$f\in C^\infty(\widetilde M)$,
\begin{equation}\label{eexactprojectionconformal}
 \mathsf P_{e^{2f}g}=\mathsf P_g.
\end{equation}
The double-form kernels also coincide. Moreover, on $(p+1)$-forms,
 the equations
$d\eta=\delta_g\eta=0$ are conformally invariant.
\end{lemma}

\begin{proof}
The subspace of relative exact forms is independent of the metric.
Since $n=2p+2$, the fibre pairing on $(p+1)$-forms satisfies
\begin{equation}\label{inproductconformal}
 \langle\cdot,\cdot\rangle_{e^{2f}g}\,dV_{e^{2f}g}
 =e^{(n-2p-2)f}\langle\cdot,\cdot\rangle_g\,dV_g
 =\langle\cdot,\cdot\rangle_g\,dV_g.
\end{equation}
Thus $g$ and $e^{2f}g$ induce the same $L^2$ inner product.
Their orthogonal projections onto the metric-independent
subspace of relative exact forms therefore coincide:
$\mathsf P_{e^{2f}g}=\mathsf P_g$. Their double-form kernels also coincide. The last assertion follows from \eqref{confomaldelta}.
\end{proof}

We next compare the kernels on the common exterior region using
only boundary data supported in the observed patch.

\begin{proposition}
\label{czprrecoverexactkernel}
 Assume \eqref{equDNintropar} and let $F_0$ be the given
collar diffeomorphism satisfying \eqref{ecollaralternatives}
or \eqref{ecollarconformalalternative}. Construct $\widetilde M_j$ and $U$ as above. After identifying the
exterior collars by $F_0$ in both kernel variables, we have
\begin{equation}\label{ekernelequalitycollar}
 \mathsf K_1(x,y)=\mathsf K_2(x,y),
 \qquad x,y\in U,\quad x\neq y.
\end{equation}
Moreover, for every $y\in U$ and
$\vartheta\in(\Lambda^{p+1}T_y^*U)^*$,
$F_0^*\eta_{2,y,\vartheta}=\eta_{1,y,\vartheta}$
in an interior collar near each point of $\Gamma_1$.
\end{proposition}

\begin{proof}
When $n=2p+2$, take $f_0$ from
\eqref{ecollarconformalalternative} and choose
$\chi_0\in C^\infty(\widetilde M_1)$ equal to one near
$\overline U$, with support contained in the domain of $f_0$.
Extend $\chi_0f_0$ by zero to $\widetilde M_1$ and replace
$g_1$ by $e^{2\chi_0f_0}g_1$.
Since $\chi_0f_0|_{\partial M_1}=0$,
\eqref{criticconformalinvariance} preserves
$\Lambda_{g_1}$ on the original boundary $\partial M_1$,
while \eqref{eexactprojectionconformal} preserves
$\mathsf K_1$ on $\widetilde M_1$.
Thus $F_0^*g_2=g_1$ near $\overline U$, directly when
$n\neq2p+2$ and after this replacement when $n=2p+2$.

Fix $y\in U$ and $\vartheta\in(\Lambda^{p+1}T_y^*U)^*$.
Since $y\notin M_1$, the form $\Omega_{1,y,\vartheta}$ is harmonic
and coclosed on $M_1$. Its boundary value satisfies
\[
 \varphi:=\tang\Omega_{1,y,\vartheta}|_{\partial M_1},\qquad \supp\varphi\subset\supp\chi
 \Subset\Gamma.
\]
Here $\chi\in C_c^\infty(\Gamma_0)$ is the function defining
the exterior extension in \eqref{elocalexteriorextension}.
Let $\widehat\Omega_2$ solve \eqref{equationrelativebvp} on $M_2$
with trace $\varphi$, using $\kappa$ to identify the boundary data.
It is coclosed by Proposition~\ref{lemmadeltazero}.
For $\widehat\eta_2=d\widehat\Omega_2$, we have
\begin{equation}\label{efulltracematcheta}
 \tang\widehat\eta_2=\tang\eta_{1,y,\vartheta},\qquad
 \normal\widehat\eta_2=\Lambda_{g_2}\varphi =\Lambda_{g_1}\varphi=\normal\eta_{1,y,\vartheta}
 \qquad\text{on }\Gamma.
\end{equation}
The middle equality uses the identity
\eqref{equDNintropar}.

On $\widetilde M_2$, define
\[
 \zeta:=
 \begin{cases}
 \widehat\eta_2-\eta_{2,y,\vartheta},  &\text{on }M_2,\\
 \eta_{1,y,\vartheta}-\eta_{2,y,\vartheta}, &\text{on }\widetilde M_2\setminus M_2.
\end{cases}\qquad
 \alpha:=
 \begin{cases}
  \widehat\Omega_2-\Omega_{2,y,\vartheta},&\text{on }M_2,\\
\Omega_{1,y,\vartheta}-\Omega_{2,y,\vartheta},&\text{on }\widetilde M_2\setminus M_2.
 \end{cases}
\]
Since $g_1=g_2$ near $y$ after identification by $F_0$,
subtracting the identities \eqref{defGrelj} for
$\Delta G_{\mathrm{rel},1}^{(p+1)}$ and
$\Delta G_{\mathrm{rel},2}^{(p+1)}$ cancels the Dirac delta terms in their kernels
at $x=y$. The remaining harmonic projection kernels are smooth.
Hence
$G_{\mathrm{rel},1}^{(p+1)}(x,y)
-G_{\mathrm{rel},2}^{(p+1)}(x,y)$ is smooth there.
By construction, $\alpha$ and $\zeta$ are therefore smooth
near $y$. They are piecewise smooth elsewhere, so
$\alpha,\zeta\in L^2(\widetilde M_2)$.

Both pieces of $\zeta$ are coclosed
for $g_2$. For every
$\psi\in\Omega_c^p(\widetilde M_2^\circ)$, applying Green's
formula on the two pieces gives
\[
 \langle\delta_{g_2}\zeta,\psi\rangle_{L^2(\widetilde M_2,g_2)}=\langle\zeta,d\psi\rangle_{L^2(\widetilde M_2,g_2)}
 =\int_{\{r=0,\chi>0\}}\left\langle
   \normal(\widehat\eta_2-\eta_{1,y,\vartheta}),\tang\psi \right\rangle\,dS=0.
\]
The last equality follows from \eqref{efulltracematcheta}.
Thus $\delta_{g_2}\zeta=0$ in the distributional sense.
The tangential traces of the two pieces of $\alpha$ agree
across the interface, since
$\tang\widehat\Omega_2=\varphi
=\tang\Omega_{1,y,\vartheta}$.
Their tangential traces vanish on the exterior boundary
$\partial\widetilde M_2$ by the relative boundary conditions
and $\varphi=0$ on
$\partial M_2\cap\partial\widetilde M_2$.
Applying Green's formula on the two pieces, the interface
terms cancel and the exterior boundary terms vanish. Hence
\[
 \langle\zeta,\psi\rangle_{L^2(\widetilde M_2,g_2)}
 =\langle\alpha,\delta_{g_2}\psi\rangle_{L^2(\widetilde M_2,g_2)},
 \qquad \psi\in\Omega^{p+1}(\widetilde M_2).
\]
Thus $d\alpha=\zeta$ in the distributional sense, and
$\alpha$ has vanishing tangential trace in the weak sense.

Equation~\eqref{erelativeenergyvanishing} now applies to
$d\alpha=\zeta$ and $\delta_{g_2}\zeta=0$.
Hence $\zeta=0$, giving
$\widehat\eta_2=\eta_{2,y,\vartheta}$ on $M_2$ and
$\eta_{1,y,\vartheta}=\eta_{2,y,\vartheta}$ on $U$.
Since $y$ and $\vartheta$ were arbitrary,
\eqref{ekernelequalitycollar} follows.

We have proved \eqref{ekernelequalitycollar} for
$x,y\in U$, $x\neq y$. We next establish the corresponding
identity for $x$ in an interior collar of $\Gamma_1$,
with $y\in U$ fixed.

Return to the original metrics $g_1,g_2$. When $n=2p+2$,
\eqref{eexactprojectionconformal} preserves
$\mathsf K_j$ and hence $\eta_{j,y,\vartheta}$.
Near any point of
$\Gamma_1$, the form
$\eta_{1,y,\vartheta}-F_0^*\eta_{2,y,\vartheta}$ is closed and
coclosed for $g_1$, and both its tangential and normal traces
vanish by \eqref{efulltracematcheta}.
Its zero extension $\omega_0$ across $\Gamma_1$ therefore satisfies
$d\omega_0=0$ and $\delta_{g_1}\omega_0=0$ weakly. Elliptic regularity and unique
continuation \cite[Theorem~1]{AronszajnKrzywickiSzarski1962} give $\eta_{1,y,\vartheta}=F_0^*\eta_{2,y,\vartheta}$ in an interior collar of $\Gamma_1$.
\end{proof}

\subsection{Interior embedding into the space of Hilbert--Schmidt operators}
Let $\mathscr H:=L^2(U_0;\Lambda^{p+1}T^*U_0)$.
When $n\neq2p+2$, we have $g_2=g_1$ on $U_0$.
When $n=2p+2$, we have $g_2=e^{2f_0}g_1$ on $U_0$.
By \eqref{inproductconformal}, in both cases, $g_1$ and $g_2$ induce the same inner product on the Hilbert space $\mathscr H$.

For $x\in\widetilde M_j\setminus\overline{U_0}$, define
the linear map
 \begin{align*}
 \mathcal{J}_j(x):(\Lambda^{p+1}T_x^*\widetilde M_j)^* &\longrightarrow \mathscr H,\\
 \vartheta&\longmapsto
 \left.\mathsf K_j(\cdot,x)\vartheta\right|_{U_0}.
 \end{align*}
Thus $(\mathcal{J}_j(x)\vartheta)(z)=\mathsf K_j(z,x)\vartheta$
for $z\in U_0$. Using the inner product on $\mathscr H$ and the
dual fibre metric induced by $g_j$, let
\[
 \mathcal{J}_j(x)^*:\mathscr H\longrightarrow (\Lambda^{p+1}T_x^*\widetilde M_j)^*
\]
denote its adjoint. The composition $\mathcal{J}_j(x)\mathcal{J}_j(x)^*$
is a nonnegative finite-rank operator on $\mathscr H$.
Whenever $\mathcal{J}_j(x)\neq0$, define
\begin{equation}\label{normalizedkernelgram}
 \mathcal T_j(x) :=\frac{\mathcal{J}_j(x)\mathcal{J}_j(x)^*}{\Tr(\mathcal{J}_j(x)\mathcal{J}_j(x)^*)}.
\end{equation}
Lemmas~\ref{lemmaexactkernelembedding} and
\ref{lemmaexactkernelboundary} prove that $\mathcal{J}_j(x)\neq0$ on $X_j$.
Wherever $\mathcal{J}_j(x)\neq0$,
\[
 \Tr(\mathcal{J}_j\mathcal{J}_j^*)=\Tr(\mathcal{J}_j^*\mathcal{J}_j)>0,\quad \ran \mathcal T_j=\ran \mathcal{J}_j\mathcal{J}_j^*=\ran \mathcal{J}_j, \quad \Tr\mathcal T_j=1.
\]
Thus $\mathcal T_j$ takes values in the real Hilbert space of
self-adjoint Hilbert--Schmidt operators on $\mathscr H$. Each value is
nonnegative and has rank at most $\binom{n}{p+1}$.

Off-diagonal regularity makes $\mathcal{J}_j(x)$ smooth in operator norm
in local coordinates. Since its domain is finite-dimensional,
$\mathcal{J}_j(x)\mathcal{J}_j(x)^*$ is smooth in the Hilbert--Schmidt norm and
its trace is smooth. Thus $\mathcal T_j$ is smooth wherever
$\mathcal{J}_j\neq0$.
If $g_j$ is analytic on an interior coordinate patch, analytic
elliptic regularity for $\Delta_{g_j,x}\mathsf K_j(z,x)=0$ gives,
on every smaller patch,
\begin{equation}\label{ekernelanalyticestimate}
 \sup_{z\in\overline{U_0}}|\partial_x^\alpha\mathsf K_j(z,x)|\leq C^{|\alpha|+1}|\alpha|!.
\end{equation}
The constant is uniform in $x$ on the smaller patch, since
$\overline{U_0}$ is compact and disjoint from it. Integrating
in $z$ yields the same estimate for $\partial_x^\alpha \mathcal{J}_j(x)$
in operator norm. The metric coefficients in $x$ are analytic,
so $\mathcal T_j$ is analytic in the Hilbert--Schmidt norm
wherever $\mathcal{J}_j\neq0$. If the metric and a boundary patch are
analytic up to the boundary, the same conclusion holds there
by analytic boundary regularity.

We next examine the singularities of the kernel.
For $\vartheta\in (\Lambda^{p+1}T_x^*\widetilde M_j)^*$, write $\vartheta^\flat$ for its
metric dual, so that
\begin{equation*}
 \vartheta(w)=\langle\vartheta^\flat,w\rangle_{g_j(x)}, \qquad w\in\Lambda^{p+1}T_x^*\widetilde M_j.
\end{equation*}
We denote the Dirac distribution at $x$ by $\boldsymbol{\delta}_x$.

\begin{lemma}
\label{lemmaexactkernelsingularities}
Fix $j$, an interior point $x\in\widetilde M_j^\circ\setminus\overline{U_0}$ and a nonzero
$\vartheta\in(\Lambda^{p+1}T_x^*\widetilde M_j)^*$.
\begin{enumerate}[label=(\roman*)]
\item The restriction of $\mathsf K_j(\cdot,x)\vartheta$ to
$\widetilde M_j^\circ\setminus\{x\}$ has no smooth extension
across $x$.
\item Extend $\vartheta$ with constant coefficients in fixed
local coordinates near $x$. If
$V\in T_x\widetilde M_j$ and
$\vartheta_1\in(\Lambda^{p+1}T_x^*\widetilde M_j)^*$ satisfy
\begin{equation*}
 \supp\bigl(D_x[\mathsf K_j(z,x)\vartheta][V]
       -\mathsf K_j(z,x)\vartheta_1\bigr)\subset\{x\}
\end{equation*}
as distributions in $z$, then $V=0$.
\end{enumerate}
\end{lemma}

\begin{proof}
Suppress $j$ and choose fixed normal coordinates
$z=(z^1,\ldots,z^n)$ centered at $x$.
Since relative harmonic forms are coclosed by Proposition~\ref{lemmadeltazero}, the definition of
$\mathsf P$ gives the distributional identities
\begin{equation}\label{ekerneldistributionidentities}
 d_z[\mathsf K(z,x)\vartheta]=0,\qquad
 \Delta_z[\mathsf K(z,x)\vartheta]
 =d_z\delta_z(\vartheta^\flat\boldsymbol{\delta}_x).
\end{equation}
For (i), suppose that the restriction of $\mathsf K(\cdot,x)\vartheta$
to $z\neq x$ admits a smooth extension across $x$. Subtracting
the extension leaves a distribution $S$ supported at $x$. By \cite[Theorem~2.3.4]{HormanderI}, $S$ is a finite sum of
derivatives of $\boldsymbol{\delta}_x$.
The smooth extension is closed and harmonic, so
\eqref{ekerneldistributionidentities} gives
\begin{equation*}
 dS=0,\qquad\Delta S=d\delta(\vartheta^\flat\boldsymbol{\delta}_x).
\end{equation*}
The right-hand side contains at most second derivatives of
$\boldsymbol{\delta}_x$. Comparing the highest derivatives,
using $\sigma_2(\Delta)=|\xi|_{g(x)}^2\Id$, therefore shows that
$S=\eta_0\boldsymbol{\delta}_x$ for some
$\eta_0\in\Lambda^{p+1}T_x^*\widetilde M$.
Then $dS=0$ implies $\eps_\xi\eta_0=0$ for every $\xi$,
and hence $\eta_0=0$ because $p+1<n$.
By the second identity,
$\eps_\xi\iota_{\xi^{\sharp_{g(x)}}}\vartheta^\flat=0$
for every $\xi$.
Pairing with $\vartheta^\flat$ yields
$\big|\iota_{\xi^{\sharp_{g(x)}}}\vartheta^\flat\big|^2_{g(x)}=0$,
which forces $\vartheta^\flat=0$ since $p+1\geq1$.
This contradiction proves that $\mathsf K(\cdot,x)\vartheta$
has a nonremovable singularity at $x$.

For (ii), by hypothesis, the distribution in $z$
\begin{equation*}
 R:=D_x[\mathsf K(z,x)\vartheta][V]  -\mathsf K(z,x)\vartheta_1
\end{equation*}
is supported at $x$. Differentiating \eqref{ekerneldistributionidentities} gives
\begin{equation}\label{dridenformula}
 dR=0,\qquad\Delta R=d\delta\bigl(D_x[\vartheta^\flat\boldsymbol{\delta}_x][V]
       -\vartheta_1^\flat\boldsymbol{\delta}_x\bigr).
\end{equation}
The right-hand side contains at most third derivatives of
$\boldsymbol{\delta}_x$, so the same comparison shows that
$R$ contains at most first derivatives. Its Fourier transform
in $z$ is therefore a polynomial of degree at most one.
Let $\widehat R_1(\xi)$ be its degree-one homogeneous part.
In the fixed local coordinates $(z^1,z^2,\dots,z^n)$ centered at $x$,
\begin{equation*}
 D_x[\vartheta^\flat\boldsymbol{\delta}_x][V]=-\vartheta^\flat\sum_{k=1}^n
       V^k\partial_{z^k}\boldsymbol{\delta}_x +D_x[\vartheta^\flat][V]\,\boldsymbol{\delta}_x.
\end{equation*}
Taking the Fourier transform in $z$, the right-hand side of the equality becomes
\[
-\sqrt{-1}\left(\sum_{k=1}^{n}V^k\xi_k\right)\vartheta^\flat+D_x[\vartheta^\flat][V].
\]
Thus its degree-one term is
$-\sqrt{-1}\xi(V)\vartheta^\flat$. Comparing the highest
degrees in \eqref{dridenformula} gives
\begin{equation*}
 \eps_\xi\widehat R_1(\xi)=0,\qquad|\xi|_{g(x)}^2\widehat R_1(\xi)=-\sqrt{-1}\xi(V) \eps_\xi\iota_{\xi^{\sharp_{g(x)}}}\vartheta^\flat.
\end{equation*}
If $V\neq0$, the second identity makes the linear polynomial
$\widehat R_1$ vanish on $\xi(V)=0$. Therefore
$\widehat R_1(\xi)=\xi(V)\eta_1$ for a constant form $\eta_1$.
The first identity implies $\eps_\xi\eta_1=0$ for every $\xi$
by continuity, so $\eta_1=0$ because $p+1<n$.
The second identity then forces
$\eps_\xi\iota_{\xi^{\sharp_{g(x)}}}\vartheta^\flat=0$
for every $\xi$ by continuity. As in (i), this contradicts
$\vartheta\neq0$. Hence $V=0$.
\end{proof}

\begin{lemma}
\label{lemmaexactkernelembedding}
For every $x\in\widetilde M_j^\circ\setminus\overline{U_0}$,
$\mathcal{J}_j(x)$ is injective and has rank $\binom{n}{p+1}$.
The map $\mathcal T_j$ is a smooth injective immersion
from $\widetilde M_j^\circ\setminus\overline{U_0}$ into the
real Hilbert space of self-adjoint Hilbert--Schmidt operators
on $\mathscr H$. It is analytic in the Hilbert--Schmidt norm
on every interior coordinate patch where $g_j$ is analytic.
\end{lemma}

\begin{proof}
Fix $j$ and suppress it throughout the proof.
If $\mathcal{J}(x)\vartheta=0$, then
$\mathsf K(z,x)\vartheta=0$ for $z\in U_0$.
This form is closed and coclosed away from $x$ by \eqref{eOmegaandKcolumn} and \eqref{coclosedKj}.
Since $n\geq2$, the punctured interior is connected, and
smooth unique continuation
\cite[Theorem~1]{AronszajnKrzywickiSzarski1962} gives
\begin{equation*}
 \mathsf K(z,x)\vartheta=0,\qquad z\in\widetilde M^\circ\setminus\{x\}.
\end{equation*}
Lemma~\ref{lemmaexactkernelsingularities}(i) forces
$\vartheta=0$. Thus $\mathcal{J}(x)$ is injective and has rank
$\dim\bigl(\Lambda^{p+1}T_x^*\widetilde M\bigr)=\binom{n}{p+1}$. Since $\Tr(\mathcal{J}^*\mathcal{J})>0$, the regularity of $\mathcal{J}$ ensures the
asserted smoothness and analyticity of $\mathcal T$.

Suppose $\mathcal T(x)=\mathcal T(x')$ with $x\neq x'$.
Their ranges agree, so for any nonzero $\vartheta$ there is
$\vartheta'$ such that
\begin{equation*}
 \mathsf K(z,x)\vartheta=\mathsf K(z,x')\vartheta',\qquad z\in U_0.
\end{equation*}
Both sides are closed and coclosed on
$\widetilde M^\circ\setminus\{x,x'\}$. The same smooth unique continuation extends the
identity to
$z\in\widetilde M^\circ\setminus\{x,x'\}$.
The right side is smooth near $x$, whereas the left side
has a nonremovable singularity there by
Lemma~\ref{lemmaexactkernelsingularities}(i).
This contradiction proves injectivity of $\mathcal T$.

For immersion, suppose $D\mathcal T(x)[V]=0$. Here $D\mathcal T(x)$ denotes the Fr\'echet differential
of $\mathcal T$ at $x$, taken in the Hilbert--Schmidt norm.
By \eqref{normalizedkernelgram},
\begin{equation}\label{egramrangefactorization}
 \mathcal{J}=\mathcal T\,\Tr(\mathcal{J}^*\mathcal{J})\,\mathcal{J}(\mathcal{J}^*\mathcal{J})^{-1}.
\end{equation}
Differentiating along $V$ at $x$ and using
$D\mathcal T(x)[V]=0$ gives
\begin{equation*}
 D\mathcal{J}(x)[V] =\mathcal T(x)\,D\!\left[\Tr(\mathcal{J}^*\mathcal{J})\,\mathcal{J}(\mathcal{J}^*\mathcal{J})^{-1}\right](x)[V].
\end{equation*}
Thus $\ran(D\mathcal{J}(x)[V])\subset\ran \mathcal{J}(x)$.
Choose a nonzero $\vartheta$ and extend it with constant
coefficients near $x$. There is $\vartheta_1$ such that
\begin{equation*}
 D_x[\mathsf K(z,x)\vartheta][V]=\mathsf K(z,x)\vartheta_1,\qquad z\in U_0.
\end{equation*}
The operators $d_z$ and $\delta_{g,z}$ commute with $D_x$, so both sides are closed and coclosed in $z$
away from $x$.
Unique continuation extends this identity to $\widetilde M^\circ\setminus\{x\}$.
Lemma~\ref{lemmaexactkernelsingularities}(ii) therefore gives
$V=0$, proving that $\mathcal T$ is an immersion.
\end{proof}

\subsection{Boundary rank drop and spectral boundary coordinates}

We next examine the boundary, where the rank decreases.
The small eigenvalues will provide a normal coordinate for
extending the map $F=\mathcal T_2^{-1}\circ\mathcal T_1$ to
the unchanged part of the original boundary $\partial M_1\setminus\supp\chi$.

Suppress $j$ and write $x=(r,y)$ in boundary normal coordinates,
with $r\geq0$ inside $\widetilde M$.
As in Section~\ref{sectiontangeormal}, use the collar map
to identify tangential input forms on $Y_r$ with forms on $Y$,
regarding $r$ as a parameter.
The metric dual below is
taken at $(r,y)$, and all fibre adjoints use the metric there.
For $\alpha\in\Lambda^{p+1}T_y^*\partial\widetilde M$ and
$\beta\in\Lambda^pT_y^*\partial\widetilde M$, write
$\vartheta^\flat=\alpha+dr\wedge\beta$ and define $\mathcal{J}_T,\mathcal{J}_N$ by
\[
 \mathcal{J}(r,y)\vartheta=\mathcal{J}_T(r,y)\alpha+\mathcal{J}_N(r,y)\beta.
\]
The tangential trace of $\mathsf K$ in its first variable
vanishes by the relative boundary condition. The same holds in
its second variable by symmetry. Hence $\mathcal{J}_T(0,y)=0$.

For $\varphi\in\Omega^p(\partial\widetilde M)$, let
$\omega_\varphi$ solve \eqref{eoriginalproblem} on
$\widetilde M$. The form $d\omega_\varphi$ is independent
of the harmonic ambiguity in $\omega_\varphi$.
For distributional $\varphi$, extend $B$ in \eqref{eBformula}
to distributions and apply the full Dirichlet Poisson operator
for \eqref{eulldirichlet} to $\varphi+dr\wedge B\varphi$.
This defines $\omega_\varphi$ and $d\omega_\varphi$, both smooth
in $\widetilde M^\circ$.
Therefore we have $\tang\omega_\varphi=\varphi$
in the sense of distributions. Since tangential trace
commutes with $d$, we obtain
\begin{equation}\label{eboundaryharmonictrace}
 \tang(d\omega_\varphi)=d_{\partial\widetilde M}\varphi
\end{equation}
distributionally.
\begin{lemma}
\label{lemmaboundarypoissonrepresentation}
At every $y\in\partial\widetilde M$, we have
\begin{equation}\label{eboundarykernelpoisson}
  \mathcal{J}_N(0,y)\beta =\left.d\omega_{\beta\boldsymbol{\delta}_y}\right|_{U_0},\qquad
  \partial_r\mathcal{J}_T(0,y)\alpha =\left.d\omega_{\delta_{\partial\widetilde M}(\alpha\boldsymbol{\delta}_y)}\right|_{U_0}.
\end{equation}
Moreover, for any distributional boundary $p$-form $\varphi$,
\begin{equation}\label{eboundarypoissondetection}
 \left.d\omega_\varphi\right|_{U_0}=0\quad\Longrightarrow\quad
 d_{\partial\widetilde M}\varphi=0.
\end{equation}
Here $\boldsymbol{\delta}_y$ denotes the Dirac distribution
at $y$ on
$\partial\widetilde M$.
\end{lemma}

\begin{proof}
For smooth $\varphi$ and
$u\in C_c^\infty(\widetilde M^\circ;
\Lambda^{p+1}T^*\widetilde M)$, the definition of $\mathsf P$
and Green's formula give
\[
  \langle\delta(u-\mathsf Pu),\alpha\rangle_{L^2(\widetilde M)}=\langle u-\mathsf Pu,d\alpha\rangle_{L^2(\widetilde M)}=0,\qquad \langle d\omega_\varphi,d\alpha\rangle_{L^2(\widetilde M)}=\langle\delta d\omega_\varphi,\alpha   \rangle_{L^2(\widetilde M)}=0,
\]
for every smooth $p$-form $\alpha$ with $\tang\alpha=0$.
Here the first identity uses $d\alpha\in\ran\mathsf P$, and
the second uses
$\delta d\omega_\varphi
=\Delta\omega_\varphi-d\delta\omega_\varphi=0$.
Thus we have
\[\delta\mathsf Pu=\delta u,\qquad
d\omega_\varphi\perp\ran\mathsf P.\]
Applying Green's formula once more gives
\begin{equation}\label{eboundarypoissonpairing}
  \langle d\omega_\varphi,u\rangle_{L^2(\widetilde M)}=\langle\omega_\varphi,\delta\mathsf Pu\rangle_{L^2(\widetilde M)} =-\langle\varphi,\normal(\mathsf Pu)  \rangle_{L^2(\partial\widetilde M)}.
\end{equation}

For $u\in C_c^\infty(U_0;\Lambda^{p+1}T^*\widetilde M)$,
symmetry of $\mathsf K$, the identity $\nu_g=-\partial_r$,
and \eqref{eboundarypoissonpairing} give
\[
  \langle \mathcal{J}_N(0,y)\beta,u\rangle_{L^2(U_0)}=-\langle\beta,\normal(\mathsf Pu)(y)\rangle_{g(y)}=\langle d\omega_{\beta\boldsymbol{\delta}_y},u\rangle_{L^2(U_0)}.
\]
Taking the coefficient of $dr$ in
$d(\mathsf Pu)=0$, using \eqref{ecollard} in this coordinate
coframe, gives
\[
 \tang(\partial_r\mathsf Pu)
 =-d_{\partial\widetilde M}\bigl(\normal(\mathsf Pu)\bigr).
\]
Keeping the coordinate coefficients of $\alpha$ fixed, by symmetry of $\mathsf K$,
\begin{equation*}
 \langle \mathcal{J}_T(r,y)\alpha,u\rangle_{L^2(U_0)}=\langle\alpha,(\mathsf Pu)(r,y)\rangle_{g(r,y)}.
\end{equation*}
We differentiate this identity in $r$ at $r=0$.
Since $\tang(\mathsf Pu)=0$, the terms obtained by
differentiating the metric vanish.
Integration by parts on the boundary and
\eqref{eboundarypoissonpairing} then give
 \begin{align*}
  \langle\partial_r\mathcal{J}_T(0,y)\alpha,u\rangle_{L^2(U_0)} &=\langle\alpha,\tang(\partial_r\mathsf Pu)(y)\rangle_{g(y)}\\
  &=-\langle\alpha,
       d_{\partial\widetilde M}\bigl(\normal(\mathsf Pu)\bigr)(y) \rangle_{g(y)}\\
  &=-\big\langle\delta_{\partial\widetilde M}(\alpha\boldsymbol{\delta}_y),
       \normal(\mathsf Pu) \big\rangle_{\partial\widetilde M}\\
  &=\big\langle d\omega_{\delta_{\partial\widetilde M}(\alpha\boldsymbol{\delta}_y)},u  \big\rangle_{L^2(U_0)}.
 \end{align*}
Here $\langle\cdot,\cdot\rangle_{\partial\widetilde M}$ denotes
the distributional pairing on the boundary.
Since $u$ is arbitrary, this proves
\eqref{eboundarykernelpoisson}.

Finally, $d\omega_\varphi$ is closed and coclosed in the
connected interior $\widetilde M^\circ$. If it vanishes on $U_0$, unique
continuation makes it vanish throughout the interior.
Taking its distributional trace in
\eqref{eboundaryharmonictrace} proves
\eqref{eboundarypoissondetection}.
\end{proof}

\begin{lemma}
\label{lemmaexactkernelboundary}
For each $j$ and every $y\in\partial\widetilde M_j$,
\begin{equation}\label{eexactkernelboundaryrank}
 \rank \mathcal{J}_j(y)=\rank\mathcal T_j(y)=\binom{n-1}{p}<\binom{n}{p+1}.
\end{equation}
The maps $\mathcal T_j$ restrict to topological embeddings on
$X_j$, and their restrictions to $\partial\widetilde M_j$ are
smooth injective immersions. In the boundary notation above,
\begin{equation}\label{eboundarycolumntransversality}
 \partial_r\mathcal{J}_T(0,y)\alpha\in\ran \mathcal{J}_N(0,y)\quad\Longrightarrow\quad\alpha=0.
\end{equation}
\end{lemma}

\begin{proof}
Fix $j$ and suppress it.
If $\mathcal{J}_N(0,y)\beta=0$, Lemma~
\ref{lemmaboundarypoissonrepresentation} gives
$d_{\partial\widetilde M}(\beta\boldsymbol{\delta}_y)=0$.
Comparing the coefficients of the first derivatives of
$\boldsymbol{\delta}_y$ yields $\eps_\xi\beta=0$ for every
$\xi\in T_y^*\partial\widetilde M$.
Since $p\leq n-2$, this forces $\beta=0$.
Thus $\mathcal{J}_N(0,y)$ is injective. Together with $\mathcal{J}_T(0,y)=0$,
this proves \eqref{eexactkernelboundaryrank}.

If $\partial_r\mathcal{J}_T(0,y)\alpha=\mathcal{J}_N(0,y)\beta$, Lemma~
\ref{lemmaboundarypoissonrepresentation}
implies
\begin{equation*}
 d_{\partial\widetilde M}\delta_{\partial\widetilde M}(\alpha\boldsymbol{\delta}_y)=d_{\partial\widetilde M}(\beta\boldsymbol{\delta}_y).
\end{equation*}
The second derivatives of $\boldsymbol{\delta}_y$ on the left
must vanish, since the right-hand side has order at most one.
Hence $\eps_\xi\iota_{\xi^\sharp}\alpha=0$ for every $\xi$.
Pairing with $\alpha$ gives
$\big|\iota_{\xi^\sharp}\alpha\big|_{g(y)}^2=0$, so $\alpha=0$.
This proves \eqref{eboundarycolumntransversality}.

Off-diagonal boundary regularity and $\mathcal{J}_N(0,y)\neq0$ make
$\mathcal T$ smooth up to $\partial\widetilde M$ in the
Hilbert--Schmidt norm.
If $\mathcal T(0,y)=\mathcal T(0,y')$ with $y\neq y'$, their
ranges agree. For $\beta\neq0$, there is $\beta'$ such that
$\mathcal{J}_N(0,y)\beta=\mathcal{J}_N(0,y')\beta'$.
Lemma~\ref{lemmaboundarypoissonrepresentation} then gives
\begin{equation*}
 d_{\partial\widetilde M}(\beta\boldsymbol{\delta}_y)=d_{\partial\widetilde M}(\beta'\boldsymbol{\delta}_{y'}).
\end{equation*}
The two sides have disjoint supports, so both vanish, contradicting
$\beta\neq0$. Thus $\mathcal T$ is injective on the boundary.
Lemma~\ref{lemmaexactkernelembedding} and the different
interior and boundary ranks show that $\mathcal T$ is continuous
and injective on $X$. Since $X$ is compact and the target Hilbert space is Hausdorff,
$\mathcal T|_X$ is a homeomorphism onto its image and hence a
topological embedding.

To prove boundary immersion, suppose
$D_y\mathcal T(0,y)[V]=0$ for
$V\in T_y\partial\widetilde M$.
Differentiating \eqref{egramrangefactorization} with
$\mathcal{J}_N$ in place of $\mathcal{J}$ gives
\begin{equation*}
 D_y[\mathcal{J}_N(0,y)\beta][V]\in\ran \mathcal{J}_N(0,y).
\end{equation*}
Here $\beta$ is arbitrary and is extended with constant
coefficients in boundary coordinates.
Choose $\beta_1$ with
$D_y[\mathcal{J}_N(0,y)\beta][V]=\mathcal{J}_N(0,y)\beta_1$.
Since $\varphi\mapsto d\omega_\varphi$ is continuous and linear,
differentiating the first identity in
\eqref{eboundarykernelpoisson} gives
\[
 D_y[\mathcal{J}_N(0,y)\beta][V]=\left.d\omega_{D_y[\beta\boldsymbol{\delta}_y][V]}\right|_{U_0}.
\]
Applying the same identity to $\beta_1$ gives
\[
 \mathcal{J}_N(0,y)\beta_1=\left.d\omega_{\beta_1\boldsymbol{\delta}_y}\right|_{U_0}.
\]
Subtracting these equalities, we obtain
\[
 \left.d\omega_{D_y[\beta\boldsymbol{\delta}_y][V]-\beta_1\boldsymbol{\delta}_y}\right|_{U_0}=0.
\]
Thus \eqref{eboundarypoissondetection} yields
\[
 d_{\partial\widetilde M} \bigl(D_y[\beta\boldsymbol{\delta}_y][V]-\beta_1\boldsymbol{\delta}_y\bigr)=0.
\]
Comparing the second derivatives of
$\boldsymbol{\delta}_y$ yields
$\xi(V)\eps_\xi\beta=0$ for every $\xi$ and $\beta$.
If $V\neq0$, one can choose $\xi,\beta$ with
$\xi(V)\eps_\xi\beta\neq0$, since $p\leq n-2$.
Thus $V=0$, proving boundary immersion.
\end{proof}

Wherever $\mathcal{J}_j(x)\neq0$, the operator $\mathcal T_j(x)$ is
positive semidefinite with trace one on $\mathscr H$, so
its eigenvalues $\lambda_k(\mathcal T_j(x))$ are nonnegative and sum to one.
Order them decreasingly, counting multiplicities and including
zeros, and define
\begin{equation}\label{egramboundarydefiningfunction}
 \rho_j(x):=\left(1-\sum_{k=1}^{\binom{n-1}{p}}\lambda_k(\mathcal T_j(x))
 \right)^{\frac12}.
\end{equation}

\begin{lemma}
\label{lemmagramboundarydefiningfunction}
Near every point of $\partial\widetilde M_j$, there is a
smooth positive function $c_j(r,y)$ such that
\begin{equation}\label{egramsmalleigenvalueexpansion}
 \rho_j(r,y)^2=r^2c_j(r,y) =c_j(0,y)r^2+O(r^3),\qquad c_j(0,y)>0.
\end{equation}
Consequently, $\rho_j$ has a smooth signed extension across
the boundary with $\rho_j(0,y)=0$ and
$\partial_r\rho_j(0,y)>0$. If the metric and boundary patch
are analytic, this extension is analytic.
\end{lemma}

\begin{proof}
Fix $j$ and suppress it.
By Lemma~\ref{lemmaexactkernelembedding} and
Lemma~\ref{lemmaexactkernelboundary}, $\mathcal{J}^*\mathcal{J}$ has
$\binom{n}{p+1}$ positive eigenvalues for $r>0$, whereas
exactly $\binom{n-1}{p}$ eigenvalues are positive at $r=0$.
With eigenvalues ordered decreasingly, continuity therefore
shows that the first $\binom{n-1}{p}$ eigenvalues remain
bounded away from zero, while the remaining
$\binom{n-1}{p+1}$ eigenvalues tend to zero as $r\to0$.

Since $\mathcal{J}\mathcal{J}^*$ and $\mathcal{J}^*\mathcal{J}$ have the same nonzero eigenvalues
and $\mathcal T=\mathcal{J}\mathcal{J}^*/\Tr(\mathcal{J}^*\mathcal{J})$, the definition of $\rho$
gives
\[
 \rho(r,y)^2 =\frac{\displaystyle
   \sum_{\ell=\binom{n-1}{p}+1}^{\binom{n}{p+1}}
   \lambda_\ell\bigl(\mathcal{J}(r,y)^*\mathcal{J}(r,y)\bigr)}
   {\Tr\bigl(\mathcal{J}(r,y)^*\mathcal{J}(r,y)\bigr)}.
\]
We therefore compute the second-order terms of these last
$\binom{n-1}{p+1}$ eigenvalues.

Let $\Pi_N(y)$ be the orthogonal projection in $\mathscr H$
onto $\ran \mathcal{J}_N(0,y)$.
With all blocks evaluated at $(r,y)$, we have
\[
 \mathcal{J}^*\mathcal{J} =\begin{pmatrix}\mathcal{J}_T^*\\ \mathcal{J}_N^*\end{pmatrix}
  \begin{pmatrix}\mathcal{J}_T&\mathcal{J}_N\end{pmatrix}
 =\begin{pmatrix}
   \mathcal{J}_T^*\mathcal{J}_T&\mathcal{J}_T^*\mathcal{J}_N\\
   \mathcal{J}_N^*\mathcal{J}_T&\mathcal{J}_N^*\mathcal{J}_N
  \end{pmatrix}.
\]
At $r=0$, $\mathcal{J}_T=0$ and $\mathcal{J}_N^*\mathcal{J}_N$ is positive definite.
Thus the normal block remains uniformly invertible for small
$r$, while the min--max principle and $\mathcal{J}_T=O(r)$ show that
the remaining eigenvalues are $O(r^2)$.

To determine these small eigenvalues, we eliminate the normal
component $\beta$ from the eigenvalue equation
$\mathcal{J}^*\mathcal{J}\binom{\alpha}{\beta}
=\lambda\binom{\alpha}{\beta}$, obtaining
\[
 \bigl[ \mathcal{J}_T^*\mathcal{J}_T-\mathcal{J}_T^*\mathcal{J}_N
 (\mathcal{J}_N^*\mathcal{J}_N-\lambda\Id)^{-1}\mathcal{J}_N^*\mathcal{J}_T \bigr]\alpha=\lambda\alpha.
\]
Since $\lambda=O(r^2)$ and $\mathcal{J}_T^*\mathcal{J}_N=O(r)$, replacing
$\lambda$ by zero inside the inverse changes this operator
by $O(r^4)$.
Hence the second-order terms of the small eigenvalues are
determined by
\[
 \mathcal{J}_T^*\mathcal{J}_T-\mathcal{J}_T^*\mathcal{J}_N(\mathcal{J}_N^*\mathcal{J}_N)^{-1}\mathcal{J}_N^*\mathcal{J}_T=\mathcal{J}_T^*\bigl[\Id-\mathcal{J}_N(\mathcal{J}_N^*\mathcal{J}_N)^{-1}\mathcal{J}_N^*\bigr]\mathcal{J}_T.
\]
Using $\mathcal{J}_T(r,y)=r\,\partial_r\mathcal{J}_T(0,y)+O(r^2)$, this becomes
\[
 r^2\bigl(\partial_r\mathcal{J}_T(0,y)\bigr)^*(\Id-\Pi_N(y))\partial_r\mathcal{J}_T(0,y)+O(r^3).
\]
The sum of the second-order coefficients of the small
eigenvalues is therefore the trace of this operator.
Dividing by
$\Tr(\mathcal{J}^*\mathcal{J})=\Tr(\mathcal{J}_N(0,y)^*\mathcal{J}_N(0,y))+O(r)$ gives
\begin{equation}\label{egramtransversecoefficient}
 \rho(r,y)^2=r^2 \frac{\Tr\bigl(
   (\partial_r\mathcal{J}_T(0,y))^*(\Id-\Pi_N(y))
   \partial_r\mathcal{J}_T(0,y)\bigr)}
      {\Tr\bigl(\mathcal{J}_N(0,y)^*\mathcal{J}_N(0,y)\bigr)}
 +O(r^3).
\end{equation}
The denominator is positive because $\mathcal{J}_N(0,y)$ is injective.
By \eqref{eboundarycolumntransversality},
$(\Id-\Pi_N(y))\partial_r\mathcal{J}_T(0,y)$ is also injective,
so the numerator is positive.
This proves
\eqref{egramsmalleigenvalueexpansion}.

Near any fixed boundary point, choose $c>0$ such that
\begin{equation*}
 \lambda_{\binom{n-1}{p}+1}(\mathcal T)<c <\lambda_{\binom{n-1}{p}}(\mathcal T).
\end{equation*}
The spectral gap makes the projection
$\Pi_{>c}:=\mathbf{1}_{(c,\infty)}(\mathcal T)$ smooth. Since
\begin{equation*}
 \sum_{k=1}^{\binom{n-1}{p}}\lambda_k(\mathcal T) =\operatorname{Tr}(\mathcal T \Pi_{>c}),\qquad\rho^2=1-\operatorname{Tr}(\mathcal T \Pi_{>c}),
\end{equation*}
$\rho^2$ is smooth even where some of these eigenvalues coincide.

Restoring the subscript $j$, equation~
\eqref{egramtransversecoefficient} gives
\[
 \rho_j^2|_{r=0}=0,\qquad
 \partial_r(\rho_j^2)|_{r=0}=0,\qquad
 \frac12\partial_r^2(\rho_j^2)(0,y)>0.
\]
Taylor's formula gives the exact factorization
 \begin{align}
 \rho_j(r,y)^2&=r^2c_j(r,y),\notag\\
 c_j(r,y)&=\int_0^1(1-t)\,
       \partial_r^2(\rho_j^2)(tr,y)\,dt,\qquad c_j(0,y)>0.\label{esmoothgramdefiningfunction}
 \end{align}
Here $c_j$ is smooth, and its boundary value is the positive
coefficient in \eqref{egramsmalleigenvalueexpansion}.
After extending $c_j$ smoothly and positively across $r=0$,
the function $ \rho_j(r,y)=r\sqrt{c_j(r,y)}$
is a smooth signed extension with
$\partial_r\rho_j(0,y)=\sqrt{c_j(0,y)}>0$.

If the boundary and $g_j$ are analytic near the patch, the
regularity following \eqref{ekernelanalyticestimate} makes
$\mathcal T_j$ analytic up to the boundary in the
Hilbert--Schmidt norm. The
separated spectral projection and its trace are then analytic,
so $\rho_j^2$ is analytic. The factorization
\eqref{esmoothgramdefiningfunction} makes $c_j$ analytic,
and $r\sqrt{c_j(r,y)}$ is the asserted analytic signed extension.
\end{proof}

\begin{lemma}
\label{lemmakernelmapboundaryregularity}
Let $F:M_1\to M_2$ be a homeomorphism satisfying
$\mathcal T_2\circ F=\mathcal T_1$.
Suppose that $y_2=F(y_1)$ and
$y_j\in\partial M_j\cap\partial\widetilde M_j$.
Then $F$ and $F^{-1}$ are smooth near these points.
If the metrics and the boundary patches are analytic there,
$F$ and $F^{-1}$ are analytic there as well.
\end{lemma}

\begin{proof}
Lemma~\ref{lemmaexactkernelboundary} makes $\mathcal T_j$ smooth
there, with injective differential on $T_{y_j}\partial M_j$.
We can choose the same continuous linear functionals
$L_1,\ldots,L_{n-1}$ on the real Hilbert space of self-adjoint
Hilbert--Schmidt operators such that, for $j=1,2$,
\begin{equation}\label{equation::rankcondition}
 \rank\left(
 D\bigl(L_1(\mathcal T_j),\ldots,L_{n-1}(\mathcal T_j)\bigr)(y_j)
 \Big|_{T_{y_j}\partial M_j}\right)=n-1.
\end{equation}
Indeed, for each $j$, the choices of $L_1,\ldots,L_{n-1}$ satisfying
\eqref{equation::rankcondition} form an open dense set.
These two sets have nonempty intersection, so we can choose
the same functionals for both $j=1,2$.
Lemma~\ref{lemmagramboundarydefiningfunction} gives smooth
signed extensions of $\rho_j$ with
\begin{equation*}
 \rho_j|_{\partial M_j}=0,\qquad\partial_r\rho_j(y_j)>0.
\end{equation*}
Together with \eqref{equation::rankcondition}, this shows
that
\begin{equation*}
 x\longmapsto
 \bigl(\rho_j(x),L_1(\mathcal T_j(x)),\ldots,L_{n-1}(\mathcal T_j(x))\bigr)
\end{equation*}
has invertible differential at $y_j$. These functions are
therefore smooth local coordinates by the inverse function
theorem. Since $\mathcal T_2\circ F=\mathcal T_1$,
it follows from the definition \eqref{egramboundarydefiningfunction} that
\begin{equation*}
 \rho_2\circ F=\rho_1,\qquad L_i(\mathcal T_2\circ F)=L_i(\mathcal T_1),\quad 1\leq i\leq n-1.
\end{equation*}
Thus $F$ is the identity in these coordinates, proving
smoothness of $F$ and $F^{-1}$.
If the metrics and boundary patches are analytic, the
regularity following \eqref{ekernelanalyticestimate} and
Lemma~\ref{lemmagramboundarydefiningfunction} make these
coordinates analytic, which proves the last assertion.
\end{proof}

The kernel comparison also identifies the
operators $\mathcal T_j$ in the smooth setting.
Proposition~\ref{czprrecoverexactkernel} identifies
$\mathsf K_j(z,x)$ for $j=1,2$, where $z\in U_0$ and $x$ lies in the added region
or in an interior collar near $\Gamma_1$.
If $n\neq2p+2$, $F_0$ also identifies the dual fibre metrics and $\mathscr H$,
so $\mathcal{J}_1(x)\mathcal{J}_1(x)^*=\mathcal{J}_2(F_0(x))\mathcal{J}_2(F_0(x))^*$.
If $n=2p+2$, the pairing on $\mathscr H$ is conformally invariant and
the dual fibre metric scales by $e^{nf_0(x)}$. Thus
\[
 \mathcal{J}_1(x)\mathcal{J}_1(x)^*
 =e^{nf_0(x)}\mathcal{J}_2(F_0(x))\mathcal{J}_2(F_0(x))^*.
\]
Dividing by the trace shows that, in either case,
\begin{equation}\label{enormalizedgramcollaridentity}
 \mathcal T_1(x)=\mathcal T_2(F_0(x)).
\end{equation}
This holds on $U\cap X_1$ and in interior collars near every
point of $\Gamma_1$.

\subsection{Global continuation and extension to the boundary}
For metrics analytic in the interior, we use the following
local projection construction. Suppose $x_j\in M_j^\circ$ satisfy
\begin{equation}\label{ekerneltangentimages}
 \mathcal T_1(x_1)=\mathcal T_2(x_2),\qquad\mathcal E:=\ran D\mathcal T_1(x_1) =\ran D\mathcal T_2(x_2).
\end{equation}
By Lemma~\ref{lemmaexactkernelembedding},
$\dim\mathcal E=n$.
Let $\Pi_{\mathrm{tan}}$ be the orthogonal projection onto
this fixed subspace. Then
\begin{equation*}
 D(\Pi_{\mathrm{tan}}\mathcal T_j)(x_j) =D\mathcal T_j(x_j):T_{x_j}M_j\longrightarrow\mathcal E
\end{equation*}
is an isomorphism. The analytic inverse function theorem
makes the projected maps local analytic diffeomorphisms.
After restricting to suitable neighborhoods, define
\begin{equation}\label{ekernellocalcontinuation}
 \widehat F :=(\Pi_{\mathrm{tan}}\mathcal T_2)^{-1}\circ(\Pi_{\mathrm{tan}}\mathcal T_1),
\qquad \widehat F(x_1)=x_2,
\end{equation}
where the inverse is taken near $x_2$.
Any smooth local map taking $x_1$ to $x_2$ and satisfying
$\mathcal T_2\circ F=\mathcal T_1$ agrees with this formula
near $x_1$, and is therefore analytic.

\begin{proposition}
\label{rsmoothboundauation}
Let $(M_j,g_j)$, $j=1,2$, be compact connected oriented smooth
$n$-manifolds with nonempty smooth boundaries, where $n\geq2$
and $0\leq p\leq n-2$.
Assume that $g_j$ are smooth up to the boundary and real-analytic
on $M_j^\circ$.
Let $\Gamma_j\subset\partial M_j$ be nonempty open sets
identified by a smooth diffeomorphism
$\kappa:\Gamma_1\to\Gamma_2$, satisfying 
\eqref{equDNintropar}.
Suppose that $\kappa$ extends to a smooth collar diffeomorphism
$F_0$ near $\Gamma_1$.

Then $F_0$ extends to a smooth global diffeomorphism
$F:M_1\to M_2$, real-analytic in the interior, in the following cases:
\begin{enumerate}[label=(\roman*)]
\item If $F_0$ satisfies
\eqref{ecollaralternatives}, then
$ F^*g_2=g_1.$
\item If $n=2p+2$ and $F_0$ satisfies
\eqref{ecollarconformalalternative} with smooth $f_0$, then
there is $f\in C^\infty(M_1)$, real-analytic in $M_1^\circ$, such that
\begin{equation*}
 F^*g_2=e^{2f}g_1,\qquad f|_{\Gamma_1}=0.
\end{equation*}
\end{enumerate}
\end{proposition}
\begin{proof}
\smallskip\noindent
\emph{Step 1. The initial collar map.}
Construct $\widetilde M_j,X_j,U,U_0$ as in
Section~\ref{section::defiitionsext}, extending $F_0$ and the
metrics so that \eqref{ecollaralternatives} or
\eqref{ecollarconformalalternative} also holds on the
added collars.
Proposition~\ref{czprrecoverexactkernel} and
\eqref{enormalizedgramcollaridentity} give
$\mathcal T_1=\mathcal T_2\circ F_0$ in interior collars
near every point of $\Gamma_1$, and on $\Gamma_1$ by
continuity.
Lemmas~\ref{lemmaexactkernelembedding} and
\ref{lemmaexactkernelboundary} make $\mathcal T_j$
continuous and injective on $X_j$, and analytic immersions
on $M_j^\circ$.
Differentiating $\mathcal T_1=\mathcal T_2\circ F_0$ gives
\begin{equation*}
 D\mathcal T_1(x)=D\mathcal T_2(F_0(x))\circ DF_0(x).
\end{equation*}
Thus \eqref{ekerneltangentimages} holds at $x,F_0(x)$,
and \eqref{ekernellocalcontinuation} makes $F_0$
analytic on the observed interior collars.

\smallskip\noindent
\emph{Step 2. Continuation through the interior.} We adapt the proof of \cite[Theorem~3.7]{KrupchykLassasUhlmann}.
Let $\mathcal O\subset M_1^\circ$ be the maximal connected
open set containing a connected observed interior collar
on which
\begin{equation}\label{olddefinitionF}
 \mathcal T_1(\mathcal O)\subset\mathcal T_2(M_2^\circ),\qquad F:=\bigl(\mathcal T_2|_{M_2^\circ}\bigr)^{-1}\circ\mathcal T_1: \mathcal O\longrightarrow M_2^\circ
\end{equation}
is an analytic local diffeomorphism, isometric when
\eqref{ecollaralternatives} holds and conformal when \eqref{ecollarconformalalternative} holds.
The inverse is taken on the image of $\mathcal T_2$.

Suppose $x_k\in\mathcal O$ converges to
$x_1\in\partial\mathcal O\cap M_1^\circ$.
Compactness and continuity give, after passing to a subsequence,
\begin{equation*}
 F(x_k)\longrightarrow x_2\in M_2,\qquad
 \mathcal T_1(x_1)=\mathcal T_2(x_2).
\end{equation*}
If $x_2\in\partial M_2\cap\partial\widetilde M_2$, the
interior and boundary ranks give the contradiction
\begin{equation*}
 \binom{n}{p+1}=\rank\mathcal T_1(x_1)=\rank\mathcal T_2(x_2)=\binom{n-1}{p}.
\end{equation*}
If $x_2\in\partial M_2\cap\widetilde M_2^\circ$, it lies
on the attached observed patch. By
\eqref{enormalizedgramcollaridentity} and injectivity,
\begin{equation*}
 \mathcal T_1(x_1)=\mathcal T_1(F_0^{-1}(x_2))\quad\Longrightarrow\quad
 x_1=F_0^{-1}(x_2)\in\partial M_1,
\end{equation*}
again a contradiction. Hence $x_2\in M_2^\circ$.

Differentiating $\mathcal T_1=\mathcal T_2\circ F$ on
$\mathcal O$ gives
\begin{equation*}
 \ran D\mathcal T_1(x_k)
 =\ran D\mathcal T_2(F(x_k)).
\end{equation*}
The differentials are continuous in the Hilbert--Schmidt
norm and have constant rank $n$. Their ranges therefore
converge to the corresponding ranges at $x_1,x_2$, so
\eqref{ekerneltangentimages} holds there.
Formula~\eqref{ekernellocalcontinuation} defines an
analytic local diffeomorphism $\widehat F$ on a small
connected neighborhood of $x_1$, with $\widehat F(x_1)=x_2$.

For sufficiently large $k$, uniqueness of the projected
local inverse implies $\widehat F=F$ near $x_k$.
Thus
\begin{equation*}
 \mathcal T_2\circ\widehat F=\mathcal T_1,
 \end{equation*}
and
 \begin{equation*}
 \widehat F^*g_2=
 \begin{cases}
 g_1,&\text{in Case (i)},\\[2pt]
 \dfrac1n\tr_{g_1}(\widehat F^*g_2)\,g_1,&\text{in Case (ii)},
 \end{cases}
\end{equation*}
hold on a nonempty open set and hence throughout that
neighborhood by analyticity.
By injectivity of $\mathcal T_2$, $\widehat F=F$ on the
entire overlap with $\mathcal O$.
This extends $F$ across $x_1$, contradicting maximality.
Connectedness implies $\mathcal O=M_1^\circ$.

Interchanging $j=1,2$ and starting from $F_0^{-1}$ gives
the reverse inclusion of the interior images. Thus
\begin{equation*}
 \mathcal T_1(M_1^\circ)=\mathcal T_2(M_2^\circ),
\end{equation*}
and $F:M_1^\circ\to M_2^\circ$ is an analytic
diffeomorphism. Continuity and compactness then give
\begin{equation*}
 \mathcal T_1(M_1)=\overline{\mathcal T_1(M_1^\circ)}=\overline{\mathcal T_2(M_2^\circ)}
 =\mathcal T_2(M_2).
\end{equation*}
Since $\mathcal T_j|_{M_j}$ are topological embeddings,
\eqref{olddefinitionF} extends $F$ to a homeomorphism
$M_1\to M_2$. The collar identity and injectivity give
$F=F_0$ near $\Gamma_1$, so $F|_{\Gamma_1}=\kappa$.

\smallskip\noindent
\emph{Step 3. Boundary regularity and the metric identity.}
Near $\Gamma_1$, smoothness follows from $F=F_0$.
If $y_1\in\partial M_1\setminus\Gamma_1$ and $y_2=F(y_1)$,
then $y_2\notin\Gamma_2$, since $F^{-1}=F_0^{-1}$ there.
Both points lie on unchanged original boundary patches, outside $\supp\chi$.
Lemma~\ref{lemmakernelmapboundaryregularity} therefore
makes $F$ and $F^{-1}$ smooth at these points.

The interior metric identity extends to $M_1$ by continuity, giving a global isometry when
\eqref{ecollaralternatives} holds. When \eqref{ecollarconformalalternative} holds, set
\begin{equation*}
 f=\frac12\log\left(\frac1n\tr_{g_1}(F^*g_2)\right).
\end{equation*}
The trace is positive and smooth up to the boundary, and
analytic in the interior. Thus $f$ has the asserted
regularity and $F^*g_2=e^{2f}g_1$.
Near $\Gamma_1$, $F=F_0$ and hence $f=f_0$, so
$f|_{\Gamma_1}=0$.
\end{proof}

We now prove Theorem~\ref{thmsnalyticcalderon} by combining the preceding propositions.
\begin{proof}[Proof of Theorem~\ref{thmsnalyticcalderon}]
Proposition~\ref{analyticcollarrigidity} supplies an
analytic collar map $F_0$ satisfying
\eqref{ecollaralternatives} or
\eqref{ecollarconformalalternative}. Since the metrics are analytic in the interior and
smooth up to the boundary,
Proposition~\ref{rsmoothboundauation}
extends it to a smooth global isometric or conformal
diffeomorphism $F$, analytic in the interior.

Near $\Gamma_1$, we have $F=F_0$, so $F$ is
real-analytic up to the observed boundary $\Gamma_1$.
In the critical case, the conformal factor $f$ is
smooth globally and analytic in the interior.
Since $f=f_0$ near $\Gamma_1$, it is also
analytic up to $\Gamma_1$, with $f|_{\Gamma_1}=0$.
\end{proof}

\section{Boundary unique continuation and rigidity}
\label{sectioncurvatureinverse}

This section establishes the boundary unique continuation results
used in the proof of Theorem~\ref{thmharmoniccurvaturein}.
We first prove a unique continuation lemma and
then apply it to metrics with harmonic curvature whose boundary
jets agree.
Throughout this section $n\geq3$.

\subsection{A common boundary unique continuation lemma}

The following lemma treats an elliptic inequality coupled
to an ordinary differential inequality in the inward direction.
The geometric arguments below
will specify the arrays $\mathbf{u},\mathbf{v}$ and verify these inequalities.

\begin{lemma}\label{lemmacoupledcurvatureuc}
Let $g$ be a smooth Riemannian metric near $(0,0)$ in local
boundary coordinates $(r,y)$, with $r\geq0$ on the interior side.
Assume that $\partial_r$ is the inward $g$-unit normal on $r=0$.
Set $L:=g^{\mu\nu}\partial_\mu\partial_\nu$.
Suppose that finite arrays $\mathbf{u},\mathbf{v}$, smooth for $r\geq0$ up to
$r=0$, satisfy
\begin{equation}\label{ecoupledcurvatureucsystem}
 |L\mathbf{u}|\leq C(|\mathbf{u}|+|\partial \mathbf{u}|+|\mathbf{v}|),\qquad
 |\partial_r\mathbf{v}|\leq C(|\mathbf{v}|+|\mathbf{u}|+|\partial \mathbf{u}|)
\end{equation}
and
\begin{equation}\label{ecoupledcurvatureuctraces}
 \mathbf{u}|_{r=0}=\partial_r\mathbf{u}|_{r=0}=\mathbf{v}|_{r=0}=0.
\end{equation}
Here $\partial$ denotes all first coordinate derivatives, and $C>0$ is uniform in $(r,y)$ but may depend on $g,\mathbf u,\mathbf v$.
Then $\mathbf{u}=\mathbf{v}=0$ in a smaller neighborhood of $(0,0)$ with $r\geq0$.
\end{lemma}

\begin{proof}
Extend the boundary coordinate neighborhood across $r=0$ and
extend $g$ smoothly to this larger neighborhood. After shrinking,
the extended metric is positive definite. Retain the notation
$g$ and $L$ for the extended metric and operator.

Extend $\mathbf{u},\mathbf{v}$ by zero to $r<0$. The traces in
\eqref{ecoupledcurvatureuctraces} make these extensions belong to
$H^2_{\mathrm{loc}}$ and $H^1_{\mathrm{loc}}$, respectively,
and ensure that no boundary distributions occur in $L\mathbf{u}$ or
$\partial_r\mathbf{v}$. Since both extensions vanish for $r<0$,
\eqref{ecoupledcurvatureucsystem} holds almost everywhere
on the enlarged coordinate neighborhood.
For a sufficiently large fixed $\lambda$, set
\[
 \psi(r,y):=-r-|y|^2+\frac{\lambda}{2}r^2.
\]
At $(0,0)$, with respect to the extended metric $g$,
\[
 \nabla\psi=-\partial_r,\qquad
 \nabla^2\psi(\partial_r,\partial_r)=\lambda+O(1),
\]
where the $O(1)$ term is independent of $\lambda$. The quantity
$\nabla^2\psi(V,W)$ is also independent of $\lambda$
for vectors $V,W$ tangent to the boundary at $(0,0)$.
Choose $\lambda$ large and shrink the neighborhood so that
\begin{equation}\label{equation::hessianineq}
 \nabla^2\psi(V,V)+\nabla^2\psi(\nabla\psi,\nabla\psi)>0,\quad V\perp\nabla\psi,\quad |V|_g=|\nabla\psi|_g,\qquad \partial_r\psi\leq-\frac12.
\end{equation}
For $L=g^{\mu\nu}\partial_\mu\partial_\nu$, the quadratic
symbol of \[-e^{\tau\psi}Le^{-\tau\psi}=-L
 +2\tau g^{\mu\nu}(\partial_\mu\psi)\partial_\nu-\tau^2g^{\mu\nu}
   (\partial_\mu\psi)(\partial_\nu\psi)+\tau L\psi\] in $(\xi,\tau)$ is
\begin{align*}
 q_\tau:&=g^{\mu\nu}\xi_\mu\xi_\nu+2\sqrt{-1}\tau g^{\mu\nu}
     (\partial_\mu\psi)\xi_\nu-\tau^2g^{\mu\nu}(\partial_\mu\psi)(\partial_\nu\psi)\\
 &=|\xi|_{g}^2-\tau^2|\nabla\psi|_{g}^2 +2\sqrt{-1}\tau\, g(\xi^{\sharp_{g}},\nabla\psi).
 \end{align*}
On $q_\tau=0$ with $\tau>0$, the vector
$V=\xi^{\sharp_{g}}/\tau$ satisfies
$V\perp\nabla\psi$ and $|V|_{g}=|\nabla\psi|_{g}$.
The inequality \eqref{equation::hessianineq} therefore implies
\[
 \{\Re q_\tau,\Im q_\tau\}=4\tau^3\bigl( \nabla^2\psi(V,V) +\nabla^2\psi(\nabla\psi,\nabla\psi) \bigr)>0,
\]
where $\{\cdot,\cdot\}$ denotes the Poisson bracket.
Since $\partial_r\psi\leq-\frac12$, the Carleman weight $\psi$ has no critical
points in this neighborhood. The bracket positivity on $q_\tau=0$
for $\tau>0$ thus verifies the sub-ellipticity condition of
\cite[Definition~3.2]{LeRousseauLebeauRobbiano2022} for
\[
 -L=-\partial_\mu(g^{\mu\nu}\partial_\nu)+(\partial_\mu g^{\mu\nu})\partial_\nu.
\]
The Carleman
estimate \cite[Theorem~3.11]{LeRousseauLebeauRobbiano2022},
together with its $H^2$ extension
\cite[Remark~3.12]{LeRousseauLebeauRobbiano2022}, gives
\begin{equation}\label{eharmoniccurvaturecarleman}
 \tau^{\frac32}\|e^{\tau\psi}w\|_{L^2(drdy)}+\tau^{\frac12}\|e^{\tau\psi}\partial w\|_{L^2(drdy)}\leq C\|e^{\tau\psi}Lw\|_{L^2(drdy)}
\end{equation}
for compactly supported $w\in H^2$ and sufficiently large $\tau$,
by dropping the nonnegative second-derivative term in the squared
estimate and taking square roots.
Since $L$ acts componentwise, estimate \eqref{eharmoniccurvaturecarleman} holds for finite arrays by summing over their components.
For compactly supported $w\in H^1$, integration by parts also gives
\[
 \int e^{2\tau\psi}(\partial_rw)\cdot w\,drdy
 =-\tau\int\partial_r\psi\,|e^{\tau\psi}w|^2\,drdy
 \geq\frac{\tau}{2}\|e^{\tau\psi}w\|_{L^2(drdy)}^2.
\]
It follows from the Cauchy--Schwarz inequality that
\begin{equation}\label{eharmoniccurvatureodecarleman}
 \frac{\tau}{2}\|e^{\tau\psi}w\|_{L^2(drdy)}\leq\|e^{\tau\psi}\partial_rw\|_{L^2(drdy)}.
\end{equation}

Choose a product cutoff $\chi_\varepsilon$ for $\varepsilon>0$, equal to one when
$|r|\leq\varepsilon$ and $|y|\leq\varepsilon$, and supported
where $|r|<2\varepsilon$ and $|y|<2\varepsilon$.
For sufficiently small $\varepsilon$,
\[
 \psi\leq-\frac r2-|y|^2\quad(r\geq0).
\]
On the support of a derivative of $\chi_\varepsilon$ with $r\geq0$, we therefore have $\psi\leq-\min\{\varepsilon/2,\varepsilon^2\}$. Set $c_*:=\min\{\varepsilon/2,\varepsilon^2\}$.
Since $\mathbf{u}=\mathbf{v}=0$ for $r<0$, and every derivative of $\chi_\varepsilon$ is supported where $\psi\leq-c_*$ on $r\geq0$, we have the estimate
\begin{equation}\label{equation::cutofferr}
\|e^{\tau\psi}[L,\chi_\eps]\mathbf{u}\|_{L^2(drdy)}+\|e^{\tau\psi}(\partial_r\chi_\eps)\mathbf{v}\|_{L^2(drdy)}\leq Ce^{-\tau c_*},
\end{equation}
where $C$ is independent of $\tau$.
Apply \eqref{eharmoniccurvaturecarleman} to $\chi_\varepsilon \mathbf{u}$ and
\eqref{eharmoniccurvatureodecarleman} to $\chi_\varepsilon \mathbf{v}$.
The identities
\begin{equation*}
 L(\chi_\varepsilon \mathbf{u})=\chi_\varepsilon L\mathbf{u}+[L,\chi_\varepsilon]\mathbf{u},\qquad\partial_r(\chi_\varepsilon \mathbf{v})
 =\chi_\varepsilon\partial_r\mathbf{v}+(\partial_r\chi_\varepsilon)\mathbf{v},
\end{equation*}
together with \eqref{equation::cutofferr} and
\eqref{ecoupledcurvatureucsystem}, give
 \begin{align*}
 &\|e^{\tau\psi}L(\chi_\varepsilon \mathbf{u})\|_{L^2(drdy)}+\|e^{\tau\psi}\partial_r(\chi_\varepsilon \mathbf{v})\|_{L^2(drdy)}\\
 &\quad\leq\|e^{\tau\psi}\chi_\varepsilon L\mathbf{u}\|_{L^2(drdy)}
 +\|e^{\tau\psi}\chi_\varepsilon\partial_r\mathbf{v}\|_{L^2(drdy)}
 +Ce^{-\tau c_*}\\
 &\quad\leq C\bigl(\|e^{\tau\psi}\chi_\varepsilon \mathbf{u}\|_{L^2(drdy)} +\|e^{\tau\psi}\chi_\varepsilon\partial \mathbf{u}\|_{L^2(drdy)}+\|e^{\tau\psi}\chi_\varepsilon \mathbf{v}\|_{L^2(drdy)}
 \bigr)
 +Ce^{-\tau c_*}.
 \end{align*}
To control the derivative term, use
\begin{equation*}
 \chi_\varepsilon\partial \mathbf{u}=\partial(\chi_\varepsilon \mathbf{u})-(\partial\chi_\varepsilon)\mathbf{u}.
\end{equation*}
Since $\mathbf{u}=0$ for $r<0$ and $\psi\leq-c_*$ on
$\operatorname{supp}(\partial\chi_\varepsilon)\cap\{r\geq0\}$,
the estimate \eqref{equation::cutofferr} yields
\begin{equation*}
 \|e^{\tau\psi}\chi_\varepsilon\partial \mathbf{u}\|_{L^2(drdy)}\leq\|e^{\tau\psi}\partial(\chi_\varepsilon \mathbf{u})\|_{L^2(drdy)}+Ce^{-\tau c_*}.
\end{equation*}
Substituting this bound and adding the two Carleman estimates,
we obtain
\begin{equation*}
 \begin{aligned}
 &\tau^{\frac32}\|e^{\tau\psi}\chi_\varepsilon \mathbf{u}\|_{L^2(drdy)}+\tau^{\frac12}
 \|e^{\tau\psi}\partial(\chi_\varepsilon \mathbf{u})\|_{L^2(drdy)} +\frac{\tau}{2}
 \|e^{\tau\psi}\chi_\varepsilon \mathbf{v}\|_{L^2(drdy)}\\
 &\quad\leq C\bigl(\|e^{\tau\psi}\chi_\varepsilon \mathbf{u}\|_{L^2(drdy)}
 +\|e^{\tau\psi}\partial(\chi_\varepsilon \mathbf{u})\|_{L^2(drdy)}+\|e^{\tau\psi}\chi_\varepsilon \mathbf{v}\|_{L^2(drdy)}
 \bigr)+Ce^{-\tau c_*},
 \end{aligned}
\end{equation*}
where $C$ is independent of $\tau$.
For sufficiently large $\tau$, the three norm terms on the
right can be absorbed into the left. After adjusting $C$,
this implies
\begin{equation}\label{eharmoniccurvaturelocalvanishing}
\tau^{\frac32}
 \|e^{\tau\psi}\chi_\varepsilon \mathbf{u}\|_{L^2(drdy)}+\tau^{\frac12}
 \|e^{\tau\psi}\partial(\chi_\varepsilon \mathbf{u})\|_{L^2(drdy)} +\tau\|e^{\tau\psi}\chi_\varepsilon \mathbf{v}\|_{L^2(drdy)}\leq Ce^{-\tau c_*}.
\end{equation}
On a smaller neighborhood with $r\geq0$, $\chi_\varepsilon=1$ and
$\psi>-c_*/2$. Restricting
\eqref{eharmoniccurvaturelocalvanishing} there gives,
with the following norms taken over this smaller neighborhood,
\[
 \|\mathbf{u}\|_{L^2(drdy)}+\|\mathbf{v}\|_{L^2(drdy)} \leq Ce^{-\tau c_*/2}.
\]
Letting $\tau\to\infty$ proves the lemma.
\end{proof}

\subsection{Boundary rigidity for metrics with harmonic curvature}

We write $\nabla^g$ for the Levi--Civita connection of $g$
and its induced connections on tensor bundles, omitting the
superscript when the metric is clear from the context.

A metric has harmonic curvature if its Riemannian curvature tensor $R$ is divergence-free:
\begin{equation}\label{eharmoniccurvature}
 \nabla^iR_{ik\ell m}=0.
\end{equation}
Such a metric has a Codazzi Ricci tensor and constant scalar curvature
on each connected component (see \cite[Section~2.6]{DerdzinskiHarmonicCurvature1985}).
The second Bianchi identity reads
\[
 \nabla_sR_{ik\ell m} +\nabla_iR_{ks\ell m} +\nabla_kR_{si\ell m}=0.
\]
The curvature tensor therefore satisfies a second-order equation.
Indeed,
\eqref{eharmoniccurvature} and antisymmetry in the first
two curvature indices imply
\[
 \nabla^sR_{ks\ell m}=-\nabla^sR_{sk\ell m}=0.
\]
Applying $\nabla^s$ to the second Bianchi identity yields
 \begin{align}
 \nabla^s\nabla_sR_{ik\ell m} &=-\nabla^s\nabla_iR_{ks\ell m}-\nabla^s\nabla_kR_{si\ell m}\notag\\
 &=-[\nabla^s,\nabla_i]R_{ks\ell m}-\nabla_i\nabla^sR_{ks\ell m}-[\nabla^s,\nabla_k]R_{si\ell m} -\nabla_k\nabla^sR_{si\ell m}\notag\\
 &=-[\nabla^s,\nabla_i]R_{ks\ell m}-[\nabla^s,\nabla_k]R_{si\ell m}.\label{eharmoniccurvatureclosedequation}
 \end{align}
The covariant derivative commutation formula is
\[
 [\nabla_a,\nabla_b]R_{ik\ell m} =-R_{abi}{}^qR_{qk\ell m} -R_{abk}{}^qR_{iq\ell m}
    -R_{ab\ell}{}^qR_{ikqm} -R_{abm}{}^qR_{ik\ell q}.
\]
Thus the right-hand side of
\eqref{eharmoniccurvatureclosedequation}
is quadratic in the curvature components and contains no
derivatives of curvature.

We next apply Lemma~\ref{lemmacoupledcurvatureuc} to show
that matching boundary metric jets through order two imply
equality of the metrics near the boundary. Related boundary
unique continuation results for Einstein metrics were proved in
\cite{AndersonHerzlich2008,AndersonHerzlich2010Erratum}.

\begin{proposition}\label{lemmacurvatureboundaryuc}
Let $g_j=dr^2+h_{j,r}$, $j=1,2$, be smooth metrics on a common
collar as in Section~\ref{section::commoncollarset}, satisfying
\eqref{eharmoniccurvature}. If, on an open boundary patch,
\begin{equation}\label{eharmoniccurvaturefourjets}
 \left.\partial_r^kh_{1,r}\right|_{r=0}
 =\left.\partial_r^kh_{2,r}\right|_{r=0},
 \qquad k=0,1,2,
\end{equation}
then $g_1=g_2$ near that patch.
\end{proposition}

\begin{proof}
We work near a fixed point of the boundary patch, in common
coordinates $(r,y)$, with $r\geq0$ inside the collar.
The proof couples an elliptic equation for the curvature to
transport equations for the metric, following the method of
\cite[Sections~2--3]{WongYu2012}.

\smallskip\noindent
\emph{Step 1. The curvature equation and normal transport.}
For this transport argument we use normal-parallel orthonormal
frames, while the collar coordinates $(r,y)$ remain fixed.
For each metric, extend the same tangential orthonormal frame at
$r=0$ by parallel transport along the normal geodesics, and adjoin
$e_{j,n}=\partial_r$. Write
\[
 e_{j,i}=e_{j,i}^{\mu}\partial_\mu,\qquad\gamma_{j,ik}^{\ell} =g_j(\nabla^{g_j}_{e_{j,i}}e_{j,k},e_{j,\ell}),\qquad  R_{j,ik\ell m} =g_j(R^{g_j}(e_{j,i},e_{j,k})e_{j,\ell},e_{j,m}).
\]
Frame indices range from $1$ to $n$, and $\mu$ ranges over the
coordinates $(r,y)$. Suppressing $j$, recall that $e_n=\partial_r$ and
$\nabla_{\partial_r}e_i=0$. Torsion freeness therefore gives
\[
 [\partial_r,e_i]=\nabla_{\partial_r}e_i-\nabla_{e_i}\partial_r=-\gamma_{in}^{s}e_s.
\]
Since $e_i=e_i^\mu\partial_\mu$, we obtain
\[
 (\partial_r e_i^\mu)\partial_\mu=[\partial_r,e_i]=-\gamma_{in}^{s}e_s^\mu\partial_\mu.
\]
For the connection coefficients, using metric compatibility,
normal parallel transport, and the curvature identity, we compute
 \begin{align*}
 \partial_r\gamma_{ik}^{\ell}&=\partial_r\,g(\nabla_{e_i}e_k,e_\ell)\\
 &=g(\nabla_{\partial_r}\nabla_{e_i}e_k,e_\ell)+g(\nabla_{e_i}e_k,\nabla_{\partial_r}e_\ell)\\
 &=g\bigl(R(\partial_r,e_i)e_k +\nabla_{[\partial_r,e_i]}e_k,e_\ell\bigr)\\
 &=R_{nik\ell}-\gamma_{in}^{s}g(\nabla_{e_s}e_k,e_\ell)\\
 &=R_{nik\ell}-\gamma_{in}^{s}\gamma_{sk}^{\ell}.
 \end{align*}
Here we used
$\nabla_{\partial_r}e_k=\nabla_{\partial_r}e_\ell=0$
and the bracket identity above.
Consequently, the frame and connection coefficients satisfy
\begin{equation}\label{eharmoniccurvaturetransport}
 \partial_r e_i^\mu=-\gamma_{in}^{s}e_s^\mu,\qquad\partial_r\gamma_{ik}^{\ell}=R_{nik\ell}-\gamma_{in}^{s}\gamma_{sk}^{\ell}.
\end{equation}
Thus curvature determines the normal evolution of the connection
and frame coefficients.

\smallskip\noindent
\emph{Step 2. A closed system for the differences.}
Define the full component arrays
\begin{align*}
 \mathbf{u}
 &:=\bigl(R_{1,ik\ell m}-R_{2,ik\ell m}\bigr)_{i,k,\ell,m},\\
 \mathbf{v}
 &:=\Bigl(
 \bigl(e_{1,i}^{\mu}-e_{2,i}^{\mu}\bigr)_{i,\mu},
 \bigl(\gamma_{1,ik}^{\ell}-\gamma_{2,ik}^{\ell}\bigr)_{i,k,\ell},\\
 &\qquad
 \bigl(\partial_{y^\alpha}(e_{1,i}^{\mu}-e_{2,i}^{\mu})
       \bigr)_{\alpha,i,\mu},
 \bigl(\partial_{y^\alpha}(\gamma_{1,ik}^{\ell}
                         -\gamma_{2,ik}^{\ell})
       \bigr)_{\alpha,i,k,\ell}
 \Bigr).
\end{align*}
Each indexed parenthesis denotes the array of all indicated
components, without summation. The frame indices range independently
from $1$ to $n$, $\mu$ ranges over all coordinate directions,
and $1\leq\alpha\leq n-1$.
We write $\partial_y$ for the collection of first tangential
derivatives and $\partial$ for the collection of all first
coordinate derivatives.

Use Euclidean norms for these arrays. We work on a smaller collar with compact closure, so that all
coefficients and derivatives appearing below are bounded. Set $L:=g_1^{\mu\nu}\partial_\mu\partial_\nu$. We claim that
\begin{equation}\label{eharmoniccurvaturedifferencesystem}
   |L\mathbf{u}|\leq C(|\mathbf{u}|+|\partial \mathbf{u}|+|\mathbf{v}|),\qquad  |\partial_r\mathbf{v}|\leq C(|\mathbf{v}|+|\mathbf{u}|+|\partial \mathbf{u}|).
\end{equation}

Since the frames are orthonormal,
$ g_j^{\mu\nu}
 =\sum_i e_{j,i}^{\mu}e_{j,i}^{\nu}.$
Thus
\[
 |g_1^{\mu\nu}-g_2^{\mu\nu}|=\left| \sum_i\bigl[(e_{1,i}^{\mu}-e_{2,i}^{\mu})e_{1,i}^{\nu} +e_{2,i}^{\mu}(e_{1,i}^{\nu}-e_{2,i}^{\nu})\bigr]\right|\leq C|\mathbf{v}|.
\]
In particular,
\begin{equation}\label{leadingtermcur}
 \bigl|(g_1^{\mu\nu}-g_2^{\mu\nu})\partial_\mu\partial_\nu R_{2,ik\ell m}\bigr|\leq C|\mathbf{v}|.
\end{equation}
Subtracting \eqref{eharmoniccurvaturetransport} for the two metrics, we obtain
\[
 |\partial_r(e_{1,i}^{\mu}-e_{2,i}^{\mu})|\leq C|\mathbf{v}|,\qquad |\partial_r(\gamma_{1,ik}^{\ell}-\gamma_{2,ik}^{\ell})|
 \leq C(|\mathbf{u}|+|\mathbf{v}|).
\]
Together with the tangential derivatives included in $\mathbf{v}$,
these estimates bound the differences of the frame and
connection coefficients and their first derivatives by
$C(|\mathbf{u}|+|\mathbf{v}|)$.

In components, the second-order term on the left of
\eqref{eharmoniccurvatureclosedequation} is
$g_j^{\mu\nu}\partial_\mu\partial_\nu R_{j,ik\ell m}$.
The lower-order terms are linear in the curvature components
and their first derivatives.
Since the differences of their coefficients are bounded
by $C(|\mathbf{u}|+|\mathbf{v}|)$, the differences of these terms are bounded
by $C(|\mathbf{u}|+|\partial \mathbf{u}|+|\mathbf{v}|)$.
The right-hand side is quadratic in curvature, so its
difference is bounded by $C|\mathbf{u}|$.

Since
 \begin{align*}
 (L\mathbf{u})_{ik\ell m}&=g_1^{\mu\nu}\partial_\mu\partial_\nu (R_{1,ik\ell m}-R_{2,ik\ell m})\\
 &=g_1^{\mu\nu}\partial_\mu\partial_\nu R_{1,ik\ell m} -g_2^{\mu\nu}\partial_\mu\partial_\nu R_{2,ik\ell m}-(g_1^{\mu\nu}-g_2^{\mu\nu}) \partial_\mu\partial_\nu R_{2,ik\ell m},
 \end{align*}
these bounds, together with \eqref{leadingtermcur}, give
\[
 |L\mathbf{u}|\leq C(|\mathbf{u}|+|\partial \mathbf{u}|+|\mathbf{v}|).
\]

Finally, subtracting
\eqref{eharmoniccurvaturetransport} for the two metrics
and differentiating tangentially introduces $\partial_y\mathbf{u}$,
while the other differences are controlled by $\mathbf{v}$.
Together with the undifferentiated equations, this yields
\[
 |\partial_r\mathbf{v}|\leq C(|\mathbf{v}|+|\mathbf{u}|+|\partial \mathbf{u}|),
\]
completing the proof of the claim
\eqref{eharmoniccurvaturedifferencesystem}.

The metric derivatives in \eqref{eharmoniccurvaturefourjets}
determine the boundary curvature and connection coefficients.
Their tangential derivatives agree as well, so
\[
 \mathbf{u}|_{r=0}=0,\qquad \mathbf{v}|_{r=0}=0.
\]

To apply Lemma~\ref{lemmacoupledcurvatureuc}, it remains to prove
$\partial_r\mathbf{u}|_{r=0}=0$.
The boundary conditions $\mathbf{u}|_{r=0}=\mathbf{v}|_{r=0}=0$
ensure equality of the curvature components, frame coefficients,
and connection coefficients on the boundary patch.
Tangential derivatives of the curvature component functions
therefore agree as well. Since $e_{j,s}$ is tangent to the
boundary for $s<n$, it follows that
\[
 \left.(\nabla^{g_1}_{e_{1,s}}R_1)_{ik\ell m}\right|_{r=0}
 =
 \left.(\nabla^{g_2}_{e_{2,s}}R_2)_{ik\ell m}\right|_{r=0},
 \qquad s<n.
\]
The second Bianchi identity and the harmonic curvature give, respectively,
\begin{align*}
 \partial_rR_{j,ik\ell m}
 &=
 -(\nabla^{g_j}_{e_{j,i}}R_j)_{kn\ell m}
 -(\nabla^{g_j}_{e_{j,k}}R_j)_{ni\ell m},
 &&i,k<n,\\
 \partial_rR_{j,ni\ell m}
 &=
 -\sum_{s=1}^{n-1}
   (\nabla^{g_j}_{e_{j,s}}R_j)_{si\ell m},
 &&i<n.
\end{align*}
Every term on the right has the same boundary value for $j=1,2$. Subtracting the corresponding equations yields
\[
 \left.\partial_r\mathbf{u}_{ik\ell m}\right|_{r=0}=0,
 \qquad
 \left.\partial_r\mathbf{u}_{ni\ell m}\right|_{r=0}=0,
 \qquad i,k<n.
\]
By antisymmetry in the first two curvature indices, these
identities imply $\partial_r\mathbf{u}|_{r=0}=0$.

Lemma~\ref{lemmacoupledcurvatureuc} now applies to
\eqref{eharmoniccurvaturedifferencesystem} and gives $\mathbf{u}=\mathbf{v}=0$
near the chosen boundary point. In particular, the two frame
coefficient arrays agree. Since
$g_j^{\mu\nu}=\sum_i e_{j,i}^\mu e_{j,i}^\nu$, the metrics agree
there as well.
\end{proof}

\section{Conformal rigidity and global uniqueness}
\label{sectioncriticalinverse}

We complete the proof of Theorem~\ref{thmharmoniccurvaturein}.
In the critical dimension $n=2p+2$, this requires converting
the conformal boundary jets into a local conformal map.
We do so using conformal geodesics and connection matrices.
Subsection~\ref{subsectioncriticalglobal} combines this construction
with Section~\ref{sectioncurvatureinverse} and
Proposition~\ref{rsmoothboundauation} to obtain global uniqueness.

\subsection{Curvature identities}

Throughout the local conformal arguments, $n\geq4$.
Let $\mathrm{Sch}_g$, $W_g$, and $C_g$ denote the Schouten,
Weyl, and Cotton tensors, respectively. With our curvature convention,
\begin{align}
 \mathrm{Sch}_g
 &:=\frac1{n-2}\left(\operatorname{Ric}_g
                  -\frac{\mathrm{Scal}_g}{2(n-1)}g\right),\notag\\
 (W_g)_{iqk\ell}
 &:=R_{iqk\ell}
 -g_{i\ell}(\mathrm{Sch}_g)_{qk}
 -g_{qk}(\mathrm{Sch}_g)_{i\ell}
 +g_{ik}(\mathrm{Sch}_g)_{q\ell}
 +g_{q\ell}(\mathrm{Sch}_g)_{ik},\label{defweyltensor}\\
 (C_g)_{iqk}
 &:=\nabla_i(\mathrm{Sch}_g)_{qk}
            -\nabla_q(\mathrm{Sch}_g)_{ik}.\notag
\end{align}
Under a conformal change, the Cotton tensor transforms by (see, for instance,
\cite[Eqs.~(2.5), (3.2)]{CatinoMastroliaMonticelliPunzo2023})
\begin{equation}\label{cottontransformfor}
 (C_{e^{2f}g})_{iqk}=(C_g)_{iqk} +(W_g)_{iqk\ell}(\nabla^g f)^\ell.
\end{equation}
The covariant exterior derivative and covariant codifferential
of $W_g$ are given by
\begin{align*}
 (d^\nabla W_g)_{siqk\ell}
 &:=(\nabla_sW_g)_{iqk\ell}
   +(\nabla_iW_g)_{qsk\ell}
   +(\nabla_qW_g)_{sik\ell},\\
 (\delta^\nabla W_g)_{iqk}
 &:=-\nabla^s(W_g)_{siqk}.
\end{align*}
These are tensors of types $(0,5)$ and $(0,3)$, respectively. The same definitions apply to any $(0,4)$-tensor antisymmetric in its first two indices.
With our curvature conventions, the second
Bianchi identity for the Weyl tensor and its contraction
give (see also \cite[Eqs.~(4)--(5)]{GoverNagy2007})
\begin{equation}\label{ecriticalweylcottonequations}
 \begin{aligned}
 (d^\nabla W_g)_{siqk\ell}
 &=\sum_{\mathrm{cyc}(s,i,q)}
   \bigl[g_{q\ell}(C_g)_{sik}
         -g_{qk}(C_g)_{si\ell}\bigr],\\
 (\delta^\nabla W_g)_{iqk}
 &=(n-3)(C_g)_{qki}.
 \end{aligned}
\end{equation}
Thus both tensors are algebraic linear expressions in $C_g$.

If $g$ has harmonic curvature, $\mathrm{Sch}_g$ is Codazzi,
so $C_g=0$. Consequently,
\begin{equation}\label{eharmonicweylidentities}
 d^\nabla W_g=0,\qquad \delta^\nabla W_g=0.
\end{equation}

We write $W_g(V)$ for $W_g(V,\cdot,\cdot,\cdot)$.
By the curvature symmetries, vanishing of $W_g(V)$ is equivalent to
vanishing with $V$ inserted in any other slot.
As usual, $O(r^\infty)$ means that every boundary jet
vanishes on the boundary open patch in question.

We record two consequences of the Codazzi and Bianchi identities.
Here $(\mathrm{Sch}_g)^\sharp$ is defined by
\[
 g\bigl((\mathrm{Sch}_g)^\sharp Z,U\bigr)=\mathrm{Sch}_g(Z,U).
\]

\begin{lemma}
Let $g$ be a smooth metric with harmonic curvature near an open
boundary patch.
For all smooth vector fields $X,Y,Z,U$,
\begin{equation}\label{ecriticalschoutencompatibility}
 W_g\bigl(X,Y,(\mathrm{Sch}_g)^\sharp Z,U\bigr)
 +W_g\bigl(Y,Z,(\mathrm{Sch}_g)^\sharp X,U\bigr)
 +W_g\bigl(Z,X,(\mathrm{Sch}_g)^\sharp Y,U\bigr)=0.
\end{equation}
Along this boundary patch, every smooth vector field $V$ satisfies
\[
 W_g(V)=O(r^\infty)
 \quad\Longrightarrow\quad
 W_g\bigl((\mathrm{Sch}_g)^\sharp V\bigr)=O(r^\infty).
\]
\end{lemma}

\begin{proof}
Harmonic curvature makes $\mathrm{Sch}_g$ Codazzi.
Apply $\nabla_i$ to
$(\nabla_k\mathrm{Sch}_g)_{\ell m}
 -(\nabla_\ell\mathrm{Sch}_g)_{km}=0$
and add the two equations obtained by cyclically permuting
$(i,k,\ell)$. Expanding gives
\begin{equation}\label{schoutenciclicsum}
 \bigl([\nabla_i,\nabla_k]\mathrm{Sch}_g\bigr)_{\ell m}
 +\bigl([\nabla_k,\nabla_\ell]\mathrm{Sch}_g\bigr)_{im}
 +\bigl([\nabla_\ell,\nabla_i]\mathrm{Sch}_g\bigr)_{km}=0.
\end{equation}
The commutation formula for the Schouten tensor is
\begin{equation}\label{schoutencommutation}
 \bigl([\nabla_i,\nabla_k]\mathrm{Sch}_g\bigr)_{\ell m}
 =-R_{ik\ell}{}^q(\mathrm{Sch}_g)_{qm}
  -R_{ikm}{}^q(\mathrm{Sch}_g)_{\ell q}.
\end{equation}
Substituting \eqref{schoutencommutation} into the cyclic sum
\eqref{schoutenciclicsum}, we see that the terms arising from
the first term on the right-hand side of \eqref{schoutencommutation} cancel by the first Bianchi identity.
Using the antisymmetry of the curvature tensor in its last two
indices for the remaining terms, we obtain
\[
 R_{ikqm}(\mathrm{Sch}_g)_\ell{}^q
 +R_{k\ell qm}(\mathrm{Sch}_g)_i{}^q
 +R_{\ell iqm}(\mathrm{Sch}_g)_k{}^q=0.
\]
This is \eqref{ecriticalschoutencompatibility} with $R_g$ in place
of $W_g$.

To replace $R_g$ by $W_g$, according to \eqref{defweyltensor}
 \begin{align*}
 (R_g-W_g)\bigl(X,Y,(\mathrm{Sch}_g)^\sharp Z,U\bigr)=&g\bigl((\mathrm{Sch}_g)^\sharp Y,(\mathrm{Sch}_g)^\sharp Z\bigr)g(X,U)
 \\&-g\bigl((\mathrm{Sch}_g)^\sharp X, (\mathrm{Sch}_g)^\sharp Z\bigr)g(Y,U)\\
 &+\mathrm{Sch}_g(Y,Z)\mathrm{Sch}_g(X,U)\\&-\mathrm{Sch}_g(X,Z)\mathrm{Sch}_g(Y,U).
 \end{align*}
The cyclic sum over $(X,Y,Z)$ vanishes: the first two terms
cancel because
$g((\mathrm{Sch}_g)^\sharp Y,(\mathrm{Sch}_g)^\sharp Z)$
is symmetric in $Y,Z$, and the last two cancel because
$\mathrm{Sch}_g$ is symmetric. This proves
\eqref{ecriticalschoutencompatibility}.

Finally, set $Z=V$ in that identity. Its second and third terms
are $O(r^\infty)$ because they contain an insertion of $V$.
Hence
$W_g(X,Y,(\mathrm{Sch}_g)^\sharp V,U)=O(r^\infty)$
for arbitrary $X,Y,U$. The stated conclusion follows from the
curvature symmetries.
\end{proof}

\begin{lemma}
Let $g$ be a smooth metric with harmonic curvature near an open
boundary patch. For any smooth vector fields $U,V$,
\begin{equation}\label{ecriticalnullityderivative}
 W_g(U)=O(r^\infty),\qquad W_g(V)=O(r^\infty)
 \quad\Longrightarrow\quad
 W_g(\nabla_UV)=O(r^\infty).
\end{equation}
\end{lemma}

\begin{proof}
Assume $W_g(U)=O(r^\infty)$ and $W_g(V)=O(r^\infty)$.
The condition $O(r^\infty)$ is preserved by differentiation along
smooth vector fields and by contraction with smooth tensors.
For arbitrary smooth vector fields $X,Y,Z$, by the definition of the covariant derivative of a
covariant tensor, we have
\[
 \begin{aligned}
 (\nabla_XW_g)(Y,U,V,Z)
 ={}&X\bigl[W_g(Y,U,V,Z)\bigr]\\
 &-W_g(\nabla_XY,U,V,Z)
  -W_g(Y,\nabla_XU,V,Z)\\
 &-W_g(Y,U,\nabla_XV,Z)
  -W_g(Y,U,V,\nabla_XZ).
 \end{aligned}
\]
All terms on the right-hand side are $O(r^\infty)$. Thus $(\nabla_XW_g)(Y,U,V,Z)=O(r^\infty)$.
The same argument shows that
$(\nabla_YW_g)(U,X,V,Z)=O(r^\infty)$.

By \eqref{eharmonicweylidentities}, we have
\[
 (\nabla_UW_g)(X,Y,V,Z)
 +(\nabla_XW_g)(Y,U,V,Z)
 +(\nabla_YW_g)(U,X,V,Z)=0.
\]
Consequently, $(\nabla_UW_g)(X,Y,V,Z)=O(r^\infty)$.
Differentiating $W_g(X,Y,V,Z)=O(r^\infty)$ along $U$ now yields
\[
 U\bigl[W_g(X,Y,V,Z)\bigr]
 =(\nabla_UW_g)(X,Y,V,Z)+W_g(X,Y,\nabla_UV,Z)+O(r^\infty).
\]
It follows that $W_g(X,Y,\nabla_UV,Z)=O(r^\infty)$ for every
$X,Y,Z$. The desired implication \eqref{ecriticalnullityderivative} follows from the curvature symmetries.
\end{proof}

\begin{lemma}\label{lemmacriticalnormalweyl}
Let $g=dr^2+h_r$ have harmonic curvature.
If $V|_{r=0}=\partial_r$ and $W_g(V)=O(r^\infty)$, then
$W_g(\partial_r)=O(r^\infty)$.
\end{lemma}
\begin{proof}
Starting from $V_0:=V$, define
\[
 V_{k+1}:=V_k-\frac r{k+1}\nabla_{V_k}V_k,\qquad k\geq0.
\]
By \eqref{ecriticalnullityderivative},
$W_g(V_k)=O(r^\infty)$ for every $k$.
If $V_k=\partial_r+r^{k+1}B_k$, then
$\nabla_{\partial_r}\partial_r=0$ gives
\[
 \nabla_{V_k}V_k=(k+1)r^kB_k+O(r^{k+1}).
\]
Thus $V_{k+1}-\partial_r=O(r^{k+2})$. For any fixed order,
choose $k$ sufficiently large and compare $W_g(V_k)$ with
$W_g(\partial_r)$. This proves the assertion at every order.
\end{proof}
\subsection{Conformal geodesic coordinates and frames}

We first choose coordinates and frames by ordinary differential
equations. Their compatibility with conformal changes will let us
compare the two metrics using the same boundary data. 
\begin{proposition}\label{lemmacriticalcoordinates}
Let $g=dr^2+h_r$ be smooth near a boundary patch.
Fix boundary coordinates $y$, a $g$-orthonormal frame $e_i^0$
along the patch with $e_n^0=\partial_r$, and a smooth one-form
$\eta_0$ along the patch that annihilates boundary tangent vectors.
Then the following hold locally.
\begin{enumerate}[label=(\roman*)]
\item
These data determine coordinates $(t,y)$ (not necessarily
orthogonal) with boundary $t=0$, and smooth extensions
$\eta,e_i$ of $\eta_0,e_i^0$, such that
$V:=e_n=\partial_t$ satisfies $V|_{t=0}=\partial_r$ and
\begin{equation}\label{ecriticalconformalgeodesics}
 \nabla^g_VV=-2\eta(V)V+|V|_g^2\eta^\sharp,\qquad
 \nabla^g_V\eta
 =\eta(V)\eta-\tfrac12|\eta|_g^2V^\flat
   +\mathrm{Sch}_g(V,\cdot).
\end{equation}
The frame satisfies
\begin{equation}\label{ecriticalframetransport}
 \nabla^g_Ve_i
 =-\eta(V)e_i-\eta(e_i)V+g(V,e_i)\eta^\sharp.
\end{equation}

\item
Let $\theta^i$ be the dual coframe and set
$\varphi:=\log|V|_g$. Then
\begin{equation}\label{ecriticalnormalizedmetric}
 g=e^{2\varphi}g_{\boldsymbol{\theta}},\qquad
 g_{\boldsymbol{\theta}}
 :=\sum_i\theta^i\otimes\theta^i,\qquad
 \varphi|_{t=0}=0.
\end{equation}

\item
If $u$ is a smooth function vanishing on the patch, the
construction for $e^{2u}g$ and $\eta_0-(du)|_{t=0}$,
using the same boundary coordinates and initial frame,
leaves $(t,y)$ and $e_i$ unchanged and replaces
$\eta,\varphi$ by $\eta-du,\varphi+u$.
In particular, $g_{\boldsymbol{\theta}}$ is unchanged.
\end{enumerate}
\end{proposition}
\begin{proof}
Extend $g$ smoothly across the patch as a Riemannian
metric. In local coordinates, \eqref{ecriticalconformalgeodesics},
together with $\dot x=V$, is a smooth first-order ODE system
for $(x,V,\eta)$. By the local existence and smooth dependence theorem,
there are unique solutions starting at each boundary point, with
\[
 V|_{t=0}=\partial_r,\qquad \eta|_{t=0}=\eta_0.
\]
Extend boundary coordinates $y$ constantly along these curves $x(t)$.
Since the initial velocity is transverse to the boundary,
the inverse function theorem gives coordinates $(t,y)$
on a smaller neighborhood, with $V=\partial_t$.

For each $i<n$, \eqref{ecriticalframetransport} is a linear
ODE for $e_i$ along the curves $x(t)$, with
smooth coefficients determined by $g,V,\eta$.
The prescribed boundary values therefore determine unique smooth
solutions, depending smoothly on the starting boundary point.
For $i=n$, the equation follows from the first equation in
\eqref{ecriticalconformalgeodesics}, since $e_n=V$.
Metric compatibility and \eqref{ecriticalframetransport} give
\[
 \partial_tg(e_i,e_j)=-2\eta(V)g(e_i,e_j),\qquad
 \partial_t|V|_g^2=-2\eta(V)|V|_g^2.
\]
Thus $g(e_i,e_j)/|V|_g^2$ is constant along each curve.
Its initial value is $\delta_{ij}$, and hence
\[
 g(e_i,e_j)=|V|_g^2\delta_{ij}.
\]
In particular, the fields $e_i$ form a frame.
By the definitions of $\varphi$ and the dual coframe $\boldsymbol{\theta}$,
this proves \eqref{ecriticalnormalizedmetric}.

The equations \eqref{ecriticalconformalgeodesics} are the
conformal geodesic equations. Under $g\mapsto e^{2u}g$, the conformal transformation formulas are
\begin{align}
 \nabla^{e^{2u}g}_XY
 &=\nabla^g_XY+du(X)Y+du(Y)X-g(X,Y)\nabla^gu,\label{econnectiontrans}\\
 \mathrm{Sch}_{e^{2u}g}
 &=\mathrm{Sch}_g-\nabla^gdu
   +du\otimes du-\tfrac12|du|_g^2g.\label{eschoutentrans}
\end{align}
Substitution, using the induced connection on one-forms,
shows that $(V,\eta-du,e_i)$ satisfies
\eqref{ecriticalconformalgeodesics}--\eqref{ecriticalframetransport}
for $e^{2u}g$. See also
\cite[Eqs.~(1)--(5)]{Tod2012} for the curve tangent
vector $V$ and the one-form $\eta$.

Since $u=0$ on the patch, the initial normal and frame
are unchanged, and $\eta_0-(du)|_{t=0}$ still annihilates
boundary tangent vectors. Thus $\eta-du$ and $e_i$ satisfy the required initial conditions.
With the same boundary coordinates, uniqueness shows that
$(t,y)$ and $e_i$ are unchanged.
Finally,
\[
 \log|V|_{e^{2u}g}=\log(e^u|V|_g)=\varphi+u.
\]
The dual coframe, and hence $g_{\boldsymbol{\theta}}$, is unchanged.
\end{proof}

\subsection{Connection and curvature matrices}

We collect the coframe, connection forms, and Schouten tensor in
a matrix so that their transport equations can be compared together.
For a metric $\widehat g$ with orthonormal coframe $\boldsymbol{\theta}:=(\theta^i)$
and dual frame $e_i$, define the matrix of connection one-forms by
\[
 (\boldsymbol{\gamma}_{\widehat g})^i{}_j(X):=\widehat g(\nabla^{\widehat g}_Xe_j,e_i).
\]
We use the following notation:
 \begin{align*}
 (\mathbf{Sch}_{\widehat g})_i(X)&:=\mathrm{Sch}_{\widehat g}(X,e_i),\\
 (\mathbf{C}_{\widehat g})_i(X,Y)&:=C_{\widehat g}(X,Y,e_i),\\
 (\mathbf{W}_{\widehat g})^i{}_j(X,Y)&:=W_{\widehat g}(X,Y,e_j,e_i).
 \end{align*}
Matrix products
use ordinary matrix multiplication with wedge products of the
form entries. The metric structure equations take the form (see also the tractor connection and its curvature in \cite[Eqs.~(5), (158)]{GoverShaukatWaldron2009})
\begin{equation}\label{ecriticalmetricmatrix}
 \mathbf{B}_{\widehat g}:=
 \begin{pmatrix}
 0&-\boldsymbol{\theta}^{\mathsf T}&0\\
 \mathbf{Sch}_{\widehat g}&\boldsymbol{\gamma}_{\widehat g}&\boldsymbol{\theta}\\
 0&-\mathbf{Sch}_{\widehat g}^{\mathsf T}&0
 \end{pmatrix},\qquad  d\mathbf{B}_{\widehat g}+\mathbf{B}_{\widehat g}\wedge \mathbf{B}_{\widehat g}=
 \begin{pmatrix}
 0&0&0\\\mathbf{C}_{\widehat g}&\mathbf{W}_{\widehat g}&0\\0&-\mathbf{C}_{\widehat g}^{\mathsf T}&0
 \end{pmatrix}.
\end{equation}
For any column $\mathbf{b}$, set
\begin{equation}\label{ecriticalchangeframe}
 \mathbf{U}(\mathbf{b}):=
 \begin{pmatrix}
 1&0&0\\\mathbf{b}&\Id&0\\-\tfrac12|\mathbf{b}|^2&-\mathbf{b}^{\mathsf T}&1
 \end{pmatrix}.
\end{equation}

We now apply this construction to the smooth metric
$g=dr^2+h_r$ near an open boundary patch, using the conformal
geodesic coordinates of Proposition~\ref{lemmacriticalcoordinates}.
Fix a smooth one-form $\eta_0$ along the patch annihilating
boundary tangent vectors, and a $g$-orthonormal boundary
frame with $e_n^0=\partial_r$.
Apply Proposition~\ref{lemmacriticalcoordinates} with these initial data.
Define
\begin{equation}\label{ecriticalphysicalb}
 b_i:=\eta(e_i),\qquad \beta_i:=-b_i-e_i\varphi,\qquad
 \mathbf{B}:=\mathbf{U}(\boldsymbol{\beta})^{-1}\mathbf{B}_{g_{\boldsymbol{\theta}}}\mathbf{U}(\boldsymbol{\beta})  +\mathbf{U}(\boldsymbol{\beta})^{-1}d\mathbf{U}(\boldsymbol{\beta}).
\end{equation}
Here $\mathbf{b}:=(b_i),\boldsymbol{\beta}:=(\beta_i)$ are columns of functions with Euclidean norms. Define the curvature matrix 
$\boldsymbol{\Omega}:=d\mathbf{B}+\mathbf{B}\wedge \mathbf{B}$. Since $d(\mathbf{U}^{-1})=-\mathbf{U}^{-1}(d\mathbf{U})\mathbf{U}^{-1}$, it follows that
\begin{equation*}
 \boldsymbol\Omega
 =\mathbf U(\boldsymbol\beta)^{-1}
  \bigl(d\mathbf B_{g_{\boldsymbol\theta}}
       +\mathbf B_{g_{\boldsymbol\theta}}
 \wedge\mathbf B_{g_{\boldsymbol\theta}}\bigr)
  \mathbf U(\boldsymbol\beta).
\end{equation*}
Using \eqref{ecriticalmetricmatrix} and
$\mathbf U(\boldsymbol{\beta})^{-1}=\mathbf U(-\boldsymbol{\beta})$,
 we have
\begin{equation}\label{ecriticalconnectionmatrix}
 \boldsymbol{\Omega}=
 \begin{pmatrix}
 0&0&0\\
 \mathbf{C}_{g_{\boldsymbol{\theta}}}+\mathbf{W}_{g_{\boldsymbol{\theta}}}\boldsymbol{\beta}&\mathbf{W}_{g_{\boldsymbol{\theta}}}&0\\
 0&-(\mathbf{C}_{g_{\boldsymbol{\theta}}}+\mathbf{W}_{g_{\boldsymbol{\theta}}}\boldsymbol{\beta})^{\mathsf T}&0
 \end{pmatrix}.
\end{equation}

Now we calculate $\mathbf{B}(V)$. Block multiplication in \eqref{ecriticalphysicalb} gives
\begin{equation}\label{ecriticalmetricconnection}
 \mathbf B=
 \begin{pmatrix}
 -\boldsymbol{\beta}^{\mathsf T}\boldsymbol{\theta}&-\boldsymbol{\theta}^{\mathsf T}&0\\
 \mathbf B_{21}&
 \boldsymbol\gamma_{g_{\boldsymbol{\theta}}}
   +\boldsymbol{\beta}\boldsymbol{\theta}^{\mathsf T}-\boldsymbol{\theta}\boldsymbol{\beta}^{\mathsf T}
 &\boldsymbol{\theta}\\
 0&-\mathbf B_{21}^{\mathsf T}&\boldsymbol{\beta}^{\mathsf T}\boldsymbol{\theta}
 \end{pmatrix},
\end{equation}
where
\begin{equation}\label{defb21block}
 \mathbf B_{21}
 =\mathbf{Sch}_{g_{\boldsymbol{\theta}}}+d\boldsymbol{\beta}
  +\boldsymbol\gamma_{g_{\boldsymbol{\theta}}}\boldsymbol{\beta}
  +\boldsymbol{\beta}(\boldsymbol{\beta}^{\mathsf T}\boldsymbol{\theta})-\tfrac12|\boldsymbol{\beta}|^2\boldsymbol{\theta}.
\end{equation}
Here and below block indices refer to the displayed three-by-three
block decomposition.

Since $V=\partial_t$, metric compatibility and
\eqref{ecriticalconformalgeodesics} give
\[
 \partial_t|V|_g^2=2g(\nabla^g_VV,V)=-4\eta(V)|V|_g^2
   +2|V|_g^2g(\eta^\sharp,V)=-2b_n|V|_g^2,
\]
where $g(\eta^\sharp,V)=\eta(V)=b_n$.
By $\varphi=\log|V|_g$ and $e_n=\partial_t$, it follows that
\begin{equation}\label{ephitbettn}
     \partial_t\varphi=\frac{\partial_t|V|_g^2}{2|V|_g^2}=-b_n,\qquad
 \beta_n=-b_n-e_n(\varphi)=0.
\end{equation}
Substituting into \eqref{ecriticalmetricconnection}, one gets $\mathbf{B}_{11}(V)=0$.

Moreover, \eqref{ecriticalconformalgeodesics}
and \eqref{ecriticalframetransport} yield
 \begin{align}
 \partial_tb_i &=(\nabla^g_V\eta)(e_i)+\eta(\nabla^g_Ve_i)\notag\\
 &=\mathrm{Sch}_g(V,e_i)-b_nb_i+\tfrac12|\eta|_g^2g(V,e_i)\notag\\
 &=\mathrm{Sch}_g(V,e_i)-b_nb_i  +\tfrac12|\mathbf{b}|^2\delta_{in}.\label{ecriticalbtransport}
 \end{align}
The last equality uses
$g(V,e_i)=e^{2\varphi}\delta_{in}$ and
$|\eta|_g^2=e^{-2\varphi}|\mathbf{b}|^2$, which follow from \eqref{ecriticalnormalizedmetric}.
Using \eqref{econnectiontrans}, \eqref{ecriticalframetransport}, and \eqref{ephitbettn}, we compute
 \begin{align*}
 \nabla^{g_{\boldsymbol{\theta}}}_Ve_j&=\nabla^g_Ve_j-(\partial_t\varphi)e_j-(e_j\varphi)V+g(V,e_j)\nabla^g\varphi\\
 &=\beta_je_n-\delta_{jn}\sum_i\beta_i e_i.
 \end{align*}
Thus $\boldsymbol{\gamma}_{g_{\boldsymbol{\theta}}}(V)=-\boldsymbol{\beta}\mathbf e_n^{\mathsf T}+\mathbf e_n\boldsymbol{\beta}^{\mathsf T},$ where $\mathbf e_n:=(0,\dots,0,1)^{\mathsf T}$. Substituting into \eqref{ecriticalmetricconnection}, we have $\mathbf{B}_{22}(V)=0$.

By \eqref{defb21block},
 \begin{align}
 \mathbf B_{21}(V) =&\mathbf{Sch}_{g_{\boldsymbol{\theta}}}(V)+\partial_t\boldsymbol{\beta}+\boldsymbol\gamma_{g_{\boldsymbol{\theta}}}(V)\boldsymbol{\beta}
   +\boldsymbol{\beta}\beta_n-\tfrac12|\boldsymbol{\beta}|^2\mathbf e_n\notag\\
 =&\mathbf{Sch}_{g_{\boldsymbol{\theta}}}(V)+\partial_t\boldsymbol{\beta}+\tfrac12|\boldsymbol{\beta}|^2\mathbf e_n,\label{eqblockexpressionB21}
 \end{align}
where we used $\boldsymbol{\theta}(V)=\mathbf e_n$ and
$\boldsymbol\gamma_{g_{\boldsymbol{\theta}}}(V)\boldsymbol{\beta}
=-\boldsymbol{\beta}\beta_n+|\boldsymbol{\beta}|^2\mathbf e_n$.

Using \eqref{ecriticalframetransport} and
$\partial_t\varphi=-b_n$, we compute
\begin{align*}
 (\nabla^g d\varphi)(V,e_i)
 &=\partial_t(e_i\varphi)-d\varphi(\nabla^g_Ve_i)\\
 &=\partial_t(e_i\varphi)+b_ne_i\varphi-b_nb_i
   -\delta_{in}\sum_j b_je_j\varphi.
\end{align*}
Using $g(V,e_i)=e^{2\varphi}\delta_{in}$ and
$|d\varphi|_g^2=e^{-2\varphi}\sum_j(e_j\varphi)^2$,
we deduce from the Schouten transformation formula \eqref{eschoutentrans} that
 \begin{align*}
 \mathrm{Sch}_{g_{\boldsymbol{\theta}}}(V,e_i)
 =&\mathrm{Sch}_g(V,e_i)+\partial_t(e_i\varphi)-b_nb_i-\delta_{in}\sum_j
       \bigl(b_je_j\varphi+\tfrac12(e_j\varphi)^2\bigr)\\
 =&\mathrm{Sch}_g(V,e_i)+\partial_t(e_i\varphi)-b_nb_i
   +\tfrac12\bigl(|\mathbf{b}|^2-|\boldsymbol{\beta}|^2\bigr)\delta_{in},
 \end{align*}
where the last equality follows from
$\beta_j=-b_j-e_j\varphi$.
Substituting this expression and
$\partial_t\beta_i=-\partial_tb_i-\partial_t(e_i\varphi)$
into the block expression \eqref{eqblockexpressionB21} gives
\[
 \mathbf B_{21}(V)=\mathbf{Sch}_{g}(V)-\partial_t\mathbf{b}-b_n\mathbf{b}
   +\tfrac12|\mathbf{b}|^2\mathbf e_n=0,
\]
with the final equality following from
\eqref{ecriticalbtransport}.

Consequently,
\begin{equation}\label{ecriticalradialvalues}
 \mathbf{B}(V)=
 \begin{pmatrix}
 0&-\mathbf e_n^{\mathsf T}&0\\
 0&0&\mathbf e_n\\0&0&0
 \end{pmatrix}=:\mathbf{X}_0.
\end{equation}

\subsection{Reconstruction of the connection matrix}

If $g=dr^2+h_r$ has harmonic curvature, its Cotton tensor $C_g=0$. By \eqref{cottontransformfor}, we then have
$\mathbf{C}_{g_{\boldsymbol{\theta}}}=\mathbf{W}_{g_{\boldsymbol{\theta}}}(e_i\varphi)_i$. Thus
\begin{equation}\label{ecriticalphysicalcotton}
 \mathbf{C}_{g_{\boldsymbol{\theta}}}+\mathbf{W}_{g_{\boldsymbol{\theta}}}\boldsymbol{\beta}=-\mathbf{W}_{g_{\boldsymbol{\theta}}}\mathbf{b}.
\end{equation}

Apply Proposition~\ref{lemmacriticalcoordinates} with an initial
one-form $\eta_0$ annihilating boundary tangent vectors and define $\mathbf B$ by \eqref{ecriticalphysicalb}. Let $\mathbf B_0$ be the pullback of $\mathbf B$ to the boundary $t=0$. Extend $\mathbf B_0$ independently of $t$. 
The next proposition shows that the tangential boundary values
of $\mathbf B$ determine the full connection matrix near the patch.

\begin{proposition}\label{propcriticalradialmodel}
Let $g=dr^2+h_r$ have harmonic curvature near an open
boundary patch, and suppose that
$W_g(\partial_r)=O(r^\infty)$.
Then $\mathbf B_0$ uniquely determines $\mathbf B$
on a smaller neighborhood through the formula
\begin{equation}\label{ecriticalradialmodel}
 \mathbf B
 =\mathbf X_0\,dt
  +e^{-t\mathbf X_0}\mathbf B_0e^{t\mathbf X_0},
\end{equation}
where $\mathbf X_0$ is defined in
\eqref{ecriticalradialvalues}.
\end{proposition}

We divide the proof into three auxiliary lemmas. Define
\begin{equation}\label{ecriticalreferenceconnection}
 \overline{\mathbf B}
 :=\mathbf X_0\,dt
   +e^{-t\mathbf X_0}\mathbf B_0e^{t\mathbf X_0},
 \qquad \mathbf A:=\mathbf B-\overline{\mathbf B}.
\end{equation}
Let $\boldsymbol\Omega_0$ be the pullback of $\boldsymbol{\Omega}$ to the boundary, namely, $\boldsymbol\Omega_0:=d_Y\mathbf B_0+\mathbf B_0\wedge\mathbf B_0$.
Extend $\boldsymbol\Omega_0$ independently of $t$. We use these data and definitions in the
following three lemmas to prove
$\mathbf B=\overline{\mathbf B}$ near the boundary.

\begin{lemma}\label{lemmacriticalreferencejets} Under the hypothesis of Proposition~\ref{propcriticalradialmodel}, we have
\begin{equation*}
 \mathbf A=O(t^\infty),\qquad
 \iota_V\boldsymbol\Omega=O(t^\infty).
\end{equation*}
The curvature matrix associated with $\overline{\mathbf{B}}$ satisfies
\begin{equation}\label{ecriticalreferencecurvature}
 d\overline{\mathbf B}
   +\overline{\mathbf B}\wedge\overline{\mathbf B}
 =\boldsymbol\Omega_0=\begin{pmatrix}
 0&0&0\\0&\mathbf W_g|_{TY}&0\\0&0&0
 \end{pmatrix},\qquad
 [\mathbf X_0,\boldsymbol\Omega_0]=0.
\end{equation}
\end{lemma}

\begin{proof}
  We first prove
\begin{equation}\label{asseretionWVeta}
 W_g(V)=O(r^\infty),\qquad
 W_g(\eta^\sharp)=O(r^\infty).
\end{equation}
Since $V|_{r=0}=\partial_r$ and $\eta^\sharp|_{r=0}$ is normal,
the hypothesis $W_g(\partial_r)=O(r^\infty)$ gives
\[
 W_g(V)=O(r),\qquad W_g(\eta^\sharp)=O(r).
\]
Suppose inductively that, for some $k\geq1$,
\[
 W_g(V)=O(r^k),\qquad W_g(\eta^\sharp)=O(r^k).
\]
For $U\in\{V,\eta^\sharp\}$ and arbitrary smooth vector fields
$X,Y,Z$, we have $V-\partial_r=O(r)$ and hence
\[
 W_g(Y,V,U,Z)
 =W_g(Y,V-\partial_r,U,Z)+O(r^\infty)
 =O(r^{k+1}).
\]
Differentiating shows
\[
 (\nabla^g_XW_g)(Y,V,U,Z)=O(r^k),\qquad(\nabla^g_YW_g)(V,X,U,Z)=O(r^k).
\]
The first identity in \eqref{eharmonicweylidentities} therefore yields
\[
 (\nabla^g_VW_g)(X,Y,U,Z)=-(\nabla^g_XW_g)(Y,V,U,Z)-(\nabla^g_YW_g)(V,X,U,Z)=O(r^k).
\]

Taking $Z=V$ in \eqref{ecriticalschoutencompatibility}, we find that
\[
 W_g\bigl((\mathrm{Sch}_g)^\sharp V\bigr)=O(r^k).
\]
Consequently, \eqref{ecriticalconformalgeodesics} and the
curvature symmetries imply, for $U\in\{V,\eta^\sharp\}$,
\[
 \nabla^g_V\bigl(W_g(U)\bigr)=(\nabla^g_VW_g)(U,\cdot,\cdot,\cdot)
   +W_g(\nabla^g_VU) =O(r^k).
\]
Since $V=\partial_t$ and $r\asymp t$, integration in $t$ with
$W_g(V)|_{t=0}=W_g(\eta^\sharp)|_{t=0}=0$ gives
\[
 W_g(V)=O(r^{k+1}),\qquad
 W_g(\eta^\sharp)=O(r^{k+1}).
\]
This completes the induction and proves \eqref{asseretionWVeta}.
Since $g=e^{2\varphi}g_{\boldsymbol{\theta}}$, equations \eqref{ecriticalconnectionmatrix}, \eqref{ecriticalphysicalcotton} and \eqref{asseretionWVeta} therefore give
$\iota_V\boldsymbol\Omega=O(t^\infty)$.

At $t=0$, $b_i=0$ for $i<n$ and
$\mathbf W_{g_{\boldsymbol{\theta}}}\mathbf e_n=0$.
By \eqref{ecriticalphysicalcotton}, the $(2,1)$ block of $\boldsymbol\Omega_0$ is zero. Its middle block is
$\mathbf W_g|_{TY}$, giving the block form in
\eqref{ecriticalreferencecurvature}.
The equality $\mathbf W_g\mathbf e_n=0$ at the boundary
also implies $[\mathbf X_0,\boldsymbol\Omega_0]=0$.
Since $\partial_t\overline{\mathbf B}
=-[\mathbf X_0,\overline{\mathbf B}]$, the $dt$ terms
in the curvature matrix cancel. Thus
\begin{align*}
 d\overline{\mathbf B}+\overline{\mathbf B}\wedge\overline{\mathbf B}
 &=e^{-t\mathbf X_0}(d_Y\mathbf B_0+\mathbf B_0\wedge\mathbf B_0)e^{t\mathbf X_0}\\
 &=e^{-t\mathbf X_0}\boldsymbol\Omega_0e^{t\mathbf X_0}
 =\boldsymbol\Omega_0.
\end{align*}

It remains to prove $\mathbf A=O(t^\infty)$. 
Subtracting the definitions of $\boldsymbol{\Omega},\boldsymbol{\Omega}_0$, we have
\begin{equation}\label{diffOmegaAB}
\boldsymbol{\Omega}-\boldsymbol{\Omega}_0=d\mathbf{A}+\overline{\mathbf{B}}\wedge\mathbf{A}+\mathbf{A}\wedge\overline{\mathbf{B}}+\mathbf{A}\wedge\mathbf{A}.
\end{equation}
Contracting with
$V=\partial_t$, using
$\mathbf B(V)=\overline{\mathbf B}(V)=\mathbf X_0$ and $\iota_V\boldsymbol\Omega_0=0$, we obtain
\begin{equation}\label{ecriticalconnectiontransport}
 \mathbf A(V)=0,\qquad\partial_t\mathbf A=\iota_V\boldsymbol\Omega-[\mathbf X_0,\mathbf A].
\end{equation}
Since $\mathbf A|_{t=0}=0$ and
$\iota_V\boldsymbol\Omega=O(t^\infty)$, successive differentiation shows that $\mathbf A=O(t^\infty)$.
\end{proof}

To compare Weyl and Cotton tensors, we next construct a
reference metric $\overline g$ from $\overline{\mathbf B}$. 
Denote by $\overline{\boldsymbol{\theta}}:=(\overline{\theta}^i)$ the $(2,3)$ block of $\overline{\mathbf B}$.
 Set $\overline g:=\sum_i\overline\theta^i\otimes\overline\theta^i$. By Lemma~\ref{lemmacriticalreferencejets},
$\overline{\mathbf B}-\mathbf B=O(t^\infty)$.
In particular, $\overline{\boldsymbol{\theta}}-\boldsymbol{\theta}=O(t^\infty)$. 
Thus $\overline{\boldsymbol{\theta}}$ is a coframe near $t=0$. Let $\overline e_i$ be dual to $\overline{\boldsymbol{\theta}}$, and define $\overline{\boldsymbol{\beta}}$ by $\overline{\mathbf B}_{11}=-\overline{\boldsymbol{\beta}}^{\mathsf T}\overline{\boldsymbol{\theta}}$. It follows that 
 \begin{equation}\label{matchingjetsgebeta}
 \overline g-g_{\boldsymbol{\theta}}=O(t^\infty),\qquad
 \overline e_i-e_i=O(t^\infty),\qquad
 \overline{\boldsymbol{\beta}}-\boldsymbol{\beta}=O(t^\infty).
\end{equation}

For the reference metric $\overline g$ constructed above,
$\mathbf W_{\overline g}$, $\mathbf B_{\overline g}$, and the
other matrices associated with $\overline g$ are defined using the coframe
$\overline{\boldsymbol\theta}$ and its dual frame $\overline e_i$.

\begin{lemma}
  Under the hypothesis of Proposition~\ref{propcriticalradialmodel}, on a smaller neighborhood of the boundary patch, we have
\begin{equation}\label{Bbarrelationgbar}
 \overline{\mathbf B}=\mathbf U(\overline{\boldsymbol{\beta}})^{-1}\mathbf B_{\overline g}\mathbf U(\overline{\boldsymbol{\beta}})+\mathbf U(\overline{\boldsymbol{\beta}})^{-1}d\mathbf U(\overline{\boldsymbol{\beta}}).
\end{equation}
Moreover, $\overline e_n=V$ and the following equations hold:
\begin{equation}\label{ecriticalmodelweylconstant}
 \partial_t\mathbf W_{\overline g}=0,\qquad  \iota_V\mathbf W_{\overline g}=0.
\end{equation}
\begin{equation}\label{ecriticalmodelcotton}
  \mathbf C_{\overline g}+\mathbf W_{\overline g}\overline{\boldsymbol\beta}=0.
\end{equation}
\end{lemma}

\begin{proof}
Since $\mathbf X_0^3=0$,
\[
 e^{\pm t\mathbf X_0}=\Id\pm t\mathbf X_0+\tfrac12t^2\mathbf X_0^2.
\]
Thus \eqref{ecriticalreferenceconnection} makes the coefficients of
$\overline{\mathbf B}$ polynomial in $t$.
The coefficients of $\overline{\boldsymbol\theta}$, its dual frame
$\overline e_i$, $\overline g$, $\overline{\boldsymbol\beta}$,
and $\mathbf B_{\overline g}$ are therefore analytic in $t$,
with smooth dependence on $y$.
The matching jets \eqref{matchingjetsgebeta} and \eqref{ecriticalphysicalb} give
\[
\mathbf U(\overline{\boldsymbol{\beta}})^{-1}\mathbf B_{\overline g}
 \mathbf U(\overline{\boldsymbol{\beta}})+\mathbf U(\overline{\boldsymbol{\beta}})^{-1}d\mathbf U(\overline{\boldsymbol{\beta}})
-\overline{\mathbf B} =O(t^\infty).
\]
Every coefficient on the left is analytic in $t$ and has all
$t$-derivatives zero at $t=0$. Hence the left-hand side vanishes
near the patch, proving \eqref{Bbarrelationgbar}.

Substituting \eqref{Bbarrelationgbar} into
\eqref{ecriticalreferencecurvature}, and using
$d(\mathbf U^{-1})=-\mathbf U^{-1}(d\mathbf U)\mathbf U^{-1}$,
we obtain
\begin{align*}
 \boldsymbol\Omega_0&=d\overline{\mathbf B}
   +\overline{\mathbf B}\wedge\overline{\mathbf B}=\mathbf U(\overline{\boldsymbol\beta})^{-1}
   \bigl(d\mathbf B_{\overline g}
        +\mathbf B_{\overline g}\wedge\mathbf B_{\overline g}\bigr)
   \mathbf U(\overline{\boldsymbol\beta})\\
 &=\begin{pmatrix}
 0&0&0\\
 \mathbf C_{\overline g}
   +\mathbf W_{\overline g}\overline{\boldsymbol\beta}
   &\mathbf W_{\overline g}&0\\
 0&-\bigl(\mathbf C_{\overline g}
  +\mathbf W_{\overline g}\overline{\boldsymbol\beta}\bigr)^{\mathsf T}
 &0
 \end{pmatrix}.
\end{align*}
The last equality follows from \eqref{ecriticalmetricmatrix}
and multiplication by the matrices in
\eqref{ecriticalchangeframe}.
Comparing with \eqref{ecriticalreferencecurvature}, we find
\[
 (\boldsymbol\Omega_0)_{22}=\mathbf W_{\overline g},\qquad
 \mathbf C_{\overline g}
 =-\mathbf W_{\overline g}\overline{\boldsymbol\beta}.
\]
Also,
\[
 \overline{\mathbf B}(V)=\mathbf X_0
 \quad\Longrightarrow\quad
 \overline{\boldsymbol{\theta}}(V)=\mathbf e_n
 \quad\Longrightarrow\quad
 \overline e_n=V.
\]
Since $\boldsymbol\Omega_0$ is independent of $t$ and has no $dt$
component, the same is true of $\mathbf W_{\overline g}$.
Consequently,
\[
 \partial_t(\mathbf W_{\overline g})=0,\qquad
 \iota_V\mathbf W_{\overline g}=0.\qedhere
\]
\end{proof}

Define a $(0,4)$-tensor $\mathbf w$, with components taken in
the frame $e_i$, and a column $\mathbf z$ of two-forms by
\[
 w_{ik\ell m}
 :=W_{g_{\boldsymbol{\theta}}}(e_i,e_k,e_\ell,e_m)
   -W_{\overline g}(\overline e_i,\overline e_k,
                                  \overline e_\ell,\overline e_m),
\]
\begin{equation*}
 \mathbf{z}:=\mathbf W_{\overline g}\mathbf{b},\qquad
 z_{i\mu\nu}
 =\sum_j b_jW_{\overline g}(\partial_\mu,\partial_\nu,
                                      \overline e_j,\overline e_i).
\end{equation*}
Here $b_j=\eta(e_j)$, all frame indices range from $1$ to $n$,
and $\mu,\nu$ range over $(t,y)$.

We next combine the Weyl equations with transport equations
for $\mathbf A$ and $\mathbf z$ to obtain the inequalities
required by Lemma~\ref{lemmacoupledcurvatureuc}.

\begin{lemma}\label{lemmacriticaldifferencesystem}
  Set
\[
 \mathbf v:=(\mathbf A,\mathbf{z},\partial_y\mathbf A,\partial_y\mathbf{z}),
 \qquad L:=g_{\boldsymbol{\theta}}^{\mu\nu}\partial_\mu\partial_\nu.
\]
Under the hypothesis of Proposition~\ref{propcriticalradialmodel}, on a smaller neighborhood, the arrays $\mathbf w,\mathbf v$ satisfy
\begin{equation}\label{ecriticalclosedsystem}
 |L\mathbf{w}|\leq C(|\mathbf{w}|+|\partial \mathbf{w}|+|\mathbf v|),\qquad
 |\partial_t\mathbf v|\leq C(|\mathbf v|+|\mathbf{w}|+|\partial \mathbf{w}|).
\end{equation}
Here we use Euclidean norms. The symbols $\partial$
and $\partial_y$ denote all first coordinate derivatives and
the tangential ones, respectively.
\end{lemma}
\begin{proof}
Shrink to a neighborhood with compact closure, so that all
coefficient fields and derivatives used below are bounded.
By the definition of $\mathbf w$,
\[
 (\mathbf W_{g_{\boldsymbol{\theta}}}-\mathbf W_{\overline g})^k{}_{\ell}
 =\tfrac12 w_{ij\ell k}\theta^i\wedge\theta^j+\tfrac12 W_{\overline g}(\overline e_i,\overline e_j,\overline e_\ell,\overline e_k) (\theta^i\wedge\theta^j
-\overline\theta^i\wedge\overline\theta^j).
\]
Since
$\boldsymbol\theta-\overline{\boldsymbol\theta}=\mathbf A_{23}$,
this difference is linear in $\mathbf A,\mathbf w$, with smooth
coefficients. In particular,
\begin{equation}\label{ecriticalweylarraycomparison}
 |\mathbf W_{g_{\boldsymbol\theta}}-\mathbf W_{\overline g}|
 \leq C(|\mathbf A|+|\mathbf w|).
\end{equation}

We first control $\partial_t\mathbf{A}$. 
By \eqref{ecriticalmodelweylconstant},
we have $\iota_V\mathbf W_{\overline g}=0$.
By \eqref{ecriticalconnectionmatrix} and
\eqref{ecriticalphysicalcotton}, the $(2,2)$ and $(2,1)$ blocks of $\iota_V\boldsymbol\Omega$ are
$\iota_V\mathbf W_{g_{\boldsymbol\theta}}$ and
$-(\iota_V\mathbf W_{g_{\boldsymbol\theta}})\mathbf b$. Hence \eqref{ecriticalconnectiontransport} shows that
$\partial_t\mathbf A$ is linear in $\mathbf A,\mathbf w$,
with smooth coefficients.
To control $\mathbf z$, use
$W_g=e^{2\varphi}W_{g_{\boldsymbol\theta}}$ and
$g(e_i,e_j)=e^{2\varphi}\delta_{ij}$ in
\eqref{ecriticalschoutencompatibility}. We have
\begin{align*}
 W_g\bigl(X,Y,(\mathrm{Sch}_g)^\sharp V,e_i\bigr)&=\sum_j\mathrm{Sch}_g(V,e_j)W_{g_{\boldsymbol{\theta}}}(X,Y,e_j,e_i)\\
 &=\bigl[ \mathbf W_{g_{\boldsymbol{\theta}}}(X,Y) \mathbf{Sch}_g(V)\bigr]_i.
\end{align*}
Apply the same formula to the other two terms in
\eqref{ecriticalschoutencompatibility}.
Together with $\iota_V\mathbf W_{\overline g}=0$, this gives
\begin{align*}
 \mathbf W_{\overline g}(X,Y)\mathbf{Sch}_g(V)
 =&-(\mathbf W_{g_{\boldsymbol\theta}}-\mathbf W_{\overline g})(X,Y)\mathbf{Sch}_g(V)\\
 &-(\mathbf W_{g_{\boldsymbol\theta}}-\mathbf W_{\overline g}) (Y,V)\mathbf{Sch}_g(X)\\
 &-(\mathbf W_{g_{\boldsymbol\theta}}-\mathbf W_{\overline g})(V,X)\mathbf{Sch}_g(Y).
\end{align*}
Thus $\mathbf W_{\overline g}\mathbf{Sch}_g(V)$ is also
linear in $\mathbf A,\mathbf w$, with smooth coefficients.
Differentiating $\mathbf z=\mathbf W_{\overline g}\mathbf b$
with respect to $t$ and using
\eqref{ecriticalmodelweylconstant} and \eqref{ecriticalbtransport},
together with $\mathbf W_{\overline g}\mathbf e_n=0$, we obtain
\[
 \partial_t\mathbf z=\mathbf W_{\overline g}\mathbf{Sch}_g(V)-b_n\mathbf z.
\]
Consequently, the two transport equations take the form
\begin{equation}\label{ecriticalexacttransport}
 \partial_t\mathbf A =\mathbf T_1\mathbf A+\mathbf T_2\mathbf w,\qquad
 \partial_t\mathbf z =\mathbf T_3\mathbf A+\mathbf T_4\mathbf w-b_n\mathbf z,
\end{equation}
where the coefficient matrices
$\mathbf T_1,\ldots,\mathbf T_4$ are smooth with bounded
first derivatives.

The block $\mathbf B_{11}=-\boldsymbol\beta^{\mathsf T}\boldsymbol\theta$
determines $\boldsymbol\beta$ once $\boldsymbol\theta$ is given.
The same statements hold for $\overline{\mathbf{B}}$. Thus $\boldsymbol\theta-\overline{\boldsymbol\theta}$,
$\boldsymbol\beta-\overline{\boldsymbol\beta}$, and
$\boldsymbol\gamma_{g_{\boldsymbol\theta}}
 -\boldsymbol\gamma_{\overline g}$
are linear combinations of the entries of $\mathbf A$.
The connection coefficients in the respective orthonormal frames are
\[
 \gamma_{si}^{j}:=g_{\boldsymbol\theta}(\nabla^{g_{\boldsymbol\theta}}_{e_s}e_i,e_j)
  =(\boldsymbol\gamma_{g_{\boldsymbol\theta}})^j{}_i(e_s),
 \qquad \overline\gamma_{si}^{j}:=\overline g (\nabla^{\overline g}_{\overline e_s}\overline e_i,\overline e_j) =(\boldsymbol\gamma_{\overline g})^j{}_i(\overline e_s).
\]
Since
\begin{align*}
 e_s-\overline e_s &=-\sum_j(\theta^j-\overline\theta^j)(\overline e_s)e_j,\\
 \gamma_{si}^{j}-\overline\gamma_{si}^{j}&=(\boldsymbol\gamma_{g_{\boldsymbol\theta}} -\boldsymbol\gamma_{\overline g})^j{}_i(e_s)
   +(\boldsymbol\gamma_{\overline g})^j{}_i(e_s-\overline e_s),
\end{align*}
the differences $e_s-\overline e_s$ and 
$\gamma_{si}^{j}-\overline\gamma_{si}^{j}$ are smooth linear combinations of the entries of $\mathbf A$, with bounded coefficients and first derivatives.

The component formula for covariant
differentiation gives
\begin{align}
 (\nabla^{g_{\boldsymbol\theta}}W_{g_{\boldsymbol\theta}})_{siqk\ell}
  -(\nabla^{\overline g}W_{\overline g})_{siqk\ell}
 &=(\nabla^{g_{\boldsymbol\theta}}\mathbf w)_{siqk\ell}
  +(e_s-\overline e_s)\bigl((W_{\overline g})_{iqk\ell}\bigr)\notag\\
 &\quad-\sum_j\Bigl[(\gamma_{si}^{j}-\overline\gamma_{si}^{j})(W_{\overline g})_{jqk\ell}
 +(\gamma_{sq}^{j}-\overline\gamma_{sq}^{j}) (W_{\overline g})_{ijk\ell}\notag\\
 &\quad+(\gamma_{sk}^{j}-\overline\gamma_{sk}^{j}) (W_{\overline g})_{iqj\ell}+(\gamma_{s\ell}^{j}-\overline\gamma_{s\ell}^{j}) (W_{\overline g})_{iqkj}\Bigr].\label{diffcoWframe}
\end{align}
Thus, after subtracting $(\nabla^{g_{\boldsymbol\theta}}\mathbf w)_{siqk\ell}$,
the remainder is a smooth linear combination of the entries
of $\mathbf A$.

Likewise, by \eqref{ecriticalphysicalcotton} and
\eqref{ecriticalmodelcotton}, we have
\begin{equation}\label{ecriticalcottondifference}
 \mathbf C_{g_{\boldsymbol{\theta}}}-\mathbf C_{\overline g}
 =-\mathbf{z}-(\mathbf W_{g_{\boldsymbol{\theta}}}-\mathbf W_{\overline g})(\mathbf{b}+\boldsymbol{\beta})-\mathbf W_{\overline g}(\boldsymbol{\beta}-\overline{\boldsymbol{\beta}}).
\end{equation}
This expresses the Cotton difference as a smooth linear expression in
$\mathbf A,\mathbf w,\mathbf z$.
Passing to components in the respective frames adds only terms
containing a frame difference, hence terms linear in $\mathbf A$.
Taking cyclic sums and contractions in
\eqref{diffcoWframe}, we use
\eqref{ecriticalweylcottonequations} and
\eqref{ecriticalcottondifference} to express every component of
$d^{\nabla^{g_{\boldsymbol\theta}}}\mathbf w$ and
$\delta^{\nabla^{g_{\boldsymbol\theta}}}\mathbf w$
as a smooth linear combination of the entries of
$\mathbf A,\mathbf w,\mathbf z$.
Differentiating these component identities yields
\begin{align}
 &|d^{\nabla^{g_{\boldsymbol\theta}}}\mathbf w|
  +|\delta^{\nabla^{g_{\boldsymbol\theta}}}\mathbf w|
  +|\partial(d^{\nabla^{g_{\boldsymbol\theta}}}\mathbf w)|
  +|\partial(\delta^{\nabla^{g_{\boldsymbol\theta}}}\mathbf w)|\notag\\
 &\qquad\leq C\bigl(|\mathbf w|+|\partial\mathbf w|
 +|\mathbf A|+|\partial\mathbf A| +|\mathbf z|+|\partial\mathbf z|\bigr).\label{ecriticalfirstorderweyl}
\end{align}
Here $\partial$ acts on the scalar frame components.
The coefficient functions and their first derivatives are bounded
on the chosen neighborhood.

View $\mathbf w$ as a $\Lambda^2T^*M$-valued two-form.
The Weitzenb\"ock identity
(see \cite[Theorem~(1.30)]{EellsLemaire1983}) reads
\[
 \bigl(
 d^{\nabla^{g_{\boldsymbol\theta}}}
       \delta^{\nabla^{g_{\boldsymbol\theta}}}
 +\delta^{\nabla^{g_{\boldsymbol\theta}}}
       d^{\nabla^{g_{\boldsymbol\theta}}}
 \bigr)\mathbf w
 =-\operatorname{tr}_{g_{\boldsymbol\theta}}
       \bigl((\nabla^{g_{\boldsymbol\theta}})^2\mathbf w\bigr)
   +\mathcal R_{g_{\boldsymbol\theta}}\mathbf w.
\]
The trace is over the two differentiating slots.
The term $\mathcal R_{g_{\boldsymbol\theta}}\mathbf w$
consists of contractions of the Riemann curvature of
$g_{\boldsymbol\theta}$ with $\mathbf w$ and satisfies
$|\mathcal R_{g_{\boldsymbol\theta}}\mathbf w|\leq C|\mathbf w|$.
The connection Laplacian 
$-\operatorname{tr}_{g_{\boldsymbol\theta}}
((\nabla^{g_{\boldsymbol\theta}})^2\mathbf w)$
has principal part $-L\mathbf w$,
with remaining terms bounded by
$C(|\partial\mathbf w|+|\mathbf w|)$. 
Since $d^{\nabla^{g_{\boldsymbol\theta}}}$ and
$\delta^{\nabla^{g_{\boldsymbol\theta}}}$ are first-order
operators with bounded coefficients,
\eqref{ecriticalfirstorderweyl} also bounds
$\left|
 \bigl(d^{\nabla^{g_{\boldsymbol\theta}}}
       \delta^{\nabla^{g_{\boldsymbol\theta}}}
      +\delta^{\nabla^{g_{\boldsymbol\theta}}}
       d^{\nabla^{g_{\boldsymbol\theta}}}\bigr)\mathbf w
 \right|$ 
by the same right-hand side. 
Consequently, 
\[
 |L\mathbf w|
 \leq C\bigl(
 |\mathbf w|+|\partial\mathbf w|
 +|\mathbf A|+|\partial\mathbf A|
 +|\mathbf z|+|\partial\mathbf z|\bigr).
\]

Finally, the definition of $\mathbf v$ and
\eqref{ecriticalexacttransport} give
\[
 |\mathbf A|+|\partial\mathbf A|
 +|\mathbf z|+|\partial\mathbf z|
 \leq C(|\mathbf v|+|\mathbf w|).
\]
Differentiating \eqref{ecriticalexacttransport} tangentially, we obtain
\[
 |\partial_t\partial_y\mathbf A|
 +|\partial_t\partial_y\mathbf z|
 \leq C(|\mathbf v|+|\mathbf w|+|\partial_y\mathbf w|).
\]
Together with \eqref{ecriticalexacttransport}, this proves \eqref{ecriticalclosedsystem}.
\end{proof}

We now verify the boundary conditions in
Lemma~\ref{lemmacoupledcurvatureuc} and apply it to obtain
$\mathbf B=\overline{\mathbf B}$.

\begin{proof}[Proof of Proposition~\ref{propcriticalradialmodel}]
Let $\mathbf{w},\mathbf{z},\mathbf v$ be the arrays in
Lemma~\ref{lemmacriticaldifferencesystem}.
By Lemma~\ref{lemmacriticalreferencejets},
$\mathbf A=O(t^\infty)$.
The $(2,3)$ block then shows that
$\boldsymbol{\theta}-\overline{\boldsymbol{\theta}}=O(t^\infty)$.
Thus $g_{\boldsymbol{\theta}},\overline g$ and their respective frames have
matching boundary jets, so $\mathbf{w}=O(t^\infty)$.

Since $\mathbf A=O(t^\infty)$, by \eqref{diffOmegaAB}, we have
 $\boldsymbol\Omega-\boldsymbol\Omega_0=O(t^\infty)$.
By \eqref{ecriticalphysicalcotton}, the $(2,1)$ block of this
difference is $-\mathbf W_{g_{\boldsymbol{\theta}}}\mathbf{b}$, and hence
$\mathbf W_{g_{\boldsymbol{\theta}}}\mathbf{b}=O(t^\infty)$.
Equation~\eqref{ecriticalweylarraycomparison} now yields
$\mathbf{z}=\mathbf W_{\overline g}\mathbf{b}=O(t^\infty)$.
Consequently,
\[
 \mathbf{w}|_{t=0}=\partial_t\mathbf{w}|_{t=0}=\mathbf v|_{t=0}=0.
\]
At $t=0$, $g_{\boldsymbol{\theta}}=g$ and $V=\partial_r$, so $\partial_t$
is the inward $g_{\boldsymbol{\theta}}$-unit normal.
Apply Lemma~\ref{lemmacoupledcurvatureuc} with metric $g_{\boldsymbol{\theta}}$,
$r=t$, and $\mathbf u=(w_{ik\ell m})_{i,k,\ell,m}$ to the
inequalities in \eqref{ecriticalclosedsystem}.
Hence $\mathbf{w}=\mathbf v=0$ on a smaller neighborhood.
In particular $\mathbf A=0$, so $\mathbf B=\overline{\mathbf B}$,
which proves \eqref{ecriticalradialmodel}.
\end{proof}

\subsection{Local conformal rigidity}
 We now combine boundary determination with
Proposition~\ref{lemmacurvatureboundaryuc} and
Proposition~\ref{propcriticalradialmodel} to construct
an exact local conformal map.

\begin{proposition}\label{propcriticalconformalcollar}

Under the hypotheses of Theorem~\ref{thmharmoniccurvaturein},
suppose that $n=2p+2$. Then there exist a nonempty open patch $\Gamma_0\Subset\Gamma_1$ and a smooth diffeomorphism $F_0$ between neighborhoods of $\Gamma_0,\kappa(\Gamma_0)$, with
\[
 F_0^*g_2=e^{2f_0}g_1,\qquad
 F_0|_{\Gamma_0}=\kappa,\qquad f_0|_{\Gamma_0}=0.
\]
The boundary jets $j^\infty_{\Gamma_0}F_0$ and $j^\infty_{\Gamma_0}f_0$ can be
chosen to agree with those of the map and factor supplied by
Theorem~\ref{thboundaryintrodo}(ii).
\end{proposition}

\begin{proof}
  By Theorem~\ref{thboundaryintrodo}(ii), there are a smooth
collar diffeomorphism $F$ and a smooth function $f$ satisfying
$F|_{\Gamma_1}=\kappa$ and
\begin{equation}\label{ecriticalformalboundary}
 j_{\Gamma_1}^\infty(F^*g_2)=j_{\Gamma_1}^\infty(e^{2f}g_1),\qquad f|_{\Gamma_1}=0.
\end{equation}
Since $C_{g_1}=C_{g_2}=0$, naturality and
\eqref{cottontransformfor} give
\begin{equation}\label{ecriticalformalcotton}
 W_{g_1}(\nabla^{g_1}f)=O(r^\infty).
\end{equation}
By
\eqref{ecriticalformalboundary}, the scalar curvature
transformation formula gives
\[
 2\Delta_{g_1}f-(n-2)|df|_{g_1}^2=\frac{\mathrm{Scal}_{g_2}e^{2f}-\mathrm{Scal}_{g_1}}{n-1}+O(r^\infty).
\]
In boundary normal coordinates $g_1=dr^2+h_r$, we have
\[
 \Delta_{g_1}f =-\partial_r^2f+H_r\partial_rf+\Delta_{h_r}f.
\]
Consequently,
\begin{equation}\label{ecriticalformalscalar}
 \partial_r^2f=H_r\partial_rf+\Delta_{h_r}f
 -\frac{n-2}{2}\bigl((\partial_rf)^2+|d_{Y_r}f|_{h_r}^2\bigr)
 +\frac{\mathrm{Scal}_{g_1}-\mathrm{Scal}_{g_2}e^{2f}}{2(n-1)}+O(r^\infty).
\end{equation}
The scalar curvatures $\mathrm{Scal}_{g_j}$ are constant since both metrics have harmonic curvature.

\smallskip\noindent
\emph{Step 1. Reduction to two local constructions.}
If
\[
 \partial_rf|_{\Gamma_1}=0,\qquad \mathrm{Scal}_{g_1}=\mathrm{Scal}_{g_2},
\]
successive differentiation of \eqref{ecriticalformalscalar} gives
\[
 f=O(r^\infty),\qquad F^*g_2-g_1=O(r^\infty) \quad\text{on }\Gamma_1.
\]
Choose any sufficiently small boundary coordinate patch
$\Gamma_0\Subset\Gamma_1$ and proceed to Step 2.

In the remaining cases, choose $\Gamma_0$ as follows.
\begin{enumerate}[label=(\roman*)]
\item If $\partial_rf|_{\Gamma_1}\not\equiv0$, choose
$\Gamma_0$ so that $\partial_rf|_{r=0}$ is nowhere zero
on $\overline{\Gamma_0}$.
After shrinking the collar, $\partial_rf$ remains nonzero.

\item If $\partial_rf|_{\Gamma_1}=0$ and
$\mathrm{Scal}_{g_1}\neq\mathrm{Scal}_{g_2}$,
choose any sufficiently small $\Gamma_0$.
Since $f|_{r=0}=\partial_rf|_{r=0}=0$,
\eqref{ecriticalformalscalar} and Taylor's formula give
\[
 f=\frac{\mathrm{Scal}_{g_1}-\mathrm{Scal}_{g_2}} {4(n-1)}r^2+O(r^3).
\]
Consequently, both $\nabla^{g_1}f$ and $\partial_rf$
are divisible by $r$, and $(\partial_rf)/r$ is nonzero
at the boundary.
\end{enumerate}

In either subcase, the tangential component of
$\nabla^{g_1}f$ vanishes to one order higher than $\partial_rf$.
Hence the quotient extends smoothly to the boundary and satisfies
\[
 \left.\frac{\nabla^{g_1}f}{\partial_rf}\right|_{r=0}
 =\partial_r.
\]
Moreover, by \eqref{ecriticalformalcotton},
\[
 W_{g_1}\!\left(\frac{\nabla^{g_1}f}{\partial_rf}\right)=O(r^\infty).
\]
In subcase (ii), this uses the fact that division by $r$
preserves $O(r^\infty)$.
Lemma~\ref{lemmacriticalnormalweyl} therefore yields
$W_{g_1}(\partial_r)=O(r^\infty)$ on $\Gamma_0$.

Apply the same argument to the inverse formal relation,
with map $F^{-1}$ and factor $-f\circ F^{-1}$.
At the boundary,
\[
 dF(\partial_{r_1})=\partial_{r_2},\qquad
 \bigl[\partial_{r_2}(-f\circ F^{-1})\bigr]\circ\kappa
 =-\partial_{r_1}f.
\]
Thus, after interchanging $g_1$ and $g_2$, the inverse
relation satisfies the same subcase. We conclude that
\begin{equation}\label{ecriticalbothnormalweyl}
 W_{g_j}(\partial_{r_j})=O(r_j^\infty),\qquad j=1,2,
\end{equation}
on $\Gamma_0$ and $\kappa(\Gamma_0)$, respectively.
Proceed to Step 3.

\smallskip\noindent
\emph{Step 2. Construction from matching metric jets.}
In this step we already have $f=O(r^\infty)$ and $F^*g_2-g_1=O(r^\infty)$.
Let $\Psi$ send each $g_1$-normal geodesic to the
$F^*g_2$-normal geodesic from the same boundary point, preserving
the normal parameter. The geodesic equations give
\[
 j^\infty_{\Gamma_0}\Psi=j^\infty_{\Gamma_0}\Id,\qquad
 j^\infty_{\Gamma_0}(\Psi^*F^*g_2)=j^\infty_{\Gamma_0}g_1.
\]
Both metrics are now in common boundary normal coordinates.
Proposition~\ref{lemmacurvatureboundaryuc} gives
\[
 \Psi^*F^*g_2=g_1.
\]
Set $F_0:=F\circ\Psi$ and $f_0:=0$. Then
\[
 F_0^*g_2=g_1,\qquad
 j^\infty_{\Gamma_0}F_0=j^\infty_{\Gamma_0}F,\qquad
 j^\infty_{\Gamma_0}f_0=j^\infty_{\Gamma_0}f.
\]

\smallskip\noindent
\emph{Step 3. Construction under the Weyl condition.}
In this step we already have \eqref{ecriticalbothnormalweyl}.
Identify the initial orthonormal frames by $dF$, and apply
Proposition~\ref{lemmacriticalcoordinates} to $g_j$ and $\eta_{j,0}$ with
\[
 \eta_{1,0}=0,\qquad F^*\eta_{2,0}=-df|_{\Gamma_0}.
\]
Both initial one-forms annihilate boundary tangent vectors.
Use the same boundary coordinate labels via $\kappa$.
The two ODE coordinate systems then define a map
$F_0(t,y)=(t,y)$ between the two collars.
Write $\boldsymbol{\theta}_j,\varphi_j,\boldsymbol{\beta}_j,\mathbf B_j$ for the objects
constructed from $g_j$, and
$g_{\boldsymbol{\theta}_j}=e^{-2\varphi_j}g_j$.

For an exact conformal change $g_1\mapsto e^{2f}g_1$,
Proposition~\ref{lemmacriticalcoordinates} keeps the coordinates and
frame fixed and shows that
 \begin{align*}
 e^{-2(\varphi_1+f)}e^{2f}g_1&=g_{\boldsymbol{\theta}_1},\\
 -(\eta_1-df)(e_{1,i})-e_{1,i}(\varphi_1+f)&=\beta_{1,i}.
 \end{align*}
By the conformal invariance above, it suffices to compare
the constructions for $e^{2f}g_1$ and $F^*g_2$.
These metrics have matching boundary jets by
\eqref{ecriticalformalboundary}. Successive differentiation of the defining ODEs determines
all boundary jets of their solutions from the metric jets
and the prescribed initial data. Since these agree, we obtain
 \begin{align*}
 F_0^*\boldsymbol{\theta}_2-\boldsymbol{\theta}_1&=O(t^\infty),\\
 F_0^*g_{\boldsymbol{\theta}_2}-g_{\boldsymbol{\theta}_1}&=O(t^\infty),\\
 \boldsymbol{\beta}_2\circ F_0-\boldsymbol{\beta}_1&=O(t^\infty).
 \end{align*}
By \eqref{ecriticalphysicalb}, the connection matrices therefore
have the same tangential boundary value:
\[
 (F_0^*\mathbf B_2)|_{T\Gamma_0} =\mathbf B_1|_{T\Gamma_0}=:\mathbf B_0.
\]
Proposition~\ref{propcriticalradialmodel} and
\eqref{ecriticalbothnormalweyl} now give
\[
 F_0^*\mathbf B_2 =\mathbf X_0\,dt+e^{-t\mathbf X_0}\mathbf B_0e^{t\mathbf X_0}
 =\mathbf B_1.
\]
Taking the $(2,3)$ block yields
$F_0^*\boldsymbol{\theta}_2=\boldsymbol{\theta}_1$, hence
\[
 F_0^*g_2=e^{2f_0}g_1,\qquad
 f_0:=\varphi_2\circ F_0-\varphi_1,\qquad f_0|_{\Gamma_0}=0.
\]
The ODE jet comparison also shows that
\[
 j^\infty_{\Gamma_0}F_0=j^\infty_{\Gamma_0}F,\qquad
 j^\infty_{\Gamma_0}f_0=j^\infty_{\Gamma_0}f.
\]
\end{proof}

\subsection{Global partial-data uniqueness}
\label{subsectioncriticalglobal}

For the four-dimensional isometric conclusion in Remark~\ref{remark4dimensionweyl}, it suffices
to find a point of $\Gamma_1$ where $W_{g_1}\neq0$.
The next lemma provides such a point whenever
$W_{g_1}\not\equiv0$.
\begin{lemma}\label{lemmaweylboundaryuc}
Let $g$ be a smooth metric with harmonic curvature on a connected
$n$-manifold with smooth boundary, $n\geq4$.
If $W_g$ vanishes on a nonempty open boundary patch, then
$W_g\equiv0$.
\end{lemma}

\begin{proof}
Regard $W_g$ as a two-form with values in $\Lambda^2T^*M$.
Extend $g$ and its induced connection smoothly across the patch
as in Section~\ref{section::defiitionsext}, and let
$\widetilde W_g$ be the zero extension of $W_g$.
Since the full boundary value vanishes, no boundary distribution
is produced in either identity in \eqref{eharmonicweylidentities}.
Thus
\[
 (d^\nabla+\delta^\nabla)\widetilde W_g=0,\qquad
 (d^\nabla+\delta^\nabla)^2\widetilde W_g=0
 \quad\text{distributionally}.
\]
The operator $(d^\nabla+\delta^\nabla)^2$ has principal symbol $|\xi|_g^2\Id$, so elliptic
regularity makes $\widetilde W_g$ smooth.
Aronszajn's unique continuation theorem
\cite[Remark~3, p.~248]{Aronszajn1957}, in the vector-bundle formulation stated by B\"ar
\cite[Section~4]{Bar1997}, then forces vanishing near the patch and throughout the
connected interior.
\end{proof}

We now prove the main theorem.
\begin{proof}[Proof of Theorem~\ref{thmharmoniccurvaturein}]
By DeTurck--Goldschmidt \cite[Theorem~1]{DeTurckGoldschmidt1989}, the metrics are
real-analytic in compatible interior harmonic coordinates.
It suffices to construct the map near the observed boundary
required by Proposition~\ref{rsmoothboundauation}.

If $n\neq2p+2$, Theorem~\ref{thboundaryintrodo}(i) gives matching
metric jets in corresponding boundary normal coordinates.
Proposition~\ref{lemmacurvatureboundaryuc} supplies a local isometry
near $\Gamma_1$, and Proposition~\ref{rsmoothboundauation}(i)
extends it globally.

If $n=2p+2$, Proposition~\ref{propcriticalconformalcollar} constructs a
conformal map near the nonempty open subset $\Gamma_0$ of the original observed boundary set $\Gamma_1$.
Proposition~\ref{rsmoothboundauation}(ii), applied to $\Gamma_0$ and the measurement data \eqref{equDNintropar} restricted to $\Gamma_0$, gives a smooth
global diffeomorphism $F$ and a smooth function $f$ such that
\[
 F^*g_2=e^{2f}g_1,\qquad f|_{\Gamma_0}=0.
\]

Next, we prove that $F$ is an extension of $\kappa$. 
Let $h_j$ denote the induced boundary metrics, and set
$\Phi:=F|_{\Gamma_1}$ and $\rho:=f|_{\Gamma_1}$.
Theorem~\ref{thboundaryintrodo}(ii) and the conformal identity give
\[
 \kappa^*h_2=h_1,\qquad \Phi^*h_2=e^{2\rho}h_1.
\]
Since $F$ extends the constructed collar map,
$\Phi=\kappa$ and $\rho=0$ on $\Gamma_0$.

The boundary dimension is $n-1\geq3$.
The connection and Schouten transformation formulas
\eqref{econnectiontrans}--\eqref{eschoutentrans},
applied to $h_j$, lead to a first-order ODE system along each boundary path
for $(\Phi,d\Phi,\rho,d\rho)$, with smooth coefficients.
The map $\kappa$, with conformal factor zero, satisfies
the same system. On the open set $\Gamma_0$,
\[
 (\Phi,d\Phi,\rho,d\rho)=(\kappa,d\kappa,0,0),
\]
where all differentials are taken on the boundary.
ODE uniqueness along paths and connectedness of $\Gamma_1$
therefore imply
\[
 F|_{\Gamma_1}=\kappa,\qquad f|_{\Gamma_1}=0.
\]

Suppose that at some $x_0\in\Gamma_1$,
\begin{equation*}
 W_{g_1}(\nu_{g_1},\cdot,\cdot,\cdot)|_{x_0}\neq0.
\end{equation*}
By continuity, choose an open boundary patch
$\Gamma\subset\Gamma_1$ containing $x_0$ on which
$W_{g_1}(\nu_{g_1})$ is nowhere zero.

Since both metrics have harmonic curvature,
$C_{g_1}=C_{g_2}=0$.
Using $F^*g_2=e^{2f}g_1$, it follows from \eqref{cottontransformfor} that
\[
 0=F^*C_{g_2}=C_{F^*g_2}=C_{e^{2f}g_1}=C_{g_1}+W_{g_1}(\cdot,\cdot,\cdot,\nabla^{g_1}f).
\]
Thus the symmetries imply
$W_{g_1}(\nabla^{g_1}f)=0$ on $M_1$.

The boundary condition $f|_{\Gamma_1}=0$ implies that
all tangential derivatives of $f$ vanish on $\Gamma$.
Consequently,
\begin{equation*}
 \nabla^{g_1}f
 =(\partial_{\nu_{g_1}}f)\nu_{g_1},\qquad
 (\partial_{\nu_{g_1}}f)W_{g_1}(\nu_{g_1})=0
 \quad\text{on }\Gamma.
\end{equation*}
Since $W_{g_1}(\nu_{g_1})$ is nowhere zero there,
$\partial_{\nu_{g_1}}f=0$, and hence $df|_\Gamma=0$.
Tangential differentiation of $df|_\Gamma=0$ and symmetry
of the Hessian show that its tangential and mixed components
vanish. Therefore,
\[
 \nabla^{g_1}_{\nu_{g_1}}\nabla^{g_1}f=(\nabla^{g_1}df)(\nu_{g_1},\nu_{g_1})\nu_{g_1}
 \quad\text{on }\Gamma.
\]
Now differentiate the identity
$W_{g_1}(\nabla^{g_1}f)=0$ along $\nu_{g_1}$.
Restricting to $\Gamma$, where $\nabla^{g_1}f=0$, gives
 \begin{align*}
 0&=(\nabla^{g_1}_{\nu_{g_1}}W_{g_1})(\nabla^{g_1}f) +W_{g_1}(\nabla^{g_1}_{\nu_{g_1}}\nabla^{g_1}f)\\
 &=(\nabla^{g_1}df)(\nu_{g_1},\nu_{g_1}) W_{g_1}(\nu_{g_1}).
 \end{align*}
The nonvanishing of $W_{g_1}(\nu_{g_1})$ forces the normal
Hessian component to vanish as well. Thus
\[
 f=0,\qquad df=0,\qquad \nabla^{g_1}df=0 \quad\text{on }\Gamma.
\]

Since both scalar curvatures are constant, the transformation
formula for the scalar curvature, with our nonnegative Laplacian, gives
\[
 e^{2f}\mathrm{Scal}_{g_2} =\mathrm{Scal}_{g_1}
  +2(n-1)\Delta_{g_1}f -(n-1)(n-2)|df|_{g_1}^2.
\]
On $\Gamma$, we have $f=0$, $df=0$, and
$\Delta_{g_1}f=-\operatorname{tr}_{g_1}(\nabla^{g_1}df)=0$.
Thus $\mathrm{Scal}_{g_1}=\mathrm{Scal}_{g_2}$.
The transformation formula consequently becomes
\[
 2\Delta_{g_1}f-(n-2)|df|_{g_1}^2
 =\frac{\mathrm{Scal}_{g_1}}{n-1}(e^{2f}-1).
\]
Since $f$ and $|df|_{g_1}$ are bounded, there is a constant $C_f>0$ such that
\begin{equation}\label{ellipticineq}
 |\Delta_{g_1}f|\leq C_f(|f|+|df|_{g_1}).
\end{equation}
Extend $g_1$ smoothly across $\Gamma$ and extend $f$ by zero.
The vanishing of $f$, $df$, and $\nabla^{g_1}df$ on $\Gamma$
makes the zero extension $C^2$, and it satisfies the same
elliptic inequality \eqref{ellipticineq} across $\Gamma$.
Unique continuation \cite[pp.~235--236]{Aronszajn1957} gives $f=0$ near
$\Gamma$ and then throughout the connected interior.
Hence $f\equiv0$ on $M_1$, and $F^*g_2=g_1$.
\end{proof}

\section*{Acknowledgments}

This work was partially supported by the National Natural Science Foundation of China (Grant No. 12571064). Y. J. was supported by the CAS Project for Young Scientists in Basic Research (Grant No. YSBR-032).
\section*{Declaration of AI Use}
During the preparation of this work, OpenAI's ChatGPT was used for language editing. It also helped develop some of the proof ideas for the main theorems and assisted with technical details, particularly in the proof of Proposition~\ref{lemmacurvatureboundaryuc} and Proposition~\ref{propcriticalradialmodel}.  All AI-assisted arguments and revisions were independently checked by the authors, who take full responsibility for the final manuscript.

\end{document}